%% file: main.tex
\documentclass[10pt]{article}
\usepackage{amsfonts,color}
\usepackage{amssymb}
\usepackage{footnote}
\usepackage{mathtools}
\usepackage{amsthm,wasysym}
\usepackage{amsmath}
\usepackage[english]{babel}
\usepackage{bbm}
\usepackage{colonequals}
\usepackage{mathrsfs}

\usepackage{graphicx}
\usepackage{amsfonts}
\usepackage{hyperref}
\usepackage{subcaption}
\usepackage{nicefrac}
\usepackage{gensymb}
\usepackage{enumerate}
\usepackage[nameinlink,capitalize]{cleveref}
\usepackage{cite}
\usepackage[nottoc,numbib]{tocbibind}
\usepackage{pgf,tikz,pgfplots}
\pgfplotsset{compat=1.14}
\usepackage{mathrsfs}
\usetikzlibrary{arrows}
\usetikzlibrary{arrows.meta}
\usepackage{overpic}

\usepackage{chngcntr}

\counterwithin*{equation}{section}

\newcommand{\vb}[1]{\mathbf{#1}}
\newcommand{\bm}[1]{\boldsymbol{#1}}

\DeclareMathOperator{\orth}{\mathrm{orth}}
\DeclareMathOperator{\sgn}{\mathrm{sgn}}

\DeclareMathOperator{\sym}{\mathrm{sym}}

\DeclareMathOperator{\sph}{\mathrm{sph}}
\DeclareMathOperator{\dev}{\mathrm{dev}}
\DeclareMathOperator{\skw}{\mathrm{skw}}

\DeclareMathOperator{\tr}{\mathrm{tr}}

\DeclareMathOperator{\cof}{\mathrm{cof}}
\newcommand{\jump}[1]{\ensuremath{[\![#1]\!]} }
\newcommand{\one}{\bm{\mathbbm{1}}}
\newcommand{\con}[2]{\langle {#1} , \, {#2} \rangle}
\newcommand{\norm}[1]{\| {#1} \|}
\newcommand\subsetsim{\mathrel{%
  \ooalign{\raise0.2ex\hbox{$\subset$}\cr\hidewidth\raise-0.8ex\hbox{\scalebox{0.9}{$\sim$}}\hidewidth\cr}}}

\newcommand{\dd}{\mathrm{d}}
\DeclareMathOperator{\di}{\mathrm{div}}
\DeclareMathOperator{\Di}{\mathrm{Div}}
\DeclareMathOperator{\rot}{\mathrm{rot}}
\DeclareMathOperator{\Rot}{\mathrm{Rot}}
\DeclareMathOperator{\curl}{\mathrm{curl}}
\DeclareMathOperator{\Curl}{\mathrm{Curl}}
\DeclareMathOperator{\Inc}{\mathrm{Inc}}

\DeclareMathOperator{\airy}{\mathrm{airy}}
\DeclareMathOperator{\D}{\mathrm{D}}

\newcommand{\Reg}{\mathcal{R}}
\newcommand{\CG}{\mathcal{CG}}
\newcommand{\DG}{\mathcal{DG}}

\newcommand{\so}{\mathfrak{so}}
\renewcommand{\sl}{\mathfrak{sl}}
\newcommand{\Sym}{\mathrm{Sym}}
\newcommand{\Ned}{\mathcal{N}_{I}}

\newcommand{\PS}{\mathcal{PS}}
\newcommand{\GLS}{\mathcal{GLS}}
\newcommand{\HHJ}{\mathcal{HHJ}}
\newcommand{\HZ}{\mathcal{HZ}}
\newcommand{\HMS}{\mathcal{HMS}}
\newcommand{\RT}{\mathcal{RT}}

\newcommand{\Nedtwo}{\mathcal{N}_{II}}

\newcommand{\Po}{\mathit{P}}
\newcommand{\Le}{{\mathit{L}^2}}

\newcommand{\Hone}{\mathit{H}^1}
\newcommand{\Honez}{\mathit{H}_0^1}
\newcommand{\HsD}[1]{\mathit{H}^\mathrm{sym}(\mathrm{Div}{#1})}

\newcommand{\Hd}[1]{\mathit{H}(\mathrm{div}{#1})}

\newcommand{\HD}[1]{\mathit{H}(\mathrm{Div}{#1})}

\newcommand{\Hc}[1]{\mathit{H}(\mathrm{curl}{#1})}
\newcommand{\HC}[1]{\mathit{H}(\mathrm{Curl}{#1})}

\newcommand{\HdD}[1]{\mathit{H}(\mathrm{div}\,\mathrm{Div}{#1})}

\newcommand{\HR}[1]{\mathit{H}(\mathrm{Rot}{#1})}

\newcommand{\HcD}[1]{\mathit{H}(\mathrm{curl}\,\mathrm{Div}{#1})}
\newcommand{\HrD}[1]{\mathit{H}(\mathrm{rot}\,\mathrm{Div}{#1})}
\newcommand{\HcC}[1]{\mathit{H}(\mathrm{Inc}{#1})}
\newcommand{\HrR}[1]{\mathit{H}(\mathrm{rot}\,\mathrm{Rot}{#1})}

\newcommand{\body}{\mathfrak{B}}
\newcommand{\vol}{V}
\newcommand{\surf}{A}
\newcommand{\curv}{s}

\newcommand{\R}{\mathbb{R}}

\newcommand{\V}{\mathit{V}}
\newcommand{\Z}{\mathit{Z}}

\newcommand{\C}{\mathit{C}}
\renewcommand{\H}{\mathit{H}}
\newcommand{\tem}{\mathcal{T}}
\newcommand{\ver}{\mathcal{V}}
\newcommand{\edge}{\mathcal{E}}
\newcommand{\face}{\mathcal{F}}
\newcommand{\cell}{\mathcal{C}}

\newcommand{\tv}{\bm{\psi}}
\renewcommand{\tt}{\bm{\Psi}}

\newcommand{\ksmu}{k_s \mu}

\newcommand{\J}{\mathbb{J}}
\newcommand{\A}{\mathbb{A}}

\newcommand{\RN}[1]{%
	\textup{\uppercase\expandafter{\romannumeral#1}}%
}

\allowdisplaybreaks[1]

\makeindex
                    
\usepackage{chngcntr}              

\newtheoremstyle{break}
{\topsep}{\topsep}%
{\itshape}{}%
{\bfseries}{}%
{\newline}{}%
\theoremstyle{break}

\newtheorem{remark}{Remark}

\newtheorem{definition}{Definition}

\counterwithin{theorem}{section}
\counterwithin{lemma}{section}
\counterwithin{remark}{section}
\counterwithin{observation}{section}
\counterwithin{definition}{section}

\renewcommand{\thefootnote}{\arabic{footnote}}

   \makeatletter
\let\@fnsymbol\@arabic
\makeatother

\crefname{Problem}{Problem.}{Problem.}

\title{Formulae and transformations for simplicial tensorial finite elements via polytopal templates}
\author{\normalsize{Adam Sky}\thanks{Corresponding author: Adam Sky, Institute of Computational Engineering and Sciences, Department of Engineering, Faculty of Science, Technology and Medicine, University of Luxembourg, 6 Avenue de la Fonte, L-4364, Esch-sur-Alzette, Luxembourg, email: adam.sky@uni.lu}
    , \quad
    \normalsize{Michael Neunteufel}\thanks{Michael Neunteufel, Department of Mathematics and Statistics, Portland State University, Portland, OR, United States, email: mneunteu@pdx.edu} 
    , \quad 
	\normalsize{Jack S. Hale}\thanks{Jack S. Hale, Institute of Computational Engineering and Sciences, Department of Engineering, Faculty of Science, Technology and Medicine, University of Luxembourg, 6 Avenue de la Fonte, L-4362 Esch-sur-Alzette, Luxembourg, email: jack.hale@uni.lu}
    , \quad
    \normalsize{and} \quad
	\normalsize{Andreas Zilian}\thanks{Andreas Zilian, Institute of Computational Engineering and Sciences, Department of Engineering, Faculty of Science, Technology and Medicine, University of Luxembourg, 6 Avenue de la Fonte, L-4364, Esch-sur-Alzette, Luxembourg, email: andreas.zilian@uni.lu}
}

\begin{document}

\maketitle

\begin{abstract}
We introduce a unified method for constructing the
basis functions of a wide variety of partially continuous tensor-valued finite elements on simplices using polytopal templates. These finite element spaces are essential
for achieving well-posed discretisations of mixed formulations of partial differential equations that involve tensor-valued functions,
such as the Hellinger--Reissner formulation of linear elasticity. In
our proposed polytopal template method, the basis functions are
constructed from template tensors associated with the geometric
polytopes (vertices, edges, faces etc.) of the reference simplex and
any scalar-valued $\Hone$-conforming finite element space. From this
starting point we can construct the Regge,
Hellan--Herrmann--Johnson, Pechstein--Schöberl, Hu--Zhang, Hu--Ma--Sun and Gopalakrishnan--Lederer--Sch\"oberl elements. Because the Hu--Zhang element and the Hu--Ma--Sun
element cannot be mapped from the reference simplex
to a physical simplex via standard double Piola mappings, we also
demonstrate that the polytopal template tensors can be used to define a
consistent mapping from a reference simplex even to a non-affine simplex in the physical mesh. Finally, we discuss the implications of element regularity with two numerical examples for the Reissner--Mindlin plate problem.   


\end{abstract}

\section{Introduction}
Tensor-valued functions $\bm{P}: \body \subset \R^d \to \R^{d \times d}$ for $d =
\{ 2, 3 \}$ naturally appear in many partial differential equations. Correspondingly, the variational formulations of these partial differential equations necessitate the construction of conforming tensor-valued finite element subspaces. The construction of these spaces as well as the conditions needed for the existence of unique solutions are studied in the field of finite element exterior calculus (FEEC) \cite{arnold_finite_2006,arnold_complexes_2021}. In fact, FEEC suggests that variational formulations are to be investigated in the context of associated Hilbert space complexes. Indeed, the treatment of variational problems using Hilbert space complexes and FEEC allows for the clear discernment of the properties of the finite element
subspaces for stability. Specifically, the existence of a Hilbert complex in which the variational problem can framed, and the definition of commuting interpolants with respect to said complex, are usually sufficient to immediately assert existence and uniqueness of solutions on conforming finite element subspaces via Fortin's criterion~\cite{Bra2013}. Because of this elegant
result, much effort has been invested in the construction of finite element subspaces that define commuting interpolants, see e.g.~\cite{Demkowicz2000,arnold_mixed_2002,hu_finite_2016}.

Classical examples for variational problems that fit into a Hilbert space complex are the mixed Hellinger--Reissner and Hu--Washizu
formulations of linear elasticity
\cite{SKY2024112808,NEUNTEUFEL2021113857,sky2023reissnermindlin}. Specifically, the variational forms of the latter can be studied and understood in the a sense of the so called elasticity complex \cite{pauly_elasticity_2022,pauly_hilbert_2022}, also known as the Kr\"oner
complex \cite{Kroner1954}. The elasticity complex insinuates the definition of the tensor-valued
stress field $\bm{\sigma}: \body \subset \R^d \to \R^{d \times d}$ in $\HsD{,\body}$, being the Hilbert space of square-integrable symmetric
tensors with row-wise square-integrable divergence. The first conforming $\HsD{,\body}$-elements on simplices satisfying the commuting property are the Arnold--Winther
element~\cite{arnold_mixed_2002} on triangles and the Arnold--Awanou--Winther
element~\cite{arnold_finite_2008} on tetrahedra. Both elements are based on incomplete polynomial
spaces and were introduced without a mapping procedure from the reference simplex
to a general simplex in the physical mesh. Without a mapping, the basis functions
must be computed for each individual simplex in the physical mesh. A recent
generalised approach to transformations \cite{kirby_general_2018} derived a
transformation procedure for affine triangles that can be applied to the
Arnold--Winther element \cite{aznaran_transformations_2021}. However, the
procedure does not extend to the case of non-affinely mapped (curved) triangles.
A recent alternative to the Arnold--Winther elements are the Hu--Zhang 
\cite{hu_family_2014,hu_family_2015,hu_finite_2016,Hu2014Q,Hu2015Q} and Hu--Ma--Sun elements \cite{hu2023new}. 
Unlike, the Arnold--Winther and the Arnold--Awanou--Winther elements, these elements are based on complete polynomial spaces. Consequently, we were able to introduce a non-affine mapping from the reference triangle to the physical triangle on the mesh for the Hu--Zhang element in \cite{sky2023reissnermindlin}, using the polytopal template approach \cite{sky_polytopal_2022}.
This was achieved by exploiting the association of the basis functions with the polytopal structure of the triangle. 

It is important to note that none of the aforementioned $\HsD{,\body}$-conforming elements are minimally $\HD{,\body}$-regular. In fact, as proved in \cite{arnold_mixed_2002}, it is impossible to construct
minimally regular, yet conforming \textit{\textbf{standard}} simplex finite elements for $\HsD{,\body}$.
Consequently, the Arnold--Winther, Hu--Zhang, Arnold--Awanou--Winther and
Hu--Ma--Sun elements, all increase the
regularity of the construction on vertices in two dimensions, and at vertices
and edges in three-dimensions. The only known exception to this problem are the Johnson--Mercier \textit{\textbf{macro}} elements \cite{gopalakrishnan2024johnsonmercier}, which are defined on Clough--Tocher splits \cite{Philippe} in two dimensions, or Křížek--Alfeld splits \cite{Křížek1982,ALFELD1984169} in three dimensions.  
The lack of the minimal regularity property for standard simplices has led a number of authors to introduce approaches that avoid working directly in $\HsD{,\body}$. A recent example is the
tangential-displacement-normal-normal-stress (TDNNS) method
\cite{pechstein_analysis_2018,pechstein_anisotropic_2012,pechstein_tdnns_2017},
which is based on the Hellan--Herrmann--Johnson principle \cite{arnold_hellan--herrmann--johnson_2020,Hel67,Her67,Joh73}. The
method reinterprets the inner product of the divergence of the stress tensor
and elastic displacement field $\con{\Di \bm{\sigma}}{\vb{u}}_{\Le}$ in the Hellinger-Reissner formulation of linear
elasticity as a dual product between the spaces $\H'(\mathrm{curl},\body)$ and
$\Hc{,\body}$, $\con{\Di \bm{\sigma}}{\vb{u}}_{\H'(\mathrm{curl})\times \Hc{}}$. As a result, the tensorial stress field $\bm{\sigma}$ is redefined in the space $\HdD{,\body}$, such
that its divergence $\Di \bm{\sigma}$ lies in $\H^{-1}(\mathrm{div},\body) = \H'(\mathrm{curl},\body)$. In
two dimensions, the Hellan--Herrmann--Johnson elements
\cite{neunteufel_hellanherrmannjohnson_2019,Hel67,Her67,Joh73} are minimally regular and `almost' $\HdD{,\surf}$-conforming
$\HHJ^p(\surf) \subsetsim \HdD{,\surf}$. In three-dimensions, the analogous
elements are the Pechstein--Sch\"oberl elements \cite{PS2011,Sin2009} on
tetrahedra $\PS^p(\vol) \subsetsim \HdD{,\vol}$. Befittingly, the method uses N\'ed\'elec elements \cite{SKY2022115298,sky_hybrid_2021,Nedelec1980,Ned2,haubold2023high} to
discretise the displacement field $\vb{u} \in \Nedtwo^p(\body) \subset \Ned^p(\body) \subset \Hc{,\body}$.

Analogously to mixed linear elasticity, other variational problems can be studied using Hilbert space complexes and tensor-valued finite elements. In the case of the biharmonic equation, the associated
complex is the $\di \Di$-complex \cite{PaulyDiv,CRMECA_2023__351_S1_A8_0,dipietro2023discrete},
which can be completed into the relaxed micromorphic complex
\cite{SkyOn,SKYNOVEL,Lewintan2021,LewintanInc,LewintanInc2,NEFF20151267,Neff2012,Gmeineder1,Gmeineder2}
for the relaxed micromorphic model
\cite{sky_higher_2023,SkyPamm,SKY2022115298,sky_hybrid_2021,GOURGIOTIS2024112700,Neff2014}. Naturally, the Hellan--Herrmann--Johnson \cite{neunteufel_hellanherrmannjohnson_2019,Hel67,Her67,Joh73} and the Pechstein--Sch\"oberl elements \cite{PS2011,Sin2009} are applicable to the $\di\Di$-complex for discretisations of fields in $\HdD{,\body}$.

Regge elements originally stem from Regge calculus for solving Einstein's field
equations in a coordinate-free manner \cite{Regge61}. In \cite{Chr2011} the
concept of tangential-tangential continuous tensors has been embedded in the
FEEC setting as a (slightly) non-conforming subspace of the function space
$\HcC{,\body}$, $\Reg^p(\body) \subsetsim \HcC{,\body}$, which also naturally appears in the elasticity complex.
Higher order polynomial Regge elements on
arbitrary simplex-dimensions were developed in \cite{Li2018}. Regge elements
have been successfully used to discretise strain $\Reg(\body)\ni \bm{\varepsilon}:\body \subset \R^{d} \to \R^{d \times d}$ and metric fields in continuum mechanics, shells, and curvature approximations \cite{HH2013,NS21,GNSW2023}.

The tangential-normal continuous Gopalakrishnan--Lederer--Sch\"oberl elements
for triangles and tetrahedra have been introduced in
\cite{Led2019,gopalakrishnan2020mass,GKLS2023}. The elements are intrinsically related to the grad-curl complex \cite{arnold_complexes_2021,Kaiboqcurl} and build a non-conforming
subspace of $\HcD{,\body}$, $\GLS^p(\body) \subsetsim \HcD{,\body}$. Their divergence lies in
$\H^{-1}(\mathrm{curl},\body) = \H'(\mathrm{div},\body)$ and can be paired with
$\Hd{,\body}$-conforming Raviart--Thomas or Brezzi--Douglas--Marini elements \cite{sky_polytopal_2022,BDM,Raviart}. Due to
their intrinsically trace-free construction they are well-suited for discretisations of the (Navier-)Stokes equations with exactly divergence-free
velocity fields~\cite{Devloo2}. 

In this work we introduce a unified method of constructing tensor-valued finite
elements on simplices using polytopal templates. The polytopal template method was introduced in
\cite{sky_polytopal_2022} for the construction of the vector-valued N\'edelec,
Raviart--Thomas and Brezzi--Douglas--Marini elements \cite{Nedelec1980,Ned2,BDM,Raviart}. In this paper we take the
vectorial template sets from \cite{sky_polytopal_2022} and use them to
construct tensorial template sets. From these tensorial template sets we can
construct the well-known Regge, Hellan--Herrmann--Johnson,
Pechstein--Sch\"oberl, and Gopalakrishnan--Lederer--Sch\"oberl finite
elements. With the addition of Cartesian and interface-unique orthogonal template
sets, we also show it is possible to construct the Hu--Zhang and Hu--Ma--Sun
elements.
By virtue of the polytopal method, the clear association with the polytopes of a simplex allows us to introduce
consistent transformations of the basis functions from the reference to the physical simplex for all constructions.
We emphasise that in the case of an affine transformation from the reference simplex to a physical simplex in the mesh, the polytopal template method is independent of the underlying scalar-valued finite element basis that is chosen in the construction of the
tensorial finite elements elements. The underlying polynomial subspace must
simply be $\Hone(\body)$-conforming. Thus, any such continuous Lagrangian subspace $\CG^p(\body) \subset \Hone(\body)$ may be chosen. Accordingly, the
resulting tensorial finite element directly inherits many of its properties
from the underlying scalar polynomial subspace. For example, higher order,
heterogeneous p-refinement \cite{DEVLOO20091716}, optimal complexity \cite{AinsworthOpt}, or $\Le(T)$-orthogonality \cite{Beuchler2006} are all
possible. \textit{\textbf{For the non-affine case, it is necessary to
choose a hierarchical basis; this point is discussed further in the paper}}.

This work is organised as follows. We shortly outline the notation employed in
this work. Afterwards in \cref{sec:recap}, we recap the polytopal template method and recall the
vectorial template sets. \cref{sec:tan} is dedicated to the presentation of
tangential-continuous finite elements, where we introduce the Regge element. In \cref{sec:norm} we discuss normal-continuous finite elements and
derive polytopal template variants of the Hellan--Herrmann--Johnson,
Pechstein--Sch\"oberl, Hu--Zhang, and Hu--Ma--Sun elements. In \cref{sec:tannorm} we present the construction of the Gopalakrishnan--Lederer--Sch\"oberl
elements, which are tangential-normal continuous. In \cref{sec:appRM} we
demonstrate an application of the Hu--Zhang and Hellan--Herrmann--Johnson elements for computations of the Reissner--Mindlin plate problem.
Finally, we give our conclusions and outlook.    

\section{Preliminaries}

\subsection{Notation}
The following notation is used throughout this work.
Exceptions to these rules are made clear in the precise context.
\begin{itemize}
    \item Vectors are defined as bold lower-case letters $\vb{v}, \, \bm{\xi}$.
    \item Second order tensors are denoted with bold capital letters $\bm{T}$.
    \item Fourth-order tensors are designated by the blackboard-bold format $\mathbb{A}$.
    \item We denote the Cartesian basis as $\{\vb{e}_1, \, \vb{e}_2, \, \vb{e}_3\}$.
    \item Summation over indices follows the standard rule of repeating indices. Latin indices represent summation over the full dimension, whereas Greek indices define summation over the co-dimension.
    \item The angle-brackets are used to define scalar products of arbitrary dimensions $\con{\vb{a}}{\vb{b}} = a_i b_i$, $\con{\bm{A}}{\bm{B}} = A_{ij}B_{ij}$.
    \item The matrix product is used to indicate all partial-contractions between a higher-order and a lower-order tensor $\bm{A}\vb{v} = A_{ij} v_j \vb{e}_i$, $\mathbb{A}\bm{B} = A_{ijkl}B_{kl}\vb{e}_i \otimes \vb{e}_j$.
    \item The second-order identity tensor is defined via $\one$, such that $\one \vb{v} = \vb{v}$. Analogously, the fourth-order identity tensor $\J$ yields $\J \bm{T} = \bm{T}$. 
    \item A general physical body of some arbitrary dimension $d$ is denoted with $\body \subset \R^d$.
    \item Volumes, surfaces and curves of the physical domain are identified via $\vol \subset \R^3$, $\surf \subset \R^2$ and $\curv \subset \R$, respectively. Their counterparts on the reference simplex are $\Omega \subset \R^3$, $\Gamma \subset \R^2$ and $\mu \subset \R$. Additionally, $T \subset \body \subset \R^{d}$ represents the domain of a single element in the physical mesh. 
    \item Tangential and normal vectors on the physical domain are designated by $\vb{t}$ and $\vb{n}$, respectively. On reference domain the respective vectors read $\bm{\tau}$ and $\bm{\nu}$. 
    \item The polytopes of an element are identified using multi-indices. For example, edge $e_j$ with $j = (0,1) \in \mathcal{J} = \{(0,1),(0,2),\dots\}$, where $\{0,1,2,\dots\}$ are the indices of the vertices $\{v_0, v_1, v_2\}$.
    \item The operator $\sim$ is to be understood as "associated with".
    \item In our construction of finite element spaces we use tensor products between sets. This notation is be understood for example as $\{1-\xi-\eta,\xi\} \otimes \{\vb{e}_1 + \vb{e}_2,\vb{e}_1\} = \{(1-\xi-\eta)(\vb{e}_1 + \vb{e}_2),(1-\xi-\eta)\vb{e}_1, \xi(\vb{e}_1 + \vb{e}_2), \xi \vb{e}_1\}$.
    \item We define the constant space of symmetric second order tensors as $\Sym(d) = \{ \bm{T} \in \R^{d \times d} \; | \; \bm{T} = \bm{T}^{T} \}$ with $ d \in \{2,3\}$.
    \item Its counterpart is the space of skew-symmetric tensors $\so(d) = \{ \bm{T} \in \R^{d \times d} \; | \; \bm{T} = -\bm{T}^{T} \}$.
   \item The spaces are associated with the operators $\sym \bm{T} = (1/2)(\bm{T} + \bm{T}^T) \in \Sym(d)$ and $\skw \bm{T} = (1/2)(\bm{T} - \bm{T}^T) \in \so(d)$, respectively.
   \item The space of constant deviatoric tensors reads $\sl(d) = \{ \bm{T} \in \R^{d \times d} \; | \; \tr \bm{T} = 0 \}$, where $\tr \bm{T} = \con{\bm{T}}{\one} = T_{ii}$.
   \item It is associated with the operator $\dev \bm{T} = \bm{T} - \sph \bm{T}$ where $\sph \bm{T} = (1/d) (\tr \bm{T}) \one$.
   \item The nabla operator is used to defined as $\nabla = \vb{e}_i \partial_i$.
   \item The left-gradient is given via $\nabla$, such that $\nabla \lambda = \nabla \otimes \lambda$.
    \item The right-gradient is defined for vectors and higher order tensors via $\D$, such that $\D \vb{v} = \vb{v} \otimes \nabla$.
    \item We define the vectorial divergence as $\di \vb{v} = \con{\nabla}{\vb{v}}$.
    \item The tensor divergence is given by $\Di \bm{T} = \bm{T} \nabla$, implying a single contraction acting row-wise 
    \item The vectorial curl operator reads $\curl \vb{v} = \nabla \times \vb{v}$
    \item For tensors the operator is given by $\Curl \bm{T} = -\bm{T} \times \nabla$, acting row-wise.
    \item The composite incompatibility operator reads $\Inc \bm{T} = \curl \Curl \bm{T} = - \nabla \times \bm{T} \times \nabla$.
    \item In two dimensions the vectorial $\curl$-operator induces the scalar operator $\nabla^\perp \lambda = \bm{R} \nabla \lambda$ and vectorial operator $\rot \vb{v} = \di(\bm{R}\vb{v})$, with $\bm{R} = \vb{e}_1 \otimes \vb{e}_2 - \vb{e}_2 \otimes \vb{e}_1$.
    \item Analogously, the tensorial $\Curl$-operator induces the vectorial operator $\D^\perp \vb{v} = (\D \vb{v}) \bm{R}^T$ and the tensorial operator $\Rot \bm{T} = \Di(\bm{T}\bm{R}^T)$, acting row-wise.
    \item Consequently, the composite incompatibility operator reduces to $\rot \Rot \bm{T}$ for tensors and the Airy-operator $\airy \lambda = \D^\perp \nabla^\perp \lambda$ for scalars. 
\end{itemize}

Further, we introduce the following Hilbert spaces and their respective norms 
    \begin{align}
\Le(\vol) &= \{ u : \vol \to \mathbb{R} \; | \; \| u \|_\Le < \infty  \} \, , & \|u\|_\Le^2 &= \int_\vol u^2 \, \dd \vol \, , \notag \\
    \Hone(\vol) &= \{ u \in \Le(\vol) \; | \; \nabla u \in [\Le(\vol)]^3 \} \, , & \norm{u}^2_{\Hone} &= \norm{u}^2_\Le + \norm{\nabla u}^2_\Le \, ,  \\
    \HC{,\vol} &= \{ \bm{T} \in [\Le(\vol)]^{3 \times 3} \; | \; \Curl \bm{T} \in [\Le(\vol)]^{3 \times 3} \} \, , & \norm{\bm{T}}^2_{\HC{}} &= \norm{\bm{T}}^2_\Le + \norm{\Curl \bm{T}}^2_\Le \, , \notag \\
    \HD{,\vol} &= \{ \bm{T} \in [\Le(\vol)]^{3 \times 3} \; | \; \Di \bm{T} \in [\Le(\vol)]^{3} \} \, , & \norm{\bm{T}}^2_{\HD{}} &= \norm{\bm{T}}^2_\Le + \norm{\Di \bm{T}}^2_\Le \, , 
\end{align}
where $\vol \subset \R^3$ implies a three-dimensional domain. 
In two dimensions $\surf \subset \R^2$ we define the analogous Hilbert spaces
\begin{align}
\Le(\surf) &= \{ u : \surf \to \mathbb{R} \; | \; \| u \|_\Le < \infty  \} \, , & \|u\|_\Le^2 &= \int_\surf u^2 \, \dd \surf \, , \notag \\
    \Hone(\surf) &= \{ u \in \Le(\surf) \; | \; \nabla u \in [\Le(\surf)]^2 \} \, , & \norm{u}^2_{\Hone} &= \norm{u}^2_\Le + \norm{\nabla u}^2_\Le \, ,  \\
    \HR{,\surf} &= \{ \bm{T} \in [\Le(\surf)]^{2 \times 2} \; | \; \Rot \bm{T} \in [\Le(\surf)]^{2} \} \, , & \norm{\bm{T}}^2_{\HR{}} &= \norm{\bm{T}}^2_\Le + \norm{\Rot \bm{T}}^2_\Le \, , \notag \\
    \HD{,\surf} &= \{ \bm{T} \in [\Le(\surf)]^{2 \times 2} \; | \; \Di \bm{T} \in [\Le(\surf)]^{2} \} \, , & \norm{\bm{T}}^2_{\HD{}} &= \norm{\bm{T}}^2_\Le + \norm{\Di \bm{T}}^2_\Le \, , \notag
\end{align}
where the differential operators are adjusted accordingly. By restricting a tensor field to be symmetric we define the Hilbert spaces  
\begin{align}
    \HsD{,\vol} &= \{ \bm{T} \in \HD{,\vol} \; | \; \bm{T} = \bm{T}^T \} \, , & 
    \HsD{,\surf} &= \{ \bm{T} \in \HD{,\surf} \; | \; \bm{T} = \bm{T}^T \} \, .  
\end{align}
Hilbert spaces with vanishing Sobolev traces \cite{Hiptmair} are marked with a zero-subscript, for example $\Honez(\body)$. Scalar products pertaining to the Hilbert spaces are indicated by a subscript on the angle-brackets
\begin{align}
    \con{u}{v}_{\Le} = \int_\vol \con{u}{v} \, \dd \vol \, ,
\end{align}
where the domain is clear from context.
Finally, we define the spaces 
\begin{align}
    \HcC{,\vol} &= \{ \bm{T} \in  \Le(\vol) \otimes \Sym(3) \; | \; \Inc \bm{T} \in \H^{-1}(\vol) \otimes \Sym(3)  \} \, , \\
    \HdD{,\vol} &= \{ \bm{T} \in  \Le(\vol) \otimes \Sym(3) \; | \; \di\Di \bm{T} \in \H^{-1}(\vol)  \} \, , \\
    \HcD{,\vol} &= \{ \bm{T} \in  \Le(\vol) \otimes \sl(3) \; | \; \curl \Di \bm{T} \in [\H^{-1}(\vol)]^3  \} \, ,
\end{align}
where $\H^{-1}(\vol)$ is the dual space of $\Honez(\vol)$. Their two-dimensional counterparts read
\begin{align}
    \HrR{,\surf} &= \{ \bm{T} \in  \Le(\surf) \otimes \Sym(2) \; | \; \rot\Rot \bm{T} \in \H^{-1}(\surf)  \} \, , \\
    \HdD{,\surf} &= \{ \bm{T} \in  \Le(\surf)\otimes \Sym(2) \; | \; \rot\Rot \bm{T} \in \H^{-1}(\surf)  \} \, , \\
    \HrD{,\surf} &= \{ \bm{T} \in  \Le(\surf) \otimes \sl(2) \; | \; \rot\Di \bm{T} \in \H^{-1}(\surf) \} \, . 
\end{align}

\subsection{Templates on the reference triangle and tetrahedron} \label{sec:recap}
The polytopal construction methodology \cite{sky_polytopal_2022,sky_higher_2023} is based on the association of a tensor-valued polynomial base function to a respective polytope of a simplex, such as a triangle or a tetrahedron. The association takes into account the support of the underlying scalar polynomial function and the projection of the tensor-valued construction on the tangential or normal vectors of the polytope. \textit{\textbf{Together, these two characteristics define the connectivity of the base function in the physical finite element mesh}}. Further, since in the construction the underlying scalar polynomial basis is linearly independent and multiplied with a linearly independent set of tensors at each polytope, the linear independence of the tensorial basis is intrinsically ensured. 
The polytopes of the reference triangle
\begin{align}
    \Gamma = \{ (\xi, \eta) \in [0,1]^2 \; | \; \xi + \eta \leq 1 \} \subset \R^2 \, , 
\end{align}
or tetrahedron
\begin{align}
    \Omega = \{ (\xi, \eta, \zeta) \in [0,1]^3 \; | \; \xi + \eta + \zeta \leq 1 \} \subset \R^3 \, , 
\end{align}
are identified using multi-indices.
\begin{definition}[Polytopes on the reference simplex]
The multi-indices of each reference simplex are with respect to its vertices.
\begin{itemize}
    \item There are three vertices $v_i$ on the reference triangle $i \in \{0,1,2\}$ and four vertices $v_i$ on the reference tetrahedron $i \in \{0,1,2,3\}$
    \item The reference triangle is equipped with multi-indices for its three edges $j \in \mathcal{J} = \{(0,1),(0,2),(1,2)\}$. Analogously, the reference tetrahedron has six edges $j \in \mathcal{J} = \{(0,1),(0,2),(0,3),(1,2),(1,3),(2,3)\}$.
    \item For triangles, the cell is defined via $c_{012}$ with the multi-index $(0,1,2)$. The three-dimensional tetrahedron has four faces $k \in \mathcal{K} =  \{(0,1,2),(0,1,3),(0,2,3),(1,2,3)\}$. 
    \item Finally, the cell of the tetrahedron $c_{0123}$ is identified with the multi-index $(0,1,2,3)$.
\end{itemize}
\end{definition}

\begin{figure}
		\centering
		\begin{subfigure}{0.24\linewidth}
			\centering
			\includegraphics[width=0.95\linewidth]{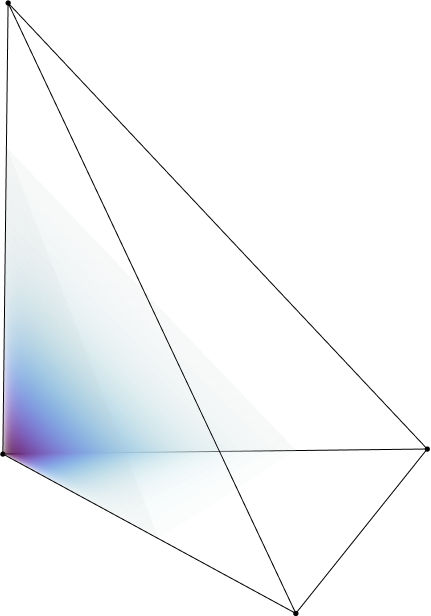}
			\caption{}
		\end{subfigure}
	    \begin{subfigure}{0.24\linewidth}
	    	\centering
	    	\includegraphics[width=0.95\linewidth]{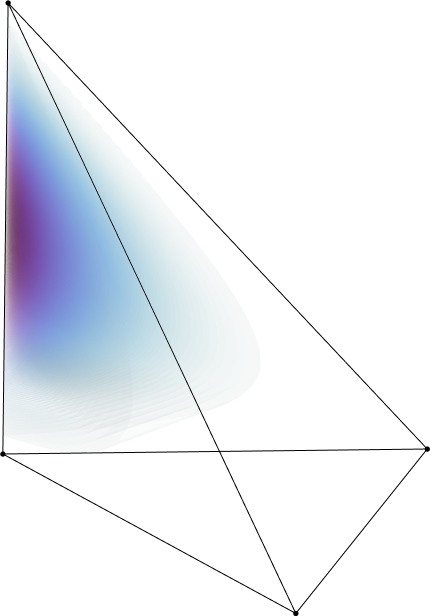}
	    	\caption{}
	    \end{subfigure}
        \begin{subfigure}{0.24\linewidth}
        	\centering
        	\includegraphics[width=0.95\linewidth]{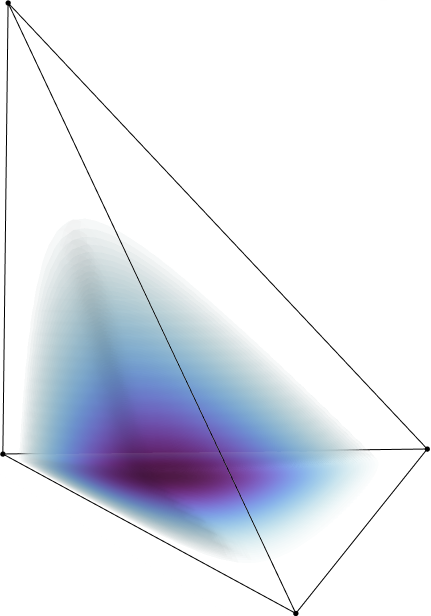}
        	\caption{}
        \end{subfigure}
        \begin{subfigure}{0.24\linewidth}
        	\centering
        	\includegraphics[width=0.95\linewidth]{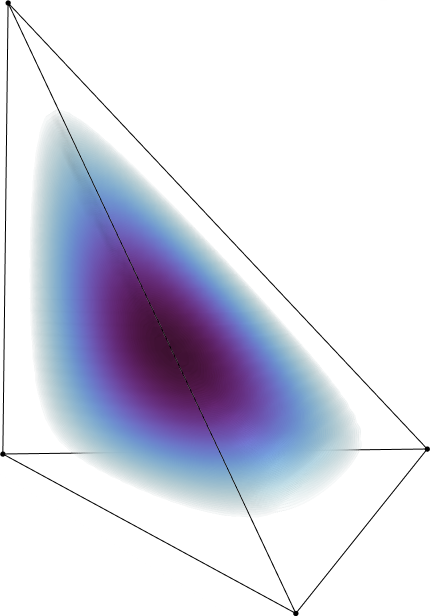}
        	\caption{}
        \end{subfigure}
		\caption{Vertex base function $n \in \ver_i^p(\Omega)$ of the vertex at the origin (a). Edge base function $n \in \edge_j^p(\Omega)$ on the left-most edge (b). Face base function $n \in \face_k^p(\Omega)$ on the bottom face (c). Finally, a cell base function $n \in \cell_{0123}^p(\Omega)$ (d). The functions belong to $\CG^p(\Omega)$ and are depicted on the reference tetrahedron $\Omega \subset \R^3$.~\ref{fn:ccby}}
		\label{fig:beziertet}
	\end{figure}

\begin{figure}
		\centering
		\input{figs/permut}
		\caption{Derivation of template vectors for the remaining edges from the first definition via permutations of the reference triangle using covariant Piola mappings. The depiction exemplifies how the first vertex-edge template vector $\tv$ corresponding to the vertex $v_0$ and edge $e_{01}$, is used to derive the vertex-edge template vectors of $v_0$-$e_{02}$ and $v_2$-$e_{12}$. Note that the permutation is always of the original reference triangle.~\ref{fn:ccby}}
		\label{fig:permut}
\end{figure}
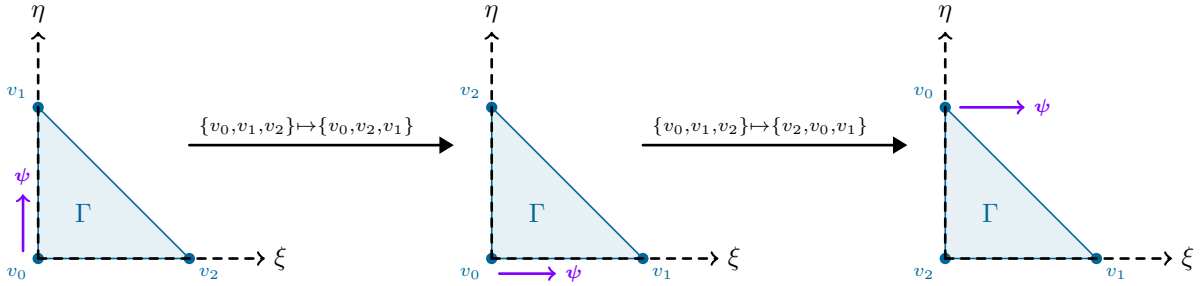

Further, we give a general definition of a polytopal base function for various trace operators.   
\begin{definition} [Simplex polytopal base functions] \label{def:poly}
	Each base function is associated with its respective polytope via a specified trace operator.
	\begin{enumerate}
			\item A vertex base function has a  vanishing trace on all other vertices and non-neighbouring edges and faces.
			\item An edge base function has a vanishing trace on all other edges and non-neighbouring faces. 
			\item A face base function has a vanishing trace on all other faces.
			\item A cell base function has a vanishing trace on the entire boundary of the element.
	\end{enumerate}
	The definition is general and the respective trace operator may change.
\end{definition}
Finally, we define for each polytope of the simplex the corresponding scalar base function based on its support. For simplicity, we assume complete polynomial spaces $\Po^p(T)$. However, the construction generalises to hierarchical anisotropic definitions as well. 
\begin{definition}[Simplex polytopal scalar spaces]
Each polytope is associated with a space of scalar base functions. Further, each polytope is identified by a multi-index of its underlying vertices. The complete space is a $\C^0(\body)$-continuous Lagrange-type $\CG^p(T)$, having the same dimension and span as $\Po^p(T)$ element-wise.
For triangles there holds $\dim \Po^p(\Gamma) = (p+2)(p+1)/2$, whereas $\dim \Po^p(\Omega) = (p+3)(p+2)(p+1)/6$ for tetrahedra.
\begin{itemize}
    \item Each vertex $v_i$ is associated with the space of its respective base function $\ver^p_i(T)$. As such, there are three spaces in 2D $i \in \{0,1,2\}$ and four spaces in 3D $i \in \{0,1,2,3\}$ and each one is of dimension one, $\dim \ver^p_i(T) = 1 \quad \forall \, i$. The index $p$ refers to the polynomial order of the space, which is always $1$ if a hierarchical polynomial basis is employed.
    \item For each edge $e_j$ there exists a space of edge functions that vanish on its respective vertices $\edge^p_{j}(T)$ with $j \in \mathcal{J} = \{(0,1),(0,2),(1,2)\}$ in 2D and $j \in \mathcal{J} = \{(0,1),(0,2),(0,3),(1,2),(1,3),(2,3)\}$ in 3D.
    The dimension of each edge space is given by $\dim \edge^p_j(T) = p-1$.
    \item In 2D the face $f_k$ is the cell of the triangle $\cell_{012}^p(T) := \face^p(T)$, such that its base functions vanish on the entire boundary of the triangle. In 3D there are four face spaces $\face_{k}^p(T)$ with $k \in \mathcal{K} =  \{(0,1,2),(0,1,3),(0,2,3),(1,2,3)\}$. The dimension of the spaces reads $\dim \face_{k}^p(T) = (p-2)(p-1)/2$. The base functions of each face vanish on its respective edges.
    \item Lastly, for three-dimensional tetrahedra the space of cell base functions is given by $\cell_{0123}^p(T)$ with the dimensionality $\dim \cell_{0123}^p(T) = (p-3)(p-2)(p-1)/6$. The cell base functions vanish on the entire boundary of the tetrahedron.
\end{itemize}
The association with a respective polytope is according to \cref{def:poly} via the trace operator $\tr f = f |_\vb{x}$, of evaluation at an interface, and is visualised in \cref{fig:beziertet}\renewcommand{\thefootnote}{*}
\footnote{\label{fn:ccby} This figure is adapted from the original under CC-BY license from \cite{sky_polytopal_2022}.}
on the reference tetrahedron $\Omega$.
\end{definition}
In essence, there are two main types of possible projections between interfacing elements, namely tangential or normal. Evidently, these also define the partial continuity of the corresponding finite element space, i.e., a tangential-continuous or a normal-continuous finite element. The definition of a vectorial basis satisfying these continuity requirements can be done for each polytope of the reference element by an initial definition on a chosen vertex, edge, face and the cell. The initial definition is then mapped to the remaining polytopes via Piola transformations, see \cref{ap:a}. This is achieved by permuting the order of the vertices on the reference element, such that multiple affine maps arise, and subsequently applying the transformation of the initial basis. The procedure is depicted in \cref{fig:permut}~\ref{fn:ccby}. We refer to \cite{sky_polytopal_2022} for an exhaustive exposition of the procedure. 
In the tensorial case the construction of the base functions is often more involved due to additional constraints, such as symmetry or tracelessness. This work focuses on tensor-valued finite elements, specifically targeting dyads.
In fact, the tensorial bases in this work are defined using dyadic products of the vectorial bases from \cite{sky_polytopal_2022}. 

We briefly expand on the ordering used in the definition of vectorial bases on each polytope, hereinafter referred to as template sets. On each lower dimensional polytope we use the minimal index ordering with respect to its intersecting higher dimensional polytopes. This ordering is directly related to the \textit{\textbf{connectivity}} of the corresponding base functions. 

We start with the ordering on two-dimensional triangles. Let $\mathcal{J}_i$ be the \textit{\textbf{index-ordered set}} of intersecting edges $e_j$ for the vertex $v_i$ with $j \in \mathcal{J}_i$
\begin{align}
    &\mathcal{J}_i = \{ j_0,j_1 \} \subset \mathcal{J} && \text{s.t.} && i \in j_0 \cap j_1 \, , \quad  j_0 < j_1 \, ,  
\end{align}
then template vectors in the \textit{\textbf{vertex-template set}} $\tem_i$ are associated with the edges of the index-pairs $j \in \mathcal{J}_i$ in the same order. 
The association is either with respect to the tangential projection or the normal projection on the edge tangent or normal vector, depending on the desired type of continuity for the construction. On the edges, the template sets $\tem_j$ are defined such that the first vector in the set is associated with the edge and the second vector is associated with the cell. The cell template set $\tem_{012}$ is associated solely with the cell.
\textit{\textbf{Let the operator $\sim$ imply "associated with" with respect to the connectivity of the vectors of a template set for the specified polytopes}}, an illustrative association for both tangential template sets $\tem^{\, \tau}$, and normal template sets $\tem^{\, \nu}$, can be given with 
\begin{align}
    \tem_0 &\sim \{ e_{01}, e_{02} \} \, , & \tem_1 &\sim \{ e_{01}, e_{12} \} \, , & \tem_2 &\sim \{ e_{02}, e_{12} \} \, , \notag \\
    \tem_{01} &\sim \{ e_{01}, c_{012} \} \, , & \tem_{02} &\sim \{ e_{02}, c_{012} \} \, , & \tem_{12} &\sim \{ e_{12}, c_{012} \} \, , \\
     & & \tem_{012} &\sim \{ c_{012},c_{012} \} \, . &  & \notag
\end{align}
Note that every template set is two-dimensional since the vectorial space is two-dimensional $\R^2$. Edge templates are associated with the edge and the cell since for one of the vectors in the set either the tangential or normal trace is zero, compare \cref{eq:tem2Dtansets} and \cref{eq:tem2Dn}. 

On the three-dimensional tetrahedron the association is more complex. For tangential continuity, the association of the vertex sets $\tem_i$ is analogous to the two-dimensional case, simply adjusted to three intersecting edges $e_j$ rather than two
\begin{align}
    &\mathcal{J}_i = \{ j_0,j_1,j_2 \} \subset \mathcal{J} && \text{s.t.} && i \in j_0 \cap j_1 \cap j_2 \, , \quad  j_0 < j_1 < j_2 \, .  
\end{align}
However, the edge template sets $\tem_j$ differ. For each edge $e_j$, which is given by the intersection of two faces $k \in \mathcal{K}_j$ of the \textit{\textbf{ordered}} set $\mathcal{K}_j$ 
\begin{align}
    &\mathcal{K}_j = \{ k_0,k_1 \} \subset \mathcal{K} && \text{s.t.} && j \in k_0 \cap k_1  \, , \quad  k_0 < k_1 \, ,  \label{eq:kjface}
\end{align}
the first edge template vector is associated with the edge itself, and the two remaining vectors are associated with the intersecting faces with the same order as $k$. On each face $f_k$ with $k \in \mathcal{K}$, the first two vectors of the set $\tem_k$ are associated with the face and the last vector is associated with the cell. Finally, the vectors of the cell template $\tem_{0123}$ are only associated with the cell. We illustrate the association for tangential template sets $\tem^{\, \tau}$ with 
\begin{align}
    \tem_0^{\, \tau} &\sim \{ e_{01}, e_{02}, e_{03} \} \, , & \tem_1^{\, \tau} &\sim \{ e_{01}, e_{12}, e_{13} \} \, , & \tem_2^{\, \tau} &\sim \{ e_{02}, e_{12}, e_{23} \} \, , 
    & \tem_3^{\, \tau} &\sim \{ e_{03}, e_{13}, e_{23} \} \, ,
    \notag \\
    \tem_{01}^{\, \tau} &\sim \{ e_{01}, f_{012}, f_{013} \} \, , & \tem_{02}^{\, \tau} &\sim \{ e_{02}, f_{012}, f_{023} \} \, , & \tem_{03}^{\, \tau} &\sim \{ e_{03}, f_{013},f_{023} \} \, , \notag \\
    \tem_{12}^{\, \tau} &\sim \{ e_{12}, f_{012}, f_{123} \} \, , & \tem_{13}^{\, \tau} &\sim \{ e_{13}, f_{013}, f_{123} \} \, , & \tem_{23}^{\, \tau} &\sim \{ e_{23}, f_{023},f_{123} \} \, ,  \\
    \tem_{012}^{\, \tau} &\sim \{ f_{012}, f_{012}, c_{0123} \} \, , & \tem_{013}^{\, \tau} &\sim \{ f_{013}, f_{013}, c_{0123} \} \, , & \tem_{023}^{\, \tau} &\sim \{ f_{023}, f_{023}, c_{0123} \} \, , 
    & \tem_{123}^{\, \tau} &\sim \{ f_{123}, f_{123}, c_{0123} \} \, ,
    \notag \\
     & & \tem_{0123}^{\, \tau} &\sim \{ 
c_{0123},c_{0123},c_{0123} \} \, . &  & \notag
\end{align}
The dimension of the template sets is now three due to the fact that vectorial space is three-dimensional $\R^3$. An edge template set is associated with itself and its two intersecting faces, since it defines a vector with a tangential trace on the edge itself, and two vectors that are orthogonal to the edge tangent but produce a tangential trace on each of its respective faces, compare \cref{eq:tem3Dtansets}.

For normal continuity, a vertex $v_i$ given by the intersection of three faces $k \in \mathcal{K}_i$ of the ordered set $\mathcal{K}_i$ is equipped with the three vector set $\tem_i$ associated with the faces in the same order as $k$ 
\begin{align}
    &\mathcal{K}_i = \{ k_0,k_1,k_2 \} \subset \mathcal{K} && \text{s.t.} && i \in k_0 \cap k_1 \cap k_2 \, , \quad  k_0 < k_1 < k_2 \, .  
\end{align}
On each edge $e_j$ given by the intersection of two faces $k \in \mathcal{K}_j$ of the ordered set $\mathcal{K}_j$ as per \cref{eq:kjface}, the first two vectors of the template set $\tem_j$ are associated with the faces in the same order as $k$, while the last vector is associated with the cell. Each face $f_k$ with $k \in \mathcal{K}$ is equipped with a set $\tem_k$ of one face vector and two cell vectors in that order. Finally, the cell template set $\tem_{0123}$ is purely associated with cell.  
We illustrate the association for normal template sets $\tem^{\, \nu}$ with
\begin{align}
    \tem_0^{\, \nu} &\sim \{ f_{012}, f_{013}, f_{023} \} \, , & \tem_1^{\, \nu} &\sim \{ f_{012}, f_{013}, f_{123} \} \, , & \tem_2^{\, \nu} &\sim \{ f_{012}, f_{023}, f_{123} \} \, , 
    & \tem_3^{\, \nu} &\sim \{ f_{013}, f_{023}, f_{123} \} \, ,
    \notag \\
    \tem_{01}^{\, \nu} &\sim \{ f_{012}, f_{013}, c_{0123} \} \, , & \tem_{02}^{\, \nu} &\sim \{ f_{012}, f_{023}, c_{0123} \} \, , & \tem_{03}^{\, \nu} &\sim \{ f_{013}, f_{023},c_{0123} \} \, , \notag \\
    \tem_{12}^{\, \nu} &\sim \{ f_{012}, f_{123}, c_{0123} \} \, , & \tem_{13}^{\, \nu} &\sim \{ f_{013}, f_{123}, c_{0123} \} \, , & \tem_{23}^{\, \nu} &\sim \{ f_{023}, f_{123},c_{0123} \} \, ,  \\
    \tem_{012}^{\, \nu} &\sim \{ f_{012}, c_{0123}, c_{0123} \} \, , & \tem_{013}^{\, \nu} &\sim \{ f_{013}, c_{0123}, c_{0123} \} \, , & \tem_{023}^{\, \nu} &\sim \{ f_{023}, c_{0123}, c_{0123} \} \, , 
    & \tem_{123}^{\, \nu} &\sim \{ f_{123}, c_{0123}, c_{0123} \} \, ,
    \notag \\
     & & \tem_{0123}^{\, \nu} &\sim \{ c_{0123},c_{0123},c_{0123} \} \, . &  & \notag
\end{align}
Here, each edge template set is associated with its two intersecting faces, such that two of its template vectors produce non-vanishing normal trace on each respective face, while the last template vector is the edge tangent vector, which produces no normal trace on either faces, compare \cref{eq:tem3Dnorsets}.  
With the association concluded, we can now move on to the definition of the template sets. We note that the orientation of each template vector is according to \cref{ap:orient}.  

\begin{figure}
		\centering
		\input{figs/tri_nii}
		\caption{Template vectors for the construction of tangential-continuous base functions on the reference triangle on their corresponding polytope. Each vertex is endowed with two template vectors, one for each of its intersecting edges. Each edge is equipped with two template vectors, one for its tangent and one normal edge-cell vector for the cell. Finally, the cell is endowed with the Cartesian basis.~\ref{fn:ccby}}
		\label{fig:tri_nii}
\end{figure}
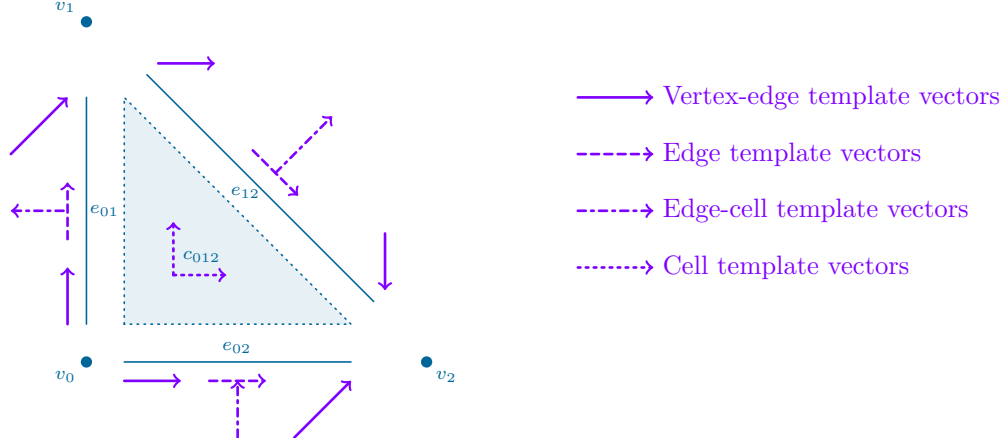

In order to construct a purely tangential-continuous finite element on triangulations we introduce a set of template vectors on the reference triangle, associated with its respective polytopes
\begin{align}
		\tem^{\, \tau}_{2D} = \{\tem_0^{\, \tau},\tem_1^{\, \tau},\tem_2^{\, \tau},\tem_{01}^{\, \tau},\tem_{02}^{\, \tau},\tem_{12}^{\, \tau},\tem_{012}^{\, \tau}\} \, , 
  \label{eq:tem2Dt}
	\end{align}
where the polytopal sets read
\begin{align}
		\tem_0^{\, \tau} &= \{\vb{e}_2,\vb{e}_1\} \, , & \tem_1^{\, \tau} &= \{\bm{\iota}^\tau_1,\vb{e}_1\} \, , & \tem_2^{\, \tau} &= \{\bm{\iota}^\tau_1,-\vb{e}_2\} \, , \notag \\
		\tem_{01}^{\, \tau} &= \{\vb{e}_2,-\vb{e}_1\} \, , & \tem_{02}^{\, \tau} &= \{\vb{e}_1,\vb{e}_2\} \, , & \tem_{12}^{\, \tau} &= \{ \bm{\iota}^\tau_2,\bm{\iota}^\tau_1\} \, , \notag \\
		\tem_{012}^{\, \tau} &= \{\vb{e}_1,\vb{e}_2\} \, ,
        \label{eq:tem2Dtansets}
	\end{align}
using the definitions $\bm{\iota}^\tau_1 = \vb{e}_1 + \vb{e}_2$ and $\bm{\iota}^\tau_2 = (1/2) (\vb{e}_1 - \vb{e}_2)$ for conciseness.
The template vectors are depicted in \cref{fig:tri_nii}~\ref{fn:ccby}. On each non-cell polytope of the reference triangle the template vectors define a basis for $\R^2$, such that their tangential projection $\con{\bm{\tau}}{\tv}$ is limited to the tangent vector of a single edge. Consequently, they are suited to define a tangential-continuous finite element space.  

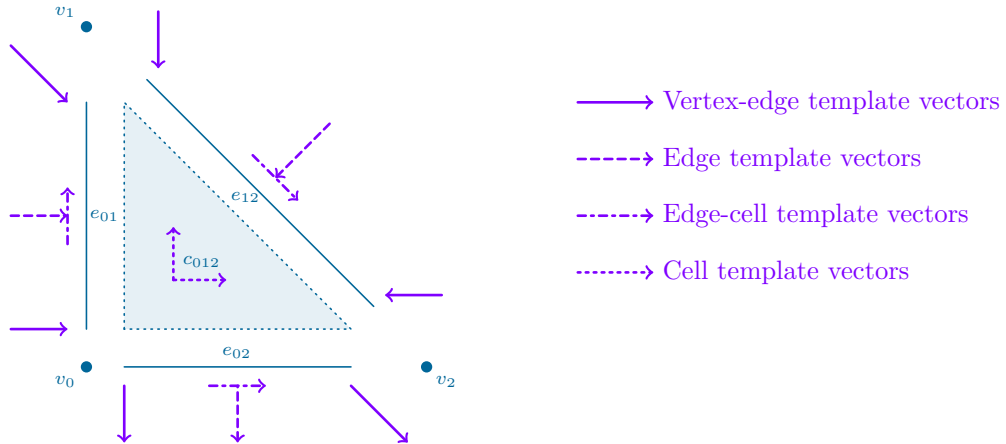
\begin{figure}
		\centering
		\input{figs/tri_bdm}
		\caption{Template vectors for the construction of normal-continuous base functions on the reference triangle on their corresponding polytope. Each vertex is endowed with two template vectors, one for each of its intersecting edges and their respective normals. Each edge is equipped with two template vectors, one for its normal and one tangent edge-cell vector for the cell. Finally, the cell is endowed with the Cartesian basis.~\ref{fn:ccby}}
		\label{fig:tri_bdm}
\end{figure}

Analogously, one can define the template super set for normal-continuity on the reference triangle 
\begin{align}
		\tem^{\, \nu}_{2D} = \{\tem_0^{\, \nu},\tem_1^{\, \nu},\tem_2^{\, \nu},\tem_{01}^{\, \nu},\tem_{02}^{\, \nu},\tem_{12}^{\, \nu},\tem_{012}^{\, \nu}\} \, , 
  \label{eq:norm2dtem}
\end{align}
such that the template vectors have a normal projection limited to one edge $\con{\bm{\nu}}{\tv}$. 
We define $\bm{\iota}^\nu_1 = \vb{e}_1 - \vb{e}_2$, $\bm{\iota}^\nu_2 = -(1/2) (\vb{e}_1 + \vb{e}_2)$ and the polytopal template sets
	\begin{align}
		\tem_0^{\, \nu} &= \{\vb{e}_1,-\vb{e}_2\} \, , & \tem_1^{\, \nu} &= \{\bm{\iota}^\nu_1,-\vb{e}_2\} \, , & \tem_2^{\, \nu} &= \{\bm{\iota}^\nu_1, -\vb{e}_1\} \, , \notag \\
		\tem_{01}^{\, \nu} &= \{\vb{e}_1,\vb{e}_2\} \, , & \tem_{02}^{\, \nu} &= \{-\vb{e}_2,\vb{e}_1\} \, , & \tem_{12}^{\, \nu} &= \{ \bm{\iota}^\nu_2,\bm{\iota}^\nu_1\} \, , \notag \\
		\tem_{012}^{\, \nu} &= \{\vb{e}_1,\vb{e}_2\} \, , \label{eq:tem2Dn}
\end{align}
which are depicted in \cref{fig:tri_bdm}~\ref{fn:ccby}.
\textbf{\textit{Clearly, these template vectors can be retrieved by rotating the tangential template vectors by 90 degrees}}. However, this does not apply to the three-dimensional tetrahedron.

\begin{figure}
		\centering
		\input{figs/tet_nii}
		\caption{Template vectors for the construction of tangential-continuous base functions on the reference tetrahedron on their corresponding polytopes. Only vectors on the visible sides of the tetrahedron are depicted.~\ref{fn:ccby}}
		\label{fig:tet_nii}
	\end{figure}
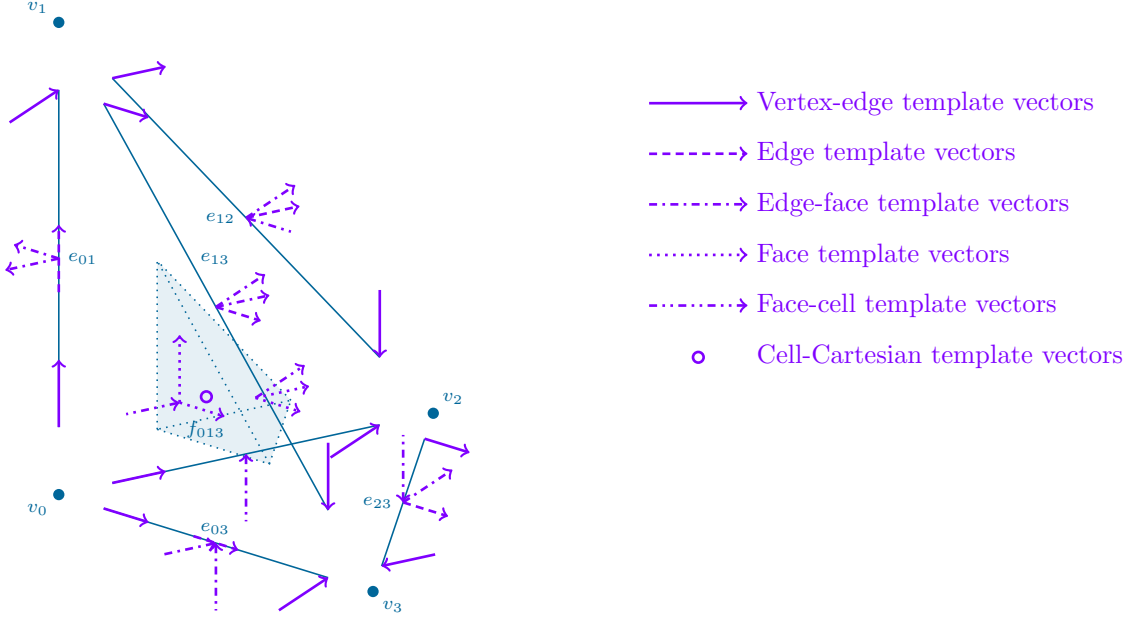

The tangential template of the reference tetrahedron reads
\begin{align}
		\tem^{\, \tau}_{3D} = \{ \tem_0^{\, \tau},\tem_1^{\, \tau},\tem_2^{\, \tau},\tem_{3}^{\, \tau},\tem_{01}^{\, \tau},\tem_{02}^{\, \tau},\tem_{03}^{\, \tau},\tem_{12}^{\, \tau},\tem_{13}^{\, \tau},\tem_{23}^{\, \tau},\tem_{012}^{\, \tau},\tem_{013}^{\, \tau},\tem_{023}^{\, \tau},\tem_{123}^{\, \tau},\tem_{0123}^{\, \tau} \} \, ,
  \label{eq:temtan3d}
\end{align}
given by the template sets
\begin{align}
    	\tem_0^{\, \tau} &= \{ \vb{e}_3,\vb{e}_2,\vb{e}_1 \} \, , & \tem_1^{\, \tau} &= \{ \bm{\iota}^\tau_1 , \vb{e}_2 , \vb{e}_1 \} \, , & \tem_2^{\, \tau} &= \{ \bm{\iota}^\tau_1, -\vb{e}_3 ,\vb{e}_1 \} \, , \notag \\
    	\tem_{3}^{\, \tau} &= \{ \bm{\iota}^\tau_1, -\vb{e}_3, -\vb{e}_2\} \, , & \tem_{01}^{\, \tau} &= \{ \vb{e}_3, -\vb{e}_2, -\vb{e}_1 \} \, , & \tem_{02}^{\, \tau} &= \{ \vb{e}_2, \vb{e}_3, -\vb{e}_1 \} \, , \notag \\
    	\tem_{03}^{\, \tau} &= \{ \vb{e}_1, \vb{e}_3, \vb{e}_2 \} \, , & \tem_{12}^{\, \tau} &= \{ \vb{e}_2, \bm{\iota}^\tau_1, -\vb{e}_1 \} \, , & \tem_{13}^{\, \tau} &= \{ \vb{e}_1, \bm{\iota}^\tau_1, \vb{e}_2 \} \, , \notag \\
    	\tem_{23}^{\, \tau} &= \{ \vb{e}_1, \bm{\iota}^\tau_1, -\vb{e}_3 \} \, , & \tem_{012}^{\, \tau} &= \{ \vb{e}_3, \vb{e}_2, -\vb{e}_1 \} \, , & \tem_{013}^{\, \tau} &= \{ \vb{e}_3, \vb{e}_1, \vb{e}_2 \} \, , \notag \\
    	\tem_{023}^{\, \tau} &= \{\vb{e}_2, \vb{e}_1, -\vb{e}_3\} \, , & \tem_{123}^{\, \tau}  &= \{ \vb{e}_2, \vb{e}_1,\bm{\iota}^\tau_1\} \, , & \tem_{0123}^{\, \tau} &= \{ \vb{e}_3, \vb{e}_2, \vb{e}_1 \} \, , 
     \label{eq:tem3Dtansets}
    \end{align}
where we employed the definition $\bm{\iota}^\tau_1 = \vb{e}_1 + \vb{e}_2 + \vb{e}_3$. The template vectors are illustrated in \cref{fig:tet_nii}~\ref{fn:ccby}. Each of the vertex template vectors produces a non-vanishing tangential projection $\con{\bm{\tau}}{\tv}$ of a single edge of the edges that intersect the vertex. For each vector set on an edge, the first vector in the set produces a tangential projection on the edge tangent vector, whereas the two remaining vectors do not. However, the two remaining vectors produce each a tangential projection on one of the respective faces that intersect at the edge. On each face, the first two vectors produce a tangential projection on the plane of the face, whereas the remaining vector is orthogonal to it. Thus, the set is suited to define a tangential-continuous finite element space.

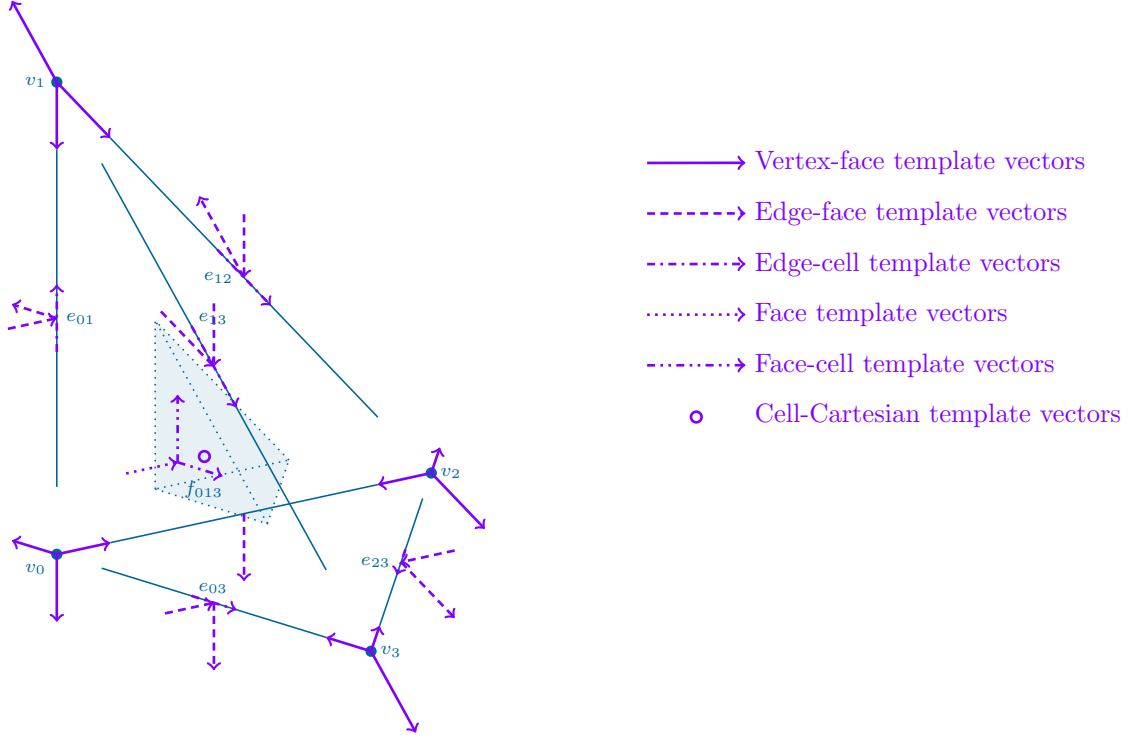
\begin{figure}
		\centering
		\input{figs/tet_bdm}
		\caption{Template vectors for the construction of normal-continuous base functions on the reference tetrahedron on their corresponding polytopes. Only vectors on the visible sides of the tetrahedron are depicted.~\ref{fn:ccby}}
		\label{fig:tet_bdm}
	\end{figure}

The super-set for normal continuity is given by
	\begin{align}
		\tem^{\, \nu}_{3D} = \{ \tem_0^{\, \nu},\tem_1^{\, \nu},\tem_2^{\, \nu},\tem_{3}^{\, \nu},\tem_{01}^{\, \nu},\tem_{02}^{\, \nu},\tem_{03}^{\, \nu},\tem_{12}^{\, \nu},\tem_{13}^{\, \nu},\tem_{23}^{\, \nu},\tem_{012}^{\, \nu},\tem_{013}^{\, \nu},\tem_{023}^{\, \nu},\tem_{123}^{\, \nu},\tem_{0123}^{\, \nu} \} \, ,
  \label{eq:tem3Dnormal}
	\end{align}
 and is composed of the following subsets
\begin{align}
    	\tem_0^{\, \nu} &= \{ -\vb{e}_1,\vb{e}_2,-\vb{e}_3 \} \, , & \tem_1^{\, \nu} &= \{ \bm{\iota}^\nu_1,  \bm{\iota}^\nu_2, -\vb{e}_3 \} \, , & \tem_2^{\, \nu} &= \{ \bm{\iota}^\nu_3, \bm{\iota}^\nu_2, -\vb{e}_2 \} \, , \notag \\
    	\tem_{3}^{\, \nu} &= \{ \bm{\iota}^\nu_3, -\bm{\iota}^\nu_1 , -\vb{e}_1 \} \, , & \tem_{01}^{\, \nu} &= \{ -\vb{e}_1, \vb{e}_2, \vb{e}_3 \} \, , & \tem_{02}^{\, \nu} &= \{ -\vb{e}_1, -\vb{e}_3, \vb{e}_2 \} \, , \notag \\
    	\tem_{03}^{\, \nu} &= \{ \vb{e}_2, -\vb{e}_3 ,\vb{e}_1 \} \, , & \tem_{12}^{\, \nu} &= \{ \bm{\iota}^\nu_1 , -\vb{e}_3, \bm{\iota}^\nu_2  \} \, , & \tem_{13}^{\, \nu} &= \{ \bm{\iota}^\nu_2, -\vb{e}_3, -\bm{\iota}^\nu_1 \} \, , \notag \\
    	\tem_{23}^{\, \nu} &= \{ \bm{\iota}^\nu_2, -\vb{e}_2, -\bm{\iota}^\nu_3 \} \, , & \tem_{012}^{\, \nu} &= \{ -\vb{e}_1, \vb{e}_3, \vb{e}_2 \} \, , & \tem_{013}^{\, \nu} &= \{ \vb{e}_2, \vb{e}_3, \vb{e}_1 \} \, , \notag \\
    	\tem_{023}^{\, \nu} &= \{ -\vb{e}_3, \vb{e}_2, \vb{e}_1 \} \, , & \tem_{123}^{\, \nu}  &= \{ -\vb{e}_3, \bm{\iota}^\nu_2 , -\bm{\iota}^\nu_1 \} \, , & \tem_{0123}^{\, \nu} &= \{\vb{e}_3, \vb{e}_2, \vb{e}_1 \} \, ,
     \label{eq:tem3Dnorsets}
    \end{align}
where we used the definitions $\bm{\iota}^\nu_1 = \vb{e}_3 - \vb{e}_1$, $\bm{\iota}^\nu_2 = \vb{e}_2 - \vb{e}_3$ and $\bm{\iota}^\nu_3 = \vb{e}_2 - \vb{e}_1$. The vectors are depicted in \cref{fig:tet_bdm}~\ref{fn:ccby}. For each vertex set, the respective vectors produce a normal projection $\con{\bm{\nu}}{\tv}$ on only one of faces that intersect at that vertex. For each edge set, the first two vectors each produce a normal projection on only one of faces that intersect at that edge. The remaining vector is tangential to both faces. In each face template set, only the first vector in the set produces a normal projection on the plane of the face, whereas the remaining vectors are tangential to it. Therefore, the bases are ideal for the construction of a normal-continuous finite element space.  

With the fundamental machinery in place, we can now move on to construct various finite element spaces. For each finite element we introduce the definition of the base functions using the scalar base functions from $\CG^p(T)$ and polytopal template sets. \textit{\textbf{As a convention, we introduce a Greek index for the association with the polytope and Latin indices for the association with the scalar functions and the template tensors. Observe that the occurrence of a Latin index implies more than one base function on the same polytope}}. Concrete examples for the base functions of two elements are given in \cref{ap:c}.       

\begin{remark}[Non-hierarchical bases with non-affine mappings]
    We note that for some of our subsequent constructions there is no guarantee that the space of constants will be captured by the element \textit{\textbf{if a non-hierarchical polynomial basis and non-affine mappings}} of the reference simplex to the physical simplex are simultaneously employed. This is so, since in these constructions the vertex base functions, edge base functions, face base functions and cell base functions undergo different transformations. Consequently, if the basis is characterised by the partition of unity property, this characteristic may be lost in the varying transformations. For a hierarchical basis there is no risk, as the vertex base functions span the constant space and transform in the same manner. 
\end{remark}

\section{Tangential-continuous elements} \label{sec:tan}

\subsection{The Regge element}
The Regge element \cite{Li2018,Chr2011,HH2013} is defined to allow for symmetric tangential-tangential continuity, and can be used to discretise fields in $\HcC{,\vol}$ or $\HrR{,\surf}$ in three dimensions or two dimensions, respectively. To clarify, the Regge element is almost conforming such that, $\Reg^p(\vol) \subsetsim \HcC{,\vol}$ or $\Reg^p(\surf)\subsetsim \HrR{,\surf}$. The element spans an element-wise polynomial space of symmetric tensors, whose tangential-tangential components are continuous across interfaces $\Xi$. All other tensorial components are allowed to jump between simplices. The corresponding trace operator for Regge elements is therefore 
\begin{align}
    \tr_{tt} \bm{P} &= \con{\vb{t}\otimes\vb{t}}{\bm{P}} = P_{tt} \, , 
\end{align}
and there must hold $\jump{\tr_{tt}\bm{P}}|_{\Xi} = 0$ on all arbitrarily defined interfaces $\Xi$ in the domain for a conforming construction.
Given a tensorial base function on the reference simplex, its tangential-tangential continuity is clearly preserved by the double covariant Piola transformation 
\begin{align}
    \con{\vb{t} \otimes \vb{t}}{\bm{Y}} = \con{\bm{J}\bm{\tau} \otimes \bm{\tau}\bm{J}^T}{\bm{Y}} = \con{\bm{\tau} \otimes \bm{\tau}}{\bm{\Upsilon}} \qquad \iff \qquad \bm{Y} = \bm{J}^{-T} \bm{\Upsilon} \bm{J}^{-1} \, . 
\end{align}
In the following we define the base functions of the Regge element on the reference simplices. The base functions are mapped to their physical counterparts in the finite element mesh using the double covariant Piola transformation.

\subsubsection{The Regge triangle element}
In order to define the construction of base functions for the Regge element on the reference triangle we use the tangential vectorial template \cref{eq:tem2Dt}. Let $\bm{\psi}$ represent template vectors from the vectorial sets, the tensorial sets are retrieved by symmetrising the dyadic products of the vectorial sets 
\begin{align}
    &\tem_\alpha^{\, \tau \tau} = \sym(\tem_\alpha^{\, \tau} \otimes \tem_\alpha^{\, \tau}) = \{ \tv_1 \otimes \tv_1,\tv_2 \otimes \tv_2,\sym(\tv_1 \otimes \tv_2) \} \, , && \dim  \tem_\alpha^{\, \tau \tau} = \dim \Sym(2) = 3 \, ,
\end{align}
where $\tem_\alpha^{\, \tau} \in \tem^{\, \tau}_{2D}$ with the multi-index $\alpha \in \{0,1,2\} \cup \mathcal{J} \cup \{(0,1,2)\}$. 
The superset $\tem^{\, \tau \tau}_{2D}$ contains all the resulting polytopal template subsets $\tem_\alpha^{\, \tau \tau}$
\begin{align}
    \tem^{\, \tau \tau}_{2D} = \{\tem_0^{\, \tau \tau},\tem_1^{\, \tau \tau},\dots, \tem_{01}^{\, \tau \tau},\dots,\tem_{012}^{\, \tau \tau}\}  \, .
\end{align}
The Regge element on the reference triangle can now be easily defined as 
\begin{align}
    &\Reg^p(\Gamma) = \left \{ \bigoplus_{i = 0}^2 \ver_i^p(\Gamma) \otimes \tem_i^{\, \tau \tau} \right \} \oplus \left \{ \bigoplus_{j \in \mathcal{J}} \edge_j^p(\Gamma) \otimes \tem_j^{\, \tau \tau} \right \} \oplus \cell_{012}^p(\Gamma) \otimes \tem_{012}^{\, \tau \tau} \, ,     
\end{align}
and the dimension of the space reads
\begin{align}
    && \dim \Reg^p(\Gamma) = \dim[\Po^p(\Gamma) \otimes \Sym(2)] = \dfrac{3(p+2)(p+1)}{2} \, ,
\end{align}
for a complete polynomial space $\Po^p(\Gamma)$.
\begin{definition}[Triangle Regge base functions]
    The base functions of the Regge triangle element are defined per polytope.
    \begin{itemize}
        \item On each edge $e_{ij}$ with multi-index $(i,j) \in \mathcal{J}$, equipped with the tangent vector $\bm{\tau}$ and trace operator $\tr_{\tau \tau} (\cdot) = \con{\bm{\tau} \otimes \bm{\tau}}{\cdot} |_{\mu_{ij}}$ we define the edge base functions
        \begin{subequations}
            \begin{align}
            &\text{Vertex-edge}: & \bm{\Upsilon}_\alpha(\xi,\eta) &= n \tt \, , & n &\in \ver^p_\alpha(\Gamma) \, , & \tt &\in \{ \tt \in \tem_\alpha^{\, \tau \tau} \; | \; \tr_{\tau \tau} \tt \neq 0 \} \, , \\
            &\text{Edge}: & \bm{\Upsilon}_{l}(\xi,\eta) &= n_l \tt \, , & n_l &\in \edge^p_{ij}(\Gamma) \, , & \tt &\in \{ \tt \in \tem_{ij}^{\, \tau \tau} \; | \; \tr_{\tau \tau} \tt \neq 0 \} \, , 
        \end{align}
        \end{subequations}
        where $\alpha \in \{i,j\}$. For each edge there are $2 \cdot 1$ vertex-edge base functions and $p-1$ edge base functions.
        \item The cell $c_{012}$ is equipped with the trace operator $\tr_{\tau \tau} (\cdot) = \con{\bm{\tau} \otimes \bm{\tau}}{\cdot} |_{\partial \Gamma}$. Its base functions read
        \begin{subequations}
        \begin{align}
            &\text{Vertex-cell}: & \bm{\Upsilon}_\alpha(\xi,\eta) &= n \tt \, , & n &\in \ver^p_\alpha(\Gamma) \, , & \tt &\in \{ \tt \in \tem_\alpha^{\, \tau \tau} \; | \; \tr_{\tau \tau} \tt = 0 \} \, , \\
            &\text{Edge-cell}: & \bm{\Upsilon}_{l  q \beta}(\xi,\eta) &= n_l \tt_q \, , & n_l &\in \edge^p_{\beta}(\Gamma) \, , & \tt_q &\in \{ \tt \in \tem_{\beta}^{\, \tau \tau} \; | \; \tr_{\tau \tau} \tt = 0 \} \, , \\
            &\text{Cell}: & \bm{\Upsilon}_{lq}(\xi,\eta) &= n_l \tt_q \, , & n_l &\in \cell^p_{012}(\Gamma) \, , & \tt_q &\in \tem_{012}^{\, \tau \tau} \, ,
        \end{align}
        \end{subequations}        
        where $\alpha \in \{0,1,2\}$ and $\beta \in \mathcal{J}$. There are $3\cdot 1$ vertex-cell base functions, $3 \cdot  2 \cdot (p-1)$ edge-cell base functions and $3 \cdot (p-2)(p-1)/2$ pure cell base functions. 
    \end{itemize}
\end{definition}

\subsubsection{The Regge tetrahedral element}
The same construction procedure is applied for the three-dimensional Regge tetrahedron. We define the template sets via
\begin{align}
    \tem_\alpha^{\, \tau \tau} &= \sym(\tem_\alpha^{\, \tau} \otimes \tem_\alpha^{\, \tau}) = \{ \tv_1 \otimes \tv_1,\tv_2 \otimes \tv_2,\tv_3 \otimes \tv_3, \sym(\tv_1 \otimes \tv_2),\sym(\tv_1 \otimes \tv_3),\sym(\tv_2 \otimes \tv_3) \} \, ,  \notag \\ \dim  \tem_\alpha^{\, \tau \tau} &= \dim \Sym(3) = 6 \, , 
    \label{eq:toreg}
\end{align}
with $\tem_\alpha^{\, \tau} \in \tem_{3D}^{\, \tau}$ using the tangential template set from \cref{eq:temtan3d}. The superset $\tem^{\, \tau \tau}_{3D}$ is constructed from all the resulting polytopal subsets with the multi-index $\alpha \in \{0,1,2,3\}\cup \mathcal{J} \cup \mathcal{K} \cup \{(0,1,2,3)\}$ 
\begin{align}
    \tem^{\, \tau \tau}_{3D} = \{\tem_0^{\, \tau \tau},\tem_1^{\, \tau \tau},\dots, \tem_{01}^{\, \tau \tau},\dots,\tem_{012}^{\, \tau \tau}, \dots,\tem_{0123}^{\, \tau \tau}\}  \, .
\end{align}
The element now reads
\begin{align}
    &\Reg^p(\Omega) = \left \{ \bigoplus_{i = 0}^3 \ver_i^p(\Omega) \otimes \tem_i^{\, \tau \tau} \right \} \oplus \left \{ \bigoplus_{j \in \mathcal{J}} \edge_j^p(\Omega) \otimes \tem_j^{\, \tau \tau} \right \} \oplus \left \{ \bigoplus_{k \in \mathcal{K}} \face_k^p(\Omega) \otimes \tem_k^{\, \tau \tau} \right \} \oplus \cell_{0123}^p(\Omega) \otimes \tem_{0123}^{\, \tau \tau} \, ,     
\end{align}
with the dimension
\begin{align}
    && \dim \Reg^p(\Omega) = \dim[\Po^p(\Omega) \otimes \Sym(3)] = (p+3)(p+2)(p+1) \, ,
\end{align}
assuming a complete polynomial space $\Po^p(\Omega)$.
\begin{definition}[Tetrahedral Regge base functions]
    The base functions of the tetrahedral Regge element are defined per polytope.
    \begin{itemize}
        \item On each edge $e_{ij}$ with multi-index $(i,j) \in \mathcal{J}$, equipped with the tangent vector $\bm{\tau}$ and trace operator $\tr_{\tau \tau} (\cdot) = \con{\bm{\tau} \otimes \bm{\tau}}{\cdot} |_{\mu_{ij}}$ we define the edge base functions
        \begin{subequations}
            \begin{align}
            &\text{Vertex-edge}: & \bm{\Upsilon}_\alpha(\xi,\eta) &= n \tt \, , & n &\in \ver^p_\alpha(\Gamma) \, , & \tt &\in \{ \tt \in \tem_\alpha^{\, \tau \tau} \; | \; \tr_{\tau \tau} \tt \neq 0 \} \, , \\
            &\text{Edge}: & \bm{\Upsilon}_{l}(\xi,\eta) &= n_l \tt \, , & n_l &\in \edge^p_{ij}(\Gamma) \, , & \tt &\in \{ \tt \in \tem_{ij}^{\, \tau \tau} \; | \; \tr_{\tau \tau} \tt \neq 0 \} \, , 
        \end{align}
        \end{subequations}
        where $\alpha \in \{i,j\}$. For each edge there are $2\cdot 1$ vertex-edge base functions and $p-1$ edge base functions.
        \item On each face $f_{ijk}$ with multi-index $(i,j,k) \in \mathcal{K}$, equipped with the unit normal vector $\bm{\nu}$ and the tangential projection operator $\tr_{\tau \tau} (\cdot) = (\one - \bm{\nu} \otimes \bm{\nu}) (\cdot)(\one - \bm{\nu} \otimes \bm{\nu}) |_{\Gamma_{ijk}}$ we define the face base functions
        \begin{subequations}
            \begin{align}
            &\text{Vertex-face}: & \bm{\Upsilon}_\alpha(\xi,\eta) &= n \tt \, , & n &\in \ver^p_\alpha(\Gamma) \, , & \tt &\in \{ \tt \in \tem_\alpha^{\, \tau \tau} \; | \; \tr_{\tau \tau} \tt \neq 0 \} \, , \\
            &\text{Edge-face}: & \bm{\Upsilon}_{\beta l q}(\xi,\eta) &= n_l \tt_q \, , & n_l &\in \edge^p_\beta(\Gamma) \, , & \tt_q &\in \{ \tt \in \tem_\beta^{\, \tau \tau} \; | \; \tr_{\tau \tau} \tt  \neq 0 \} \, , \\
            &\text{Face}: & \bm{\Upsilon}_{lq}(\xi,\eta) &= n_l \tt_q \, , & n_l &\in \face^p_{ijk}(\Gamma) \, , & \tt_q &\in \{ \tt \in \tem_{ijk}^{\, \tau \tau} \; | \; \tr_{\tau \tau} \tt \neq 0 \} \, , 
        \end{align}
        \end{subequations}
        where $\alpha \in \{i,j,k\}$, $\beta \in \mathcal{J}_{ijk} = \{(i,j),(i,k),(j,k)\} \subset \mathcal{J}$. For each face we find $3 \cdot 1$ vertex-face base functions, $3 \cdot 2 \cdot (p-1)$ edge-face base functions and $3 \cdot (p-2)(p-2)/2$ face base functions.
        \item The cell $c_{0123}$ is equipped with the tangential projection operator $\tr_{\tau \tau} (\cdot) = (\one - \bm{\nu} \otimes \bm{\nu}) (\cdot)(\one - \bm{\nu} \otimes \bm{\nu}) |_{\partial \Omega}$. Its base functions read
        \begin{subequations}
            \begin{align}
            &\text{Edge-cell}: & \bm{\Upsilon}_{\alpha l}(\xi,\eta) &= n_l \tt \, , & n_l &\in \edge^p_\alpha(\Gamma) \, , & \tt &\in \{ \tt \in \tem_\alpha^{\, \tau \tau} \; | \; \tr_{\tau \tau} \tt = 0 \} \, , \\
            &\text{Face-cell}: & \bm{\Upsilon}_{\beta l q}(\xi,\eta) &= n_l \tt_q \, , & n_l &\in \face^p_\beta(\Gamma) \, , & \tt_q &\in \{ \tt \in \tem_\beta^{\, \tau \tau} \; | \; \tr_{\tau \tau} \tt  = 0 \} \, , \\
            &\text{Cell}: & \bm{\Upsilon}_{lq}(\xi,\eta) &= n_l \tt_q \, , & n_l &\in \cell^p_{0123}(\Gamma) \, , & \tt_q &\in \tem_{0123}^{\, \tau \tau}  \, , 
        \end{align}
        \end{subequations}
        where $\alpha \in \mathcal{J}$ and $\beta \in \mathcal{K}$. There are $6 \cdot (p-1)$ edge-cell base functions, $4 \cdot 3 \cdot (p-2)(p-1)/2$ face-cell base functions and $6 \cdot (p-3)(p-2)(p-1)/6$ pure cell base functions. 
    \end{itemize}
\end{definition}
{\parindent0pt An example for the base functions of the linear Regge element is given in \cref{sec:linreg}.}

\section{Normal-continuous elements} \label{sec:norm}

\subsection{The Hellan--Herrmann--Johnson and Pechstein--Sch\"oberl elements}
The Hellan--Herrmann--Johnson element $\HHJ^p(\surf) \subsetsim \HdD{,\surf}$ is used to approximate symmetric normal-normal continuous fields on surfaces $\surf \subset \R^2$. The three-dimensional version of this construction is the Pechstein--Sch\"oberl element $\PS^p(\vol) \subsetsim \HdD{,\vol}$ with $\vol \subset \R^3$. For both elements the trace operator can be defined as
\begin{align}
    \tr_{nn} \bm{P} &= \con{\vb{n}\otimes\vb{n}}{\bm{P}} = P_{nn} \, , 
\end{align}
and a conforming construction must satisfy $\jump{\tr_{nn}\bm{P}}|_{\Xi} = 0$ on all arbitrarily defined interfaces $\Xi$ in the domain.
In the following we provide the definition of the base functions on the reference simplices. The double contravariant Piola transformation can then be employed to map the base functions to their counterparts on the physical mesh   
\begin{align}
    \con{\vb{n} \otimes \vb{n}}{\bm{Y}} = \con{(\cof\bm{J})\bm{\nu} \otimes \bm{\nu}(\cof\bm{J})^T}{\bm{Y}} = \con{\bm{\nu} \otimes \bm{\nu}}{\bm{\Upsilon}} \qquad \iff \qquad \bm{Y} = \dfrac{1}{(\det \bm{J})^2} \bm{J} \bm{\Upsilon}  \bm{J}^{T} \, , 
\end{align}
preserving its normal-normal projection properties.

\subsubsection{The Hellan--Herrmann--Johnson element} \label{sec:hhj}
We start with the Hellan--Herrmann--Johnson element \cite{Hel67,Her67,Joh73}. We note that the element can be directly constructed from the two-dimensional Regge element, since in two dimensions there hold the relations
\begin{align}
    &\vb{n} = \bm{R} \vb{t} \, , && \vb{n} \otimes \vb{n} = \bm{R} \vb{t} \otimes \vb{t}\bm{R}^T \, , && \bm{R} = \vb{e}_1 \otimes \vb{e}_2 - \vb{e}_2 \otimes \vb{e}_1 \, , 
\end{align}
implying that the Regge base functions can transformed to Hellan--Herrmann--Johnson base functions via the double rotation $ \HHJ^p(\Gamma) = \bm{R} [\Reg^p(\Gamma)] \bm{R}^T$. However, for the sake of completeness we present the full construction of the Hellan--Herrmann--Johnson element as well. \textit{\textbf{We note that the rotational property is not given for three-dimensional tetrahedra, such that Pechstein--Sch\"oberl elements cannot be defined by rotating tetrahedral Regge elements}}, $\PS^p(\Omega) \neq \bm{R} [\Reg^p(\Omega)] \bm{R}^T$.  

The construction of the Hellan--Herrmann--Johnson element follows analogously to the Regge triangle element. However, the template superset here is \cref{eq:norm2dtem}, for normal continuity
\begin{align}
    &\tem_\alpha^{\, \nu \nu} = \sym(\tem_\alpha^{\, \nu} \otimes \tem_\alpha^{\, \nu}) = \{ \tv_1 \otimes \tv_1,\tv_2 \otimes \tv_2,\sym(\tv_1 \otimes \tv_2) \} \, , && \dim  \tem_\alpha^{\, \nu \nu} = \dim \Sym(2) = 3 \, ,
    \label{eq:temhhj}
\end{align}
where $\tem_\alpha^{\, \nu} \in \tem^{\, \nu}_{2D}$ with the multi-index $\alpha \in \{0,1,2\} \cup \mathcal{J} \cup \{(0,1,2)\}$. The resulting superset $\tem^{\, \nu \nu}_{2D}$ is defined to contain all the resulting polytopal template subsets $\tem_\alpha^{\, \nu \nu}$
\begin{align}
    \tem^{\, \nu \nu}_{2D} = \{\tem_0^{\, \nu \nu},\tem_1^{\, \nu \nu},\dots, \tem_{01}^{\, \nu \nu},\dots,\tem_{012}^{\, \nu \nu}\}  \, .
\end{align}
The Hellan--Herrmann--Johnson element on the reference triangle is now given via the tensorial construction 
\begin{align}
    &\HHJ^p(\Gamma) = \left \{ \bigoplus_{i = 0}^2 \ver_i^p(\Gamma) \otimes \tem_i^{\, \nu \nu} \right \} \oplus \left \{ \bigoplus_{j \in \mathcal{J}} \edge_j^p(\Gamma) \otimes \tem_j^{\, \nu \nu} \right \} \oplus \cell_{012}^p(\Gamma) \otimes \tem_{012}^{\, \nu \nu} \, ,     
\end{align}
with the dimension
\begin{align}
    && \dim \HHJ^p(\Gamma) = \dim[\Po^p(\Gamma) \otimes \Sym(2)] = \dfrac{3(p+2)(p+1)}{2} \, ,
\end{align}
under the assumption of a complete polynomial space $\Po^p(\Gamma)$.
\begin{definition}[Hellan--Herrmann--Johnson base functions]
    The base functions of the Hellan--Herrmann--Johnson triangle element are defined per polytope.
    \begin{itemize}
        \item On each edge $e_{ij}$ with multi-index $(i,j) \in \mathcal{J}$, equipped with the normal vector $\bm{\nu} = \bm{R} \bm{\tau}$ and trace operator $\tr_{\nu \nu} (\cdot) = \con{\bm{\nu} \otimes \bm{\nu}}{\cdot} |_{\mu_{ij}}$ we define the edge base functions
        \begin{subequations}
            \begin{align}
            &\text{Vertex-edge}: & \bm{\Upsilon}_\alpha(\xi,\eta) &= n \tt \, , & n &\in \ver^p_\alpha(\Gamma) \, , & \tt &\in \{ \tt \in \tem_\alpha^{\, \nu \nu} \; | \; \tr_{\nu \nu} \tt \neq 0 \} \, , \\
            &\text{Edge}: & \bm{\Upsilon}_{l}(\xi,\eta) &= n_l \tt \, , & n_l &\in \edge^p_{ij}(\Gamma) \, , & \tt &\in \{ \tt \in \tem_{ij}^{\, \nu \nu} \; | \; \tr_{\nu \nu} \tt \neq 0 \} \, , 
        \end{align}
        \end{subequations}
        where $\alpha \in \{i,j\}$. For each edge there are $2 \cdot 1$ vertex-edge base functions and $p-1$ edge base functions.
        \item The cell $c_{012}$ is equipped with the trace operator $\tr_{\nu \nu} (\cdot) = \con{\bm{\nu} \otimes \bm{\nu}}{\cdot} |_{\partial \Gamma}$. Its base functions read
        \begin{subequations}
        \begin{align}
            &\text{Vertex-cell}: & \bm{\Upsilon}_\alpha(\xi,\eta) &= n \tt \, , & n &\in \ver^p_\alpha(\Gamma) \, , & \tt &\in \{ \tt \in \tem_\alpha^{\, \nu \nu} \; | \; \tr_{\nu \nu} \tt = 0 \} \, , \\
            &\text{Edge-cell}: & \bm{\Upsilon}_{\beta l  q}(\xi,\eta) &= n_l \tt_q \, , & n_l &\in \edge^p_{\beta}(\Gamma) \, , & \tt_q &\in \{ \tt \in \tem_{\beta}^{\, \nu \nu} \; | \; \tr_{\nu \nu} \tt = 0 \} \, , \\
            &\text{Cell}: & \bm{\Upsilon}_{lq}(\xi,\eta) &= n_l \tt_q \, , & n_l &\in \cell^p_{012}(\Gamma) \, , & \tt_q &\in \tem_{012}^{\, \nu \nu} \, ,
        \end{align}
        \end{subequations}
        where $\alpha \in \{0,1,2\}$ and $\beta \in \mathcal{J}$. There are $3\cdot 1$ vertex-cell base functions, $3 \cdot  2 \cdot (p-1)$ edge-cell base functions and $3 \cdot (p-2)(p-1)/2$ pure cell base functions. 
    \end{itemize}
\end{definition} 
{\parindent0pt We provide an example for a quadratic Hellan--Herrmann--Johnson basis in \cref{sec:quadhhj}.}

\subsubsection{The Pechstein--Sch\"oberl element}
Next we consider the Pechstein--Sch\"oberl tetrahedral element \cite{Sin2009,PS2011}. The construction of the tensorial template sets is via
\begin{align}
    \tem_\alpha^{\, \nu \nu} &= \sym(\tem_\alpha^{\, \nu} \otimes \tem_\alpha^{\, \nu}) = \{ \tv_1 \otimes \tv_1,\tv_2 \otimes \tv_2,\tv_3 \otimes \tv_3, \sym(\tv_1 \otimes \tv_2),\sym(\tv_1 \otimes \tv_3),\sym(\tv_2 \otimes \tv_3) \} \, ,  \notag \\ \dim  \tem_\alpha^{\, \nu \nu} &= \dim \Sym(3) = 6 \, , 
\end{align}
where the normal-continuous template subsets are from \cref{eq:tem3Dnormal}. The corresponding superset for the normal-normal continuous template tensors is given by
\begin{align}
    \tem^{\, \nu \nu}_{3D} = \{\tem_0^{\, \nu \nu},\tem_1^{\, \nu \nu},\dots, \tem_{01}^{\, \nu \nu},\dots,\tem_{012}^{\, \nu \nu}, \dots,\tem_{0123}^{\, \nu \nu}\}  \, .
\end{align}
The polynomial space of the Pechstein--Sch\"oberl element can now be spanned by
\begin{align}
    &\PS^p(\Omega) = \left \{ \bigoplus_{i = 0}^3 \ver_i^p(\Omega) \otimes \tem_i^{\, \nu \nu} \right \} \oplus \left \{ \bigoplus_{j \in \mathcal{J}} \edge_j^p(\Omega) \otimes \tem_j^{\, \nu \nu} \right \} \oplus \left \{ \bigoplus_{k \in \mathcal{K}} \face_k^p(\Omega) \otimes \tem_k^{\, \nu \nu} \right \} \oplus \cell_{0123}^p(\Omega) \otimes \tem_{0123}^{\, \nu \nu} \, ,     
\end{align}
and has the dimension
\begin{align}
    && \dim \PS^p(\Omega) = \dim[\Po^p(\Omega) \otimes \Sym(3)] = (p+3)(p+2)(p+1) \, ,
\end{align}
for a full polynomial space $\Po^p(\Omega)$.
\begin{definition}[Tetrahedral Pechstein--Sch\"oberl base functions]
    The base functions of the tetrahedral Pechstein--Sch\"oberl element are defined per polytope.
    \begin{itemize}
        \item On each face $f_{ijk}$ with multi-index $(i,j,k) \in \mathcal{K}$, equipped with the unit normal vector $\bm{\nu}$ and the trace operator $\tr_{\nu \nu} (\cdot) = \con{\bm{\nu} \otimes \bm{\nu}}{\cdot} |_{\Gamma_{ijk}}$ we define the face base functions
        \begin{subequations}
            \begin{align}
            &\text{Vertex-face}: & \bm{\Upsilon}_\alpha(\xi,\eta) &= n \tt \, , & n &\in \ver^p_\alpha(\Gamma) \, , & \tt &\in \{ \tt \in \tem_\alpha^{\, \nu \nu} \; | \; \tr_{\nu \nu} \tt \neq 0 \} \, , \\
            &\text{Edge-face}: & \bm{\Upsilon}_{\beta l q}(\xi,\eta) &= n_l \tt_q \, , & n_l &\in \edge^p_\beta(\Gamma) \, , & \tt_q &\in \{ \tt \in \tem_\beta^{\, \nu \nu} \; | \; \tr_{\nu \nu} \tt  \neq 0 \} \, , \\
            &\text{Face}: & \bm{\Upsilon}_{lq}(\xi,\eta) &= n_l \tt_q \, , & n_l &\in \face^p_{ijk}(\Gamma) \, , & \tt_q &\in \{ \tt \in \tem_{ijk}^{\, \nu \nu} \; | \; \tr_{\nu \nu} \tt \neq 0 \} \, , 
        \end{align}
        \end{subequations}
        where $\alpha \in \{i,j,k\}$ and $\beta \in \mathcal{J}_{ijk} = \{(i,j),(i,k),(j,k)\} \subset \mathcal{J}$. For each face we find $3 \cdot 1$ vertex-face base functions, $3 \cdot 1 \cdot (p-1)$ edge-face base functions and $1 \cdot (p-1)(p-2)/2$ face base functions.
        \item The cell $c_{0123}$ is equipped with the trace operator $\tr_{\nu \nu} (\cdot) = \con{\bm{\nu} \otimes \bm{\nu}}{\cdot} |_{\partial \Omega}$. Its base functions read
        \begin{subequations}
            \begin{align}
            &\text{Vertex-cell}: & \bm{\Upsilon}_{\alpha q}(\xi,\eta) &= n \tt_q \, , & n &\in \ver^p_\alpha(\Gamma) \, , & \tt_q &\in \{ \tt \in \tem_\alpha^{\, \nu \nu} \; | \; \tr_{\nu \nu} \tt = 0 \} \, , \\
            &\text{Edge-cell}: & \bm{\Upsilon}_{\beta l q}(\xi,\eta) &= n_l \tt_q \, , & n_l &\in \edge^p_\beta(\Gamma) \, , & \tt_q &\in \{ \tt \in \tem_\beta^{\, \nu \nu} \; | \; \tr_{\nu \nu} \tt = 0 \} \, , \\
            &\text{Face-cell}: & \bm{\Upsilon}_{\gamma l q}(\xi,\eta) &= n_l \tt_q \, , & n_l &\in \face^p_\gamma(\Gamma) \, , & \tt_q &\in \{ \tt \in \tem_\gamma^{\, \nu \nu} \; | \; \tr_{\nu \nu}  \tt  = 0 \} \, , \\
            &\text{Cell}: & \bm{\Upsilon}_{lq}(\xi,\eta) &= n_l \tt_q \, , & n_l &\in \cell^p_{0123}(\Gamma) \, , & \tt_q &\in \tem_{0123} \, , 
        \end{align}
        \end{subequations}
        where $\alpha \in \{0,1,2,3\}$, $\beta \in \mathcal{J}$ and $\gamma \in \mathcal{K}$. There are $4 \cdot 3$ vertex-cell base functions, $6 \cdot 4 \cdot (p-1)$ edge-cell base functions, $4 \cdot 5 \cdot (p-2)(p-1)/2$ face-cell base functions and $6 \cdot (p-3)(p-2)(p-1)/6$ pure cell base functions. 
    \end{itemize}
\end{definition}

\subsection{The Hu--Zhang and Hu--Ma--Sun elements}
The Hu--Zhang \cite{hu_family_2014,hu_family_2015,hu_finite_2016} and Hu--Ma--Sun \cite{hu2023new} elements are symmetric normal-continuous elements, $\HZ^p(\surf) \subset \HsD{,\surf}$ and $\HMS^p(\vol) \subset  \HsD{,\vol}$. The trace operator is therefore the same as for $\HD{,\body}$, namely 
\begin{align}
    &\tr_n \bm{P} = \bm{P} \vb{n}  = P_{1n} \vb{e}_1 + P_{2n} \vb{e}_2 \, , && \tr_n \bm{P} = \bm{P} \vb{n}  = P_{1n} \vb{e}_1 + P_{2n} \vb{e}_2 + P_{3n} \vb{e}_3 \, ,  
\end{align}
in two- and three dimensions, respectively.
In other words, the trace operator is applied row-wise, and a conforming definition must satisfy $\jump{\tr_n \bm{P}}|_\Xi = 0$ for all interfaces. 

\subsubsection{The Hu--Zhang element} \label{sec:huzhang}
As shown in \cite{arnold_mixed_2002} a conforming construction for $\HsD{,\surf}$ inevitably requires $\C^0(\surf)$-continuity of the vertex-base functions. As such, the Hu--Zhang element enforces full continuity at the vertices, and breaks this requirement on the edges. In fact, on each edge the tangential-tangential $P_{tt}$ component of a field is allowed to jump. 
Consequently, the transformation of the Hu--Zhang element introduced in \cite{sky2023reissnermindlin} leaves the vertex and cell functions unchanged, and applies an edge-wise transformation of the edge base functions
\begin{align}
     &\bm{Y} = \bm{T}_j \bm{\Upsilon} \bm{T}^T_j \, , && \bm{T}_j = \dfrac{1}{\norm{\bm{\tau}}^2}(\vb{t} \otimes \bm{\tau} + \vb{n}\otimes \bm{\nu}) \, ,   && j \in \mathcal{J} \, ,
\end{align}
for each edge $e_j$ with tangent and normal vectors $\bm{\tau}$ and $\bm{\nu} = \bm{R}\bm{\tau} \perp \bm{\tau}$ in the reference domain $\Gamma$, and corresponding counterparts $\vb{t}$ and $\vb{n} = \bm{R} \vb{t} \perp \vb{t}$ in the physical domain $\surf$. Observe that 
\begin{align}
    &\bm{T}(\bm{\tau} \otimes \bm{\tau}) \bm{T}^T = \vb{t} \otimes \vb{t} \, ,  && \bm{T}\sym(\bm{\tau} \otimes \bm{\nu}) \bm{T}^T = \sym(\vb{t} \otimes \vb{n}) \, , &&
    \bm{T}(\bm{\nu} \otimes \bm{\nu}) \bm{T}^T = \vb{n} \otimes \vb{n} \, , 
\end{align}
such that the symmetry of the tensorial basis is left intact.

We define the template sets of the Hu--Zhang element
\begin{align}
    \tem_0^{\, \mathrm{sym}, \nu} &= \tem_1^{\, \mathrm{sym}, \nu} = \tem_2^{\, \mathrm{sym}, \nu} = \tem_{012}^{\, \mathrm{sym}, \nu} = \{\vb{e}_1 \otimes \vb{e}_1 , \sym(\vb{e}_1 \otimes \vb{e}_2), \vb{e}_2 \otimes \vb{e}_2 \} \, , \notag \\ \tem_{01}^{\, \mathrm{sym}, \nu} &= \{ \vb{e}_2 \otimes \vb{e}_2 , \sym (\vb{e}_1 \otimes \vb{e}_2), \vb{e}_1 \otimes \vb{e}_1 \} \, , \notag \\ \tem_{02}^{\, \mathrm{sym}, \nu} &= \{ \vb{e}_1 \otimes \vb{e}_1 , -\sym (\vb{e}_1 \otimes \vb{e}_2), \vb{e}_2 \otimes \vb{e}_2 \} \, , \notag \\ \tem_{12}^{\, \mathrm{sym}, \nu} &= \{ (\vb{e}_1 - \vb{e}_2) \otimes (\vb{e}_1 - \vb{e}_2) , \sym [(\vb{e}_2 - \vb{e}_1) \otimes (\vb{e}_1 + \vb{e}_2)], (\vb{e}_1 + \vb{e}_2) \otimes (\vb{e}_1 + \vb{e}_2) \} \, .
\end{align}
The vertex and cell sets are given by a simple Cartesian basis for $\Sym(2)$. On edges, the $\Sym(2)$-basis is constructed specifically from the tangent $\bm{\tau}$ and normal $\bm{\nu} = \bm{R}\bm{\tau}$ vectors of the edge, such that the normal-normal, symmetric normal-tangential, and tangential-tangential tensor bases are clearly determined. The superset is given by 
\begin{align}
    \tem_{2D}^{\, \mathrm{sym}, \nu} = \{ \tem_0^{\, \mathrm{sym}, \nu} , \tem_1^{\, \mathrm{sym}, \nu}, \dots , \tem_{01}^{\, \mathrm{sym}, \nu} ,\dots , \tem_{012}^{\, \mathrm{sym}, \nu} \} \, .
    \label{eq:temnorm2Dsym}
\end{align}
The element can now be constructed via 
\begin{align}
    &\HZ^p(\Gamma) = \left \{ \bigoplus_{i = 0}^2 \ver_i^p(\Gamma) \otimes \tem_i^{\, \mathrm{sym}, \nu} \right \} \oplus \left \{ \bigoplus_{j \in \mathcal{J}} \edge_j^p(\Gamma) \otimes \tem_j^{\, \mathrm{sym}, \nu} \right \} \oplus \cell_{012}^p(\Gamma) \otimes \tem_{012}^{\, \mathrm{sym}, \nu} \, ,     
\end{align}
and has the dimension
\begin{align}
    & \dim \HZ^p(\Gamma) = \dim[\Po^p(\Gamma) \otimes \Sym(2)] = \dfrac{3(p+2)(p+1)}{2} \, ,
\end{align}
for a full polynomial space $\Po^p(\Gamma)$.
\begin{definition}[Hu--Zhang base functions]
    The base functions of the Hu--Zhang triangle element are defined per polytope.
    \begin{itemize}
        \item On each vertex $v_i$ with $i \in \{0,1,2\}$ we define the vertex base functions
        \begin{align}
            &\text{Vertex}: & \bm{\Upsilon}_{\alpha q}(\xi,\eta) &= n \tt_q \, , & n &\in \ver^p_\alpha(\Gamma) \, , & \tt_q &\in \tem_\alpha^{\, \mathrm{sym}, \nu}  \, ,
        \end{align}
        with $\alpha \in \{0,1,2\}$. On each vertex there are $3 \cdot 1$ base functions.   
        \item On each edge $e_{ij}$ with multi-index $(i,j) \in \mathcal{J}$, equipped with the normal vector $\bm{\nu} = \bm{R} \bm{\tau}$ and trace operator $\tr_{\nu} (\cdot) =  (\cdot) \bm{\nu} |_{\mu_{ij}}$ we define the edge base functions
        \begin{align}
            &\text{Edge}: & \bm{\Upsilon}_{l}(\xi,\eta) &= n_l \tt_q \, , & n_l &\in \edge^p_{ij}(\Gamma) \, , & \tt_q &\in \{ \tt \in \tem_{ij}^{\, \mathrm{sym}, \nu} \; | \; \tr_{\nu} \tt \neq 0 \} \, . 
        \end{align}
        For each edge there are $2 \cdot (p-1)$ edge base functions.
        \item The cell $c_{012}$ is equipped with the trace operator $\tr_{\nu} (\cdot) = (\cdot) \bm{\nu}|_{\partial \Gamma}$. Its base functions read
        \begin{subequations}
        \begin{align}
            &\text{Edge-cell}: & \bm{\Upsilon}_{\beta l}(\xi,\eta) &= n_l \tt \, , & n_l &\in \edge^p_{\beta}(\Gamma) \, , & \tt &\in \{ \tt \in \tem_{\beta}^{\, \mathrm{sym}, \nu} \; | \; \tr_{\nu } \tt = 0 \} \, , \\
            &\text{Cell}: & \bm{\Upsilon}_{lq}(\xi,\eta) &= n_l \tt_q \, , & n_l &\in \cell^p_{012}(\Gamma) \, , & \tt_q &\in \tem_{012}^{\, \mathrm{sym}, \nu} \, ,
        \end{align}
        \end{subequations}
        where $\beta \in \mathcal{J}$. There are $3 \cdot  1 \cdot (p-1)$ edge-cell base functions and $3 \cdot (p-2)(p-1)/2$ pure cell base functions. 
    \end{itemize}
\end{definition} 

\subsubsection{The Hu--Ma--Sun element}
The Hu--Ma--Sun tetrahedral element \cite{hu2023new} is related to the tetrahedral Hu--Zhang element \cite{hu_family_2015,hu_finite_2016}. Namely, the Hu--Ma--Sun element relaxes the continuity assumptions of the Hu--Zhang element on edges, by lifting the symmetric tangential-normal components of the tensor living on the edges from edge-type to face-type connectivity. 
Like the Hu--Zhang element, it requires full $\C^0(\vol)$-continuity of the vertex-base functions. On both the reference and physical vertices and in the cell we employ a simple Cartesian basis 
\begin{align}
    \tem_0^{\, \mathrm{sym}, \nu} &= \tem_1^{\, \mathrm{sym}, \nu} = \tem_2^{\, \mathrm{sym}, \nu} = \tem_3^{\, \mathrm{sym}, \nu} = \tem_{0123}^{\, \mathrm{sym}, \nu}  \\
    & \qquad = \{\vb{e}_1 \otimes \vb{e}_1 , \sym(\vb{e}_1 \otimes \vb{e}_2), \sym(\vb{e}_1 \otimes \vb{e}_3), \vb{e}_2 \otimes \vb{e}_2, \sym(\vb{e}_2 \otimes \vb{e}_3),\vb{e}_3 \otimes \vb{e}_3 \}  \, , \notag \\
    \dim \tem_0^{\, \mathrm{sym}, \nu} &= \dim\tem_1^{\, \mathrm{sym}, \nu} = \dim\tem_2^{\, \mathrm{sym}, \nu} = \dim\tem_3^{\, \mathrm{sym}, \nu} = \dim\tem_{0123}^{\, \mathrm{sym}, \nu} = \dim\Sym(3) = 6 \, . \notag
\end{align}
On each physical edge with the tangent vector $\vb{t} = \bm{J} \bm{\tau}$ we define the additional vectors $\vb{d}_1$, and $\vb{d}_2$, spanning an orthogonal system. A unique orthogonal vector $\vb{d}_2 \perp \vb{t}$ can be found using the vector $\vb{t}$ with the algorithmic formula \cite{Ken}
\begin{align}
    \vb{d}_2 = \orth \vb{t} = \begin{bmatrix}
        (\sgn t_1) |t_3| \\
        (\sgn t_2) |t_3| \\
        - (\sgn t_3) |t_1| - (\sgn t_3) |t_2|
     \end{bmatrix} \perp \vb{t} \, , 
     \label{eq:orthop}
\end{align}
such that $\norm{\vb{t}} \leq \norm{\vb{d}_2} \leq \sqrt{2}  \norm{\vb{t}}$. 
The specialised signum function is defined as
\begin{align}
    \sgn(t) = \left \{ \begin{matrix}
        1 & \text{for} & t \geq 0 \\
        -1 & \text{for} & t < 0 
    \end{matrix} \right . \, .
\end{align}
The $\vb{d}_1$ vector is subsequently constructed via
\begin{align}
    \vb{d}_1 = \vb{d}_2 \times \vb{t} \, ,
    \label{eq:d1}
\end{align}
such that $\norm{\vb{d}_1} = \norm{\vb{d}_2} \norm{\vb{t}}$. 
The triplet $\{\vb{t},\vb{d}_1,\vb{d}_2\}$ are the tangent vector, and two vectors that are respectively orthogonal to the tangent vector and each other. For conformity to be upheld between interfaces, the $\vb{d}_i$ vectors must be defined using solely the information of the tangent vector, which is common to all interfacing elements on an edge. 
Analogously, we can define on each edge of the reference element 
\begin{align}
    &\bm{\delta}_2 = \orth \bm{\tau} \, , &&\bm{\delta}_1 = \bm{\delta}_2 \times \bm{\tau} \, ,
\end{align}
using the corresponding edge tangent vector $\bm{\tau}$.
As such, we can define the edge-wise transformation tensor 
\begin{align}
    &\bm{T}_j = \dfrac{1}{\norm{\bm{\tau}}^2} \vb{t} \otimes \bm{\tau} + \dfrac{1}{\norm{\bm{\delta}_1}^2} \vb{d}_1 \otimes \bm{\delta}_1 + \dfrac{1}{\norm{\bm{\delta}_2}^2} \vb{d}_2 \otimes \bm{\delta}_2 \, , && j \in \mathcal{J} \, ,
    \label{eq:transdelta}
\end{align}
such that the transformation preserves symmetry of the orthogonal triplet
\begin{align}
     &\bm{T}_j: \{\bm{\tau},\bm{\delta}_1,\bm{\delta}_2\} \to \{\vb{t},\vb{d}_1,\vb{d}_2\} \, , && \bm{T}_j[\sym(\{\bm{\tau},\bm{\delta}_1,\bm{\delta}_2\} \otimes \{\bm{\tau},\bm{\delta}_1,\bm{\delta}_2\})]\bm{T}_j^T = \sym(\{\vb{t},\vb{d}_1,\vb{d}_2\} \otimes \{\vb{t},\vb{d}_1,\vb{d}_2\}) \, . 
\end{align}
In order to construct the template set of each edge we use the latter construction and the corresponding edge-face template vectors $\tv_1\,\tv_2 \in \tem_j^{\, \nu}$ from \cref{eq:tem3Dnormal}
\begin{align}
    \tem_j^{\, \mathrm{sym}, \nu} = \{ \bm{\delta}_1 \otimes \bm{\delta}_1,\sym(\bm{\delta}_1 \otimes \bm{\delta}_2), \bm{\delta}_2 \otimes \bm{\delta}_2, \sym (\bm{\tau} \otimes \bm{\psi}_1), \sym (\bm{\tau} \otimes \bm{\psi}_2), \bm{\tau} \otimes \bm{\tau} \} \, .
    \label{eq:temEdgeHuMaSun}
\end{align}
The first three template tensors are edge-type, the next two are face-type, and the last template tensor is cell type.
Due to the mix in tensorial bases, the transformation is also mixed. The vector set $ \{\bm{\tau},\bm{\delta}_1,\bm{\delta}_2\} $ is transformed via $\bm{T}_j$ of each edge $e_j$ with $j \in \mathcal{J}$. The mixed components $\{ \sym (\bm{\tau} \otimes \bm{\psi}_1), \sym (\bm{\tau} \otimes \bm{\psi}_2) \}$ are transformed via 
\begin{align}
    \bm{Y} = \dfrac{1}{\det \bm{J}} \bm{J} \bm{\Upsilon} \bm{J}^T \, ,
\end{align}
being a combination of the mapping of tangents $\bm{J}$ and the contravariant Piola transformation $(\det \bm{J})^{-1}\bm{J}$. The transformation upholds symmetry, the uniqueness of the edge tangent vector for interfacing elements, and the normal projections on neighbouring faces. 

On each reference face $f_k$ we define the 
triplet $\{\bm{\nu},\bm{\tau}_1,\bm{\tau}_2 \}$ via \cref{eq:orthop} and $\bm{\tau}_1 = \bm{\nu} \times \bm{\tau}_2$. As such, we can define the face-wise transformation tensor 
\begin{align}
    &\bm{T}_k = \dfrac{1}{\norm{\bm{\nu}}^2} \vb{n} \otimes \bm{\nu} + \dfrac{1}{\norm{\bm{\tau}_1}^2} \vb{t}_1 \otimes \bm{\tau}_1 + \dfrac{1}{\norm{\bm{\tau}_2}^2} \vb{t}_2 \otimes \bm{\tau}_2 \, , && k \in \mathcal{K} \, ,
    \label{eq:transdeltanu}
\end{align}
and there holds
\begin{align}
     &\bm{T}_k:\{\bm{\nu},\bm{\tau}_1,\bm{\tau}_2\} \to \{\vb{n},\vb{t}_1,\vb{t}_2\} \, , &&  \bm{T}_k[\sym(\{\bm{\nu},\bm{\tau}_1,\bm{\tau}_2\} \otimes \{\bm{\nu},\bm{\tau}_1,\bm{\tau}_2\})]\bm{T}_k^T = \sym(\{\vb{n},\vb{t}_1,\vb{t}_2\} \otimes \{\vb{n},\vb{t}_1,\vb{t}_2\}) \, ,
     \label{eq:facetriplet}
\end{align}
such that the transformation preserves symmetry.
The tensorial basis of each face is now given by the set
\begin{align}
    \tem_k^{\, \mathrm{sym}, \nu} &= \sym(\{\bm{\nu},\bm{\tau}_1,\bm{\tau}_2\} \otimes \{\bm{\nu},\bm{\tau}_1,\bm{\tau}_2\}) \\ & \qquad = \{\bm{\nu} \otimes \bm{\nu},\sym(\bm{\nu} \otimes \bm{\tau}_1),\sym(\bm{\nu} \otimes \bm{\tau}_2),\bm{\tau}_1 \otimes \bm{\tau}_1,\sym (\bm{\tau}_1 \otimes \bm{\tau}_2),\bm{\tau}_2 \otimes \bm{\tau}_2\} \, . \notag
\end{align}
The first three template tensors are face type, and the last three are cell type.
The template superset for the Hu--Ma--Sun element is given by all the subsets 
\begin{align}
    \tem_{3D}^{\, \mathrm{sym}, \nu} =  \{ \tem_0^{\, \mathrm{sym}, \nu} , \tem_1^{\, \mathrm{sym}, \nu}, \dots , \tem_{01}^{\, \mathrm{sym}, \nu} ,\dots , \tem_{012}^{\, \mathrm{sym}, \nu}, \dots,  \tem_{0123}^{\, \mathrm{sym}, \nu}\} \, ,
\end{align}
and the Hu--Ma--Sun space is spanned by
\begin{align}
    \HMS^p(\Omega) &= \left \{ \bigoplus_{i = 0}^3 \ver_i^p(\Omega) \otimes \tem_i^{\, \mathrm{sym}, \nu} \right \} \oplus \left \{ \bigoplus_{j \in \mathcal{J}} \edge_j^p(\Omega) \otimes \tem_j^{\, \mathrm{sym}, \nu} \right \} \oplus \left \{ \bigoplus_{k \in \mathcal{K}} \face_k^p(\Omega) \otimes \tem_k^{\, \mathrm{sym}, \nu} \right \} \oplus \cell_{0123}^p(\Omega) \otimes \tem_{0123}^{\, \mathrm{sym}, \nu} \, ,  
\end{align}
with dimension
\begin{align}
    & \dim \HMS^p(\Omega) = \dim[\Po^p(\Omega) \otimes \Sym(3)] = (p+3)(p+2)(p+1) \, ,
\end{align}
assuming a full polynomial space $\Po^p(\Omega)$.
\begin{definition}[Hu--Ma--Sun base functions]
    The base functions of the Hu--Ma--Sun tetrahedral element are defined polytope-wise.
    \begin{itemize}
        \item On each vertex $v_i$ with $i \in \{0,1,2,3\}$ we define the vertex base functions
        \begin{align}
            &\text{Vertex}: & \bm{\Upsilon}_{\alpha q}(\xi,\eta) &= n \tt_q \, , & n &\in \ver^p_\alpha(\Gamma) \, , & \tt_q &\in \tem_\alpha^{\, \mathrm{sym}, \nu}  \, ,
        \end{align}
        with $\alpha \in \{0,1,2,3\}$. On each vertex there are $6 \cdot 1$ base functions.   
        \item On each edge $e_j$ with multi-index $j \in \mathcal{J}$, equipped with the vectors $\bm{\delta}_1$ and $\bm{\delta}_2$, we define the edge base functions
        \begin{align}
            &\text{Edge}: & \bm{\Upsilon}_{l q}(\xi,\eta) &= n_l \tt_q \, , & n_l &\in \edge^p_{j}(\Gamma) \, , & \tt_q &\in \{ \bm{\delta}_1 \otimes \bm{\delta}_1,\sym(\bm{\delta}_1 \otimes \bm{\delta}_2), \bm{\delta}_2 \otimes \bm{\delta}_2 \} \subset \tem_j^{\, \mathrm{sym}, \nu} \, . 
        \end{align}
        For each edge there are $3 \cdot (p-1)$ edge base functions.
        \item On each face $f_{ijk}$ with multi-index $(i,j,k) \in \mathcal{K}$, equipped with the normal vector $\bm{\nu}$ and the trace operator $\tr_{\nu} (\cdot) = (\cdot) \bm{\nu} |_{\Gamma_{ijk}}$ we define the face base functions
        \begin{subequations}
            \begin{align}
            &\text{Edge-face}: & \bm{\Upsilon}_{\beta l}(\xi,\eta) &= n_l \tt \, , & n_l &\in \edge^p_\beta(\Gamma) \, , & \tt & \in \{ \tt \in \tem_\beta^{\, \mathrm{sym}, \nu} \setminus \tem_\beta^{\, \mathrm{sym}, \delta}  \; | \; \tr_{\nu} \tt  \neq 0 \} \, , \\
            &\text{Face}: & \bm{\Upsilon}_{lq}(\xi,\eta) &= n_l \tt_q \, , & n_l &\in \face^p_{ijk}(\Gamma) \, , & \tt_q &\in \{ \tt \in \tem_{ijk}^{\, \mathrm{sym}, \nu} \; | \; \tr_{\nu} \tt \neq 0 \} \, , 
        \end{align}
        \end{subequations}
        where $\alpha \in \{i,j,k\}$, $\beta \in \mathcal{J}_{ijk} = \{(i,j),(i,k),(j,k)\} \subset \mathcal{J}$ and $\tem_\beta^{\, \mathrm{sym}, \delta} = \{ \bm{\delta}_1 \otimes \bm{\delta}_1,\sym(\bm{\delta}_1 \otimes \bm{\delta}_2), \bm{\delta}_2 \otimes \bm{\delta}_2 \} \subset \tem_\beta^{\, \mathrm{sym}, \nu}$. For each face we find $3 \cdot 1 \cdot (p-1)$ edge-face base functions and $3 \cdot (p-1)(p-2)/2$ face base functions.
        \item The cell $c_{0123}$ is equipped with the trace operator $\tr_{\nu} (\cdot) = (\cdot) \bm{\nu}|_{\partial \Omega}$. Its base functions read
        \begin{subequations}
        \begin{align}
            &\text{Edge-cell}: & \bm{\Upsilon}_{\beta l}(\xi,\eta) &= n_l \tt \, , & n_l &\in \edge^p_{\beta}(\Gamma) \, , & \tt &\in \{ \tt \in \tem_{\beta}^{\, \mathrm{sym}, \nu} \; | \; \tr_{\nu } \tt = 0 \} \, , \\
            &\text{Face-cell}: & \bm{\Upsilon}_{\gamma l q}(\xi,\eta) &= n_l \tt_q \, , & n_l &\in \face^p_{\gamma}(\Gamma) \, , & \tt_q &\in \{ \tt \in 
            \tem_{\gamma}^{\, \mathrm{sym}, \nu} \; | \; \tr_{\nu } \tt = 0 \} \, , \\
            &\text{Cell}: & \bm{\Upsilon}_{lq}(\xi,\eta) &= n_l \tt_q \, , & n_l &\in \cell^p_{0123}(\Omega) \, , & \tt_q &\in \tem_{0123}^{\, \mathrm{sym}, \nu} \, ,
        \end{align}
        \end{subequations}
        where $\beta \in \mathcal{J}$ and $\gamma \in \mathcal{K}$. There are $6 \cdot  1 \cdot (p-1)$ edge-cell base functions, $4 \cdot 3 \cdot (p-2)(p-1)/2$ face-cell base functions and $6 \cdot (p-3)(p-2)(p-1)/6$ cell base functions. 
    \end{itemize}
\end{definition} 
We note that the sets for the edge-face basis functions can also be expressed as  
\begin{align}
    &\{ \tt \in \tem_\beta^{\, \mathrm{sym}, \nu} \setminus \tem_\beta^{\, \mathrm{sym}, \delta}  \; | \; \tr_{\nu} \tt  \neq 0 \} = \{ \tt \in \tem_\beta^{\, \mathrm{sym}, \nu}  \; | \; \tr_{\nu} \tt  \neq 0 \, , \quad \tr_{\tau} \tt  \neq 0 \} \, , && \tr_{\tau} \tt = \tt \bm{\tau} \, ,
\end{align}
making the identification of their respective components straight-forward. 

\section{Tangential-normal-continuous elements} \label{sec:tannorm}

\subsection{The Gopalakrishnan--Lederer--Sch\"oberl elements}

The Gopalakrishnan--Lederer--Sch\"oberl element \cite{Led2019,gopalakrishnan2020mass,GKLS2023} is defined to assert deviatoric tangential-normal-continuous continuity. To clarify, the space $\GLS^p(\body) \subsetsim \HcD{,\body}$ is constructed of element-wise deviatoric polynomial tensor fields $\Po^p \otimes \sl(d)$ with $d \in \{2,3\}$, such that their tangential-normal component is continuous across element interfaces. The corresponding trace operator is therefore
\begin{align}
    \tr_{t n} \bm{P} = \con{ \vb{t} \otimes \vb{n}}{\bm{P}} = P_{tn} \, ,
\end{align}
and a conforming construction must satisfy $\jump{\tr_{t n} \bm{P}}|_{\Xi} = 0$ for any arbitrary interface $\Xi$. Base functions on the reference element can be consistently mapped to the physical element via the mixed covariant-contravariant Piola transformation
\begin{align}
    \con{\vb{t} \otimes \vb{n}}{\bm{Y}} = \con{\bm{J}\bm{\tau} \otimes \bm{\nu}(\cof\bm{J})^T}{\bm{Y}} = \con{\bm{\tau} \otimes \bm{\nu}}{\bm{\Upsilon}} \qquad \iff \qquad \bm{Y} = \dfrac{1}{\det \bm{J}} \bm{J}^{-T} \bm{\Upsilon}  \bm{J}^{T} \, .
\end{align}
Observe that the transformation maintains the tracelessness of deviatoric tensors $\tr\bm{D} = 0$ since 
\begin{align}
    \dfrac{1}{\det \bm{J}}\tr(\bm{J}^{-T}\bm{D}\bm{J}^T) = \dfrac{1}{\det \bm{J}} \con{\bm{J}^{-T}\bm{D}\bm{J}^T}{\one} = \dfrac{1}{\det \bm{J}} \con{\bm{D}}{\bm{J}^{-1}\one\bm{J}} = \dfrac{1}{\det \bm{J}} \con{\bm{D}}{\one} = \dfrac{1}{\det \bm{J}} \tr\bm{D} = 0 \, . 
\end{align}

\subsubsection{Gopalakrishnan--Lederer--Sch\"oberl triangles}

In order to construct Gopalakrishnan--Lederer--Sch\"oberl triangle elements $\GLS^p(\surf) \subsetsim \HrD{,\surf}$, we make use of both the tangential template set \cref{eq:tem2Dt} and normal-template set \cref{eq:norm2dtem}. The set can now be constructed via
\begin{align}
    \tem_{\alpha}^{\, \mathrm{dev},\tau \nu} = \dev(\tem_\alpha^{\, \tau} \otimes \tem_\alpha^{\, \nu}) = \{ \tv_1^\tau \otimes \tv_1^\nu , \tv_1^\nu \otimes \tv_1^\tau,\dev(\tv_1^\tau \otimes \tv_2^\nu - \tv_2^\tau \otimes \tv_1^\nu) \} \, , && \dim  \tem_\alpha^{\, \mathrm{dev},\tau \nu} = \dim \sl(2) = 3 \, .
\end{align}
Note that the negative sign of the linear combination in $\dev(\tv_1^\tau \otimes \tv_2^\nu - \tv_2^\tau \otimes \tv_1^\nu)$ is not arbitrary but rather by choice. In fact, we have by construction $\tr(\tv_1^\tau \otimes \tv_2^\nu) = \tr(\tv_2^\tau \otimes \tv_1^\nu)$ such that for the choice of a negative sign we automatically get a deviatoric tensor and the operator $\dev$ is redundant. The only exception is the interior of the cell, where the Cartesian basis is employed. Further, the first two vectors in the tangential and normal sets are orthogonal to each other $\tv_1^\tau \perp \tv_1^\nu$ such that they already represent the deviatoric off-diagonal term of a tensor. The superset for the tangential-normal continuous template tensors can now be given by
\begin{align}
    \tem^{\, \mathrm{dev} , \tau \nu}_{2D} = \{\tem_0^{\, \mathrm{dev} , \tau \nu},\tem_1^{\, \mathrm{dev} , \tau \nu},\dots, \tem_{01}^{\, \mathrm{dev} , \tau \nu},\dots,\tem_{012}^{\, \mathrm{dev} , \tau \nu}\}  \, .
\end{align}
The polynomial space of the Gopalakrishnan--Lederer--Sch\"oberl element can now be spanned by
\begin{align}
    &\GLS^p(\Gamma) = \left \{ \bigoplus_{i = 0}^2 \ver_i^p(\Gamma) \otimes \tem_i^{\, \mathrm{dev} , \tau \nu} \right \} \oplus \left \{ \bigoplus_{j \in \mathcal{J}} \edge_j^p(\Gamma) \otimes \tem_j^{\, \mathrm{dev} , \tau \nu} \right \} \oplus \cell_{012}^p(\Gamma) \otimes \tem_{012}^{\, \mathrm{dev} , \tau \nu} \, ,     
\end{align}
and has the dimension
\begin{align}
    && \dim \GLS^p(\Gamma) = \dim[\Po^p(\Gamma) \otimes \sl(2)] = \dfrac{3(p+2)(p+1)}{2} \, ,
\end{align}
for a full polynomial space $\Po^p(\Gamma)$.
\begin{definition}[Gopalakrishnan--Lederer--Sch\"oberl triangle base functions]
    The base functions of the Gopalakrishnan--Lederer--Sch\"oberl triangle element are defined per polytope.
    \begin{itemize}
        \item On each edge $e_{ij}$ with multi-index $(i,j) \in \mathcal{J}$, equipped with the tangent $\bm{\tau}$ and normal vectors $\bm{\nu} = \bm{R} \bm{\tau}$, and trace operator $\tr_{\tau \nu} (\cdot) = \con{\bm{\tau} \otimes \bm{\nu}}{\cdot} |_{\mu_{ij}}$ we define the edge base functions
        \begin{subequations}
            \begin{align}
            &\text{Vertex-edge}: & \bm{\Upsilon}_\alpha(\xi,\eta) &= n \tt \, , & n &\in \ver^p_\alpha(\Gamma) \, , & \tt &\in \{ \tt \in \tem_\alpha^{\, \mathrm{dev} , \tau \nu} \; | \; \tr_{\tau \nu} \tt \neq 0 \} \, , \\
            &\text{Edge}: & \bm{\Upsilon}_{l}(\xi,\eta) &= n_l \tt \, , & n_l &\in \edge^p_{ij}(\Gamma) \, , & \tt &\in \{ \tt \in \tem_{ij}^{\, \mathrm{dev} , \tau \nu} \; | \; \tr_{\tau \nu} \tt \neq 0 \} \, , 
        \end{align}
        \end{subequations}
        where $\alpha \in \{i,j\}$. For each edge there are $2 \cdot 1$ vertex-edge base functions and $p-1$ edge base functions.
        \item The cell $c_{012}$ is equipped with the trace operator $\tr_{\tau \nu} (\cdot) = \con{\bm{\tau} \otimes \bm{\nu}}{\cdot} |_{\partial \Gamma}$. Its base functions read
        \begin{subequations}
        \begin{align}
            &\text{Vertex-cell}: & \bm{\Upsilon}_\alpha(\xi,\eta) &= n \tt \, , & n &\in \ver^p_\alpha(\Gamma) \, , & \tt &\in \{ \tt \in \tem_\alpha^{\, \mathrm{dev} , \tau \nu} \; | \; \tr_{\tau \nu} \tt = 0 \} \, , \\
            &\text{Edge-cell}: & \bm{\Upsilon}_{\beta l  q}(\xi,\eta) &= n_l \tt_q \, , & n_l &\in \edge^p_{\beta}(\Gamma) \, , & \tt_q &\in \{ \tt \in \tem_{\beta}^{\, \mathrm{dev} , \tau \nu} \; | \; \tr_{\tau \nu} \tt = 0 \} \, , \\
            &\text{Cell}: & \bm{\Upsilon}_{lq}(\xi,\eta) &= n_l \tt_q \, , & n_l &\in \cell^p_{012}(\Gamma) \, , & \tt_q &\in \tem_{012}^{\, \mathrm{dev} , \tau \nu} \, ,
        \end{align}
        \end{subequations}
        where $\alpha \in \{0,1,2\}$ and $\beta \in \mathcal{J}$. There are $3\cdot 1$ vertex-cell base functions, $3 \cdot  2 \cdot (p-1)$ edge-cell base functions and $3 \cdot (p-2)(p-1)/2$ pure cell base functions. 
    \end{itemize}
\end{definition} 

\subsubsection{Gopalakrishnan--Lederer--Sch\"oberl tetrahedra}
The construction of the base functions on tetrahedra $\GLS^p(\vol) \subsetsim \HcD{,\vol}$ is analogous to triangles. The tensor template sets are defined via
\begin{align}
    \tem_{\alpha}^{\, \mathrm{dev},\tau \nu} &= \dev(\tem_\alpha^{\, \tau} \otimes \tem_\alpha^{\, \nu}) = \{ \tv_1^\tau \otimes \tv_1^\nu ,\tv_2^\tau \otimes \tv_2^\nu ,\tv_3^\tau \otimes \tv_3^\nu , \tv_1^\nu \otimes \tv_1^\tau,\tv_2^\nu \otimes \tv_2^\tau,\tv_3^\nu \otimes \tv_3^\tau, \notag \\ & \qquad \qquad \qquad \qquad \qquad
    \dev(\tv_1^\tau \otimes \tv_2^\nu - \tv_2^\tau \otimes \tv_1^\nu),\dev(\tv_1^\tau \otimes \tv_3^\nu - \tv_1^\tau \otimes \tv_3^\nu) \} \, , \notag \\ \dim  \tem_\alpha^{\, \mathrm{dev},\tau \nu} &= \dim \sl(3) = 8 \, ,
\end{align}
where we made use of three-dimensional tangential \cref{eq:temtan3d} and normal \cref{eq:tem3Dnormal} vector sets. \textit{\textbf{Unlike in the two-dimensional case, the ordering of the three-dimensional vectors is not automatic (the indices $\{1,2,3\}$ are only for illustration), such that one must first check whether two vectors are parallel $\tv_i \parallel \tv_j$ or not. In general, orthogonal vectors build the off-diagonal tensors, and linear combinations of parallel vectors build the deviatoric diagonal terms}}.
The superset for the tangential-normal continuous template tensors can now be given by
\begin{align}
    \tem^{\, \mathrm{dev} , \tau \nu}_{3D} = \{\tem_0^{\, \mathrm{dev} , \tau \nu},\tem_1^{\, \mathrm{dev} , \tau \nu},\dots, \tem_{01}^{\, \mathrm{dev} , \tau \nu},\dots,\tem_{012}^{\, \mathrm{dev} , \tau \nu} , \dots , \tem_{0123}^{\, \mathrm{dev} , \tau \nu}\}  \, .
\end{align}
The polynomial space of the Gopalakrishnan--Lederer--Sch\"oberl tetrahedron is spanned by
\begin{align}
    &\GLS^p(\Omega) = \left \{ \bigoplus_{i = 0}^3 \ver_i^p(\Omega) \otimes \tem_i^{\, \mathrm{dev} , \tau \nu} \right \} \oplus \left \{ \bigoplus_{j \in \mathcal{J}} \edge_j^p(\Omega) \otimes \tem_j^{\, \mathrm{dev} , \tau \nu} \right \} \oplus \left \{ \bigoplus_{k \in \mathcal{K}} \face_k^p(\Omega) \otimes \tem_k^{\, \mathrm{dev} , \tau \nu} \right \} \oplus \cell_{0123}^p(\Omega) \otimes \tem_{0123}^{\, \mathrm{dev} , \tau \nu} \, ,     
\end{align}
and has the dimension
\begin{align}
    && \dim \GLS^p(\Omega) = \dim[\Po^p(\Omega) \otimes \sl(3)] = \dfrac{8(p+3)(p+2)(p+1)}{6} \, ,
\end{align}
for a full polynomial space $\Po^p(\Omega)$.
\begin{definition}[Tetrahedral Gopalakrishnan--Lederer--Sch\"oberl base functions]
    The base functions of the tetrahedral Gopalakrishnan--Lederer--Sch\"oberl element are defined per polytope.
    \begin{itemize}
        \item On each face $f_{ijk}$ with multi-index $(i,j,k) \in \mathcal{K}$, equipped with the unit normal vector $\bm{\nu}$ and the trace operator $\tr_{\tau \nu} (\cdot) = (\one - \bm{\nu} \otimes \bm{\nu})(\cdot)\bm{\nu} |_{\Gamma_{ijk}}$ we define the face base functions
        \begin{subequations}
            \begin{align}
            &\text{Vertex-face}: & \bm{\Upsilon}_{\alpha q}(\xi,\eta) &= n \tt_q \, , & n &\in \ver^p_\alpha(\Gamma) \, , & \tt_q &\in \{ \tt \in \tem_\alpha^{\, \nu \nu} \; | \; \tr_{\tau \nu} \tt \neq 0 \} \, , \\
            &\text{Edge-face}: & \bm{\Upsilon}_{\beta l q}(\xi,\eta) &= n_l \tt_q \, , & n_l &\in \edge^p_\beta(\Gamma) \, , & \tt_q &\in \{ \tt \in \tem_\beta^{\, \nu \nu} \; | \; \tr_{\tau \nu} \tt  \neq 0 \} \, , \\
            &\text{Face}: & \bm{\Upsilon}_{lq}(\xi,\eta) &= n_l \tt_q \, , & n_l &\in \face^p_{ijk}(\Gamma) \, , & \tt_q &\in \{ \tt \in \tem_{ijk}^{\, \nu \nu} \; | \; \tr_{\tau \nu} \tt \neq 0 \} \, , 
        \end{align}
        \end{subequations}
        where $\alpha \in \{i,j,k\}$ and $\beta \in \mathcal{J}_{ijk} = \{(i,j),(i,k),(j,k)\} \subset \mathcal{J}$. For each face we find $3 \cdot 2$ vertex-face base functions, $3 \cdot 2 \cdot (p-1)$ edge-face base functions and $2 \cdot (p-1)(p-2)/2$ face base functions.
        \item The cell $c_{0123}$ is equipped with the trace operator $\tr_{\tau \nu} (\cdot) = (\one - \bm{\nu} \otimes \bm{\nu})(\cdot)\bm{\nu} |_{\partial \Omega}$. Its base functions read
        \begin{subequations}
            \begin{align}
            &\text{Vertex-cell}: & \bm{\Upsilon}_{\alpha q}(\xi,\eta) &= n \tt_q \, , & n &\in \ver^p_\alpha(\Gamma) \, , & \tt_q &\in \{ \tt \in \tem_\alpha^{\, \nu \nu} \; | \; \tr_{\tau \nu} \tt = 0 \} \, , \\
            &\text{Edge-cell}: & \bm{\Upsilon}_{\beta l q}(\xi,\eta) &= n_l \tt_q \, , & n_l &\in \edge^p_\beta(\Gamma) \, , & \tt_q &\in \{ \tt \in \tem_\beta^{\, \nu \nu} \; | \; \tr_{\tau \nu} \tt = 0 \} \, , \\
            &\text{Face-cell}: & \bm{\Upsilon}_{\gamma l q}(\xi,\eta) &= n_l \tt_q \, , & n_l &\in \face^p_\gamma(\Gamma) \, , & \tt_q &\in \{ \tt \in \tem_\gamma^{\, \nu \nu} \; | \; \tr_{\tau \nu}  \tt  = 0 \} \, , \\
            &\text{Cell}: & \bm{\Upsilon}_{lq}(\xi,\eta) &= n_l \tt_q \, , & n_l &\in \cell^p_{0123}(\Gamma) \, , & \tt_q &\in \tem_{0123} \, , 
        \end{align}
        \end{subequations}
        where $\alpha \in \{0,1,2,3\}$, $\beta \in \mathcal{J}$ and $\gamma \in \mathcal{K}$. There are $4 \cdot 2$ vertex-cell base functions, $6 \cdot 4 \cdot (p-1)$ edge-cell base functions, $4 \cdot 6 \cdot (p-2)(p-1)/2$ face-cell base functions and $8 \cdot (p-3)(p-2)(p-1)/6$ pure cell base functions. 
    \end{itemize}
\end{definition}

\section{An application to the Reissner--Mindlin plate problem} \label{sec:appRM}
This work introduces a novel methodology for the construction of numerous tensorial finite elements, \textit{\textbf{focusing on the underlying theory}}. However, in order to demonstrate the usefulness of the methodology and the discussed elements we present here a practical application to the Reissner--Mindlin plate problem.

The primal formulation of the Reissner--Mindlin plate reads
\begin{align} 
    &\int_\surf \dfrac{t^3}{12} \langle \sym \D \delta \bm{\phi} , \, \mathbb{D} \sym \D \bm{\phi} \rangle + \ksmu \, t \langle \nabla \delta w - \delta \bm{\phi} , \, \nabla w - \bm{\phi} \rangle \, \dd \surf  =\int_\surf t \, \delta w \, f \, \dd \surf \, , \quad \forall\, \{\delta w, \delta \bm{\phi}\} \in \V^p(\surf) \, ,
\end{align}   
where the deflection $w \in \Hone(\surf)$ and the rotations $\bm{\phi} \in [\Hone(\surf)]^2$ are to be determined, and $\V^p(\surf) = \CG^p(\surf) \times [\CG^p(\surf)]^2 \subset \Hone(\surf) \times [\Hone(\surf)]^2$. The material constants of the problem are Young's modulus $E$ and the Poisson ratio $\nu$ within the plane-stress elasticity tensor $\mathbb{D} = E/(1- \nu^2)[\nu \one \otimes \one + (1-\nu) \mathbb{J}]$, as well as the weighted shear modulus $\ksmu$, and the plate-thickness $t$.
Unfortunately, this formulation is susceptible to the problem of shear-locking, which led to countless methods of alleviating it \cite{bathe_mitc7_1989,hale_simple_2018,falk_finite_2008,sky2023reissnermindlin,ORYNYAK2024104103,HuBordas2020,pechstein_tdnns_2017,daVeiga}.
Concisely, the problem of shear locking is the inability of the discrete rotation field $\bm{\phi} \in [\CG^p(\surf)]^2$ to satisfactorily approach the gradient of the discrete deflection field $\nabla w \in \nabla \CG^p(\surf) \subset \Nedtwo^{p-1}(\surf)$ in the case of a very thin plate $t \to 0$. In this scenario, the Kirchhoff--Love constraint $\bm{\phi} = \nabla w$ must be satisfied in order to avoid a shear-dominated kinematic where bending is expected. Correspondingly, locking occurs due to the incompatibility of the discrete spaces $[\CG^p(\surf)]^2 \nsubseteq \Nedtwo^{p-1}(\surf)$ and $[\CG^p(\surf)]^2 \nsupseteq \Nedtwo^{p-1}(\surf)$.

In \cite{sky2023reissnermindlin} we introduced a locking-free mixed four-field formulation for Reissner--Mindlin plate (FFSRM) using the Hu--Zhang elements (\cref{sec:huzhang})
\begin{align} 
        \int_\surf \con{\delta \bm{M}}{\A \bm{M}} + \dfrac{t^2}{\ksmu} \con{\delta \vb{q}}{ \vb{q}} + \con{\Di \delta \bm{M}}{\bm{\phi}} - (\di \delta \vb{q}) \, w  -  \con{\delta \vb{q}}{\bm{\phi}} \, \dd \surf &= 0 \quad \forall \, \{\delta \bm{M}, \delta \vb{q} \} \in \Z^p(\surf) \, ,  \\
        \int_\surf \con{\delta \bm{\phi}}{\Di \bm{M}} -  \delta  w \, (\di  \vb{q}) - \con{\delta \bm{\phi}}{\vb{q}} \, \dd \surf &= - \int_\surf \delta w \, g \, \dd \surf \qquad  \forall \, \{\delta w, \delta \bm{\phi} \} \in \mathit{D}^p(\surf) \, , \notag
\end{align} 
with the compliance tensor $\mathbb{A} = 12 \, \mathbb{D}^{-1} = (12/E) [(1+ \nu) \mathbb{J} - \nu \one \otimes \one ]$, the forces $g = t^{-2}f$, and the bending moments $\bm{M} \in \HsD{,\surf}$ and shear stress $\vb{q} \in \Hd{,\surf}$ as two additional unknowns, such that $\Z^p(\surf) = \HZ^p(\surf) \times \RT^{p-1}(\surf)$ and $\mathit{D}^p(\surf) = \DG^{p-1}(\surf) \times [\DG^{p-1}(\surf)]^2$. Another locking-free method is introduced in \cite{pechstein_tdnns_2017} using the TDNNS method via the Hellan--Herrmann--Johnson elements (\cref{sec:hhj})
\begin{align} 
        \int_\surf \con{\delta \bm{M}}{\mathbb{A} \bm{M}} \, \dd \surf +  \con{\Di \delta \bm{M}}{\bm{\phi}}_{\mathcal{T}}  &= 0 \quad \forall \, \delta \bm{M} \in \HHJ^{p-1}(\surf) \, , \\
      \con{\delta \bm{\phi}}{\Di \bm{M}}_{\mathcal{T}}  - \int_\surf \dfrac{\ksmu}{t^2} \con{\nabla \delta w - \delta \bm{\phi}}{\nabla w - \bm{\phi}} \, \dd \surf &= - \int_\surf \delta w \, g \, \dd \surf \quad \forall \, \{\delta w, \delta \bm{\phi}\} \in \CG^p(\surf) \times \Ned^{p-1}(\surf) \, , \notag
\end{align} 
where the scalar product $\con{\Di \delta \bm{M}}{\bm{\phi}}_{\mathcal{T}}$ is to be understood in the distributional sense \cite{pechstein_tdnns_2017} and includes boundary terms on each element $\con{\Di \delta \bm{M}}{\bm{\phi}}_{\mathcal{T}} = \sum_{T \in \mathcal{T}} \langle \Di \delta\bm{M} , \, \bm{\phi} \rangle_{\Le(T)} - \langle \delta \bm{M} \vb{n} , ( \vb{t}\otimes \vb{t} )\bm{\phi}  \rangle_{\Le(\partial T)}$. 
In the following we discuss the two methods over two examples.

\subsection{Numerical examples}

\subsubsection{Minimal regularity}
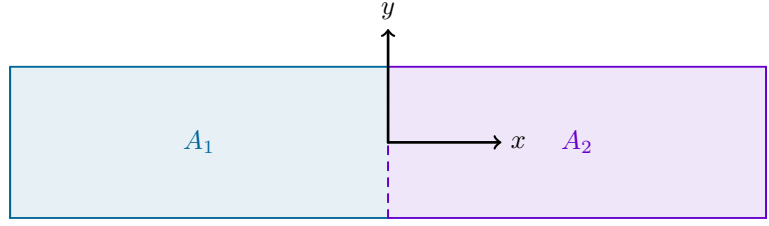
\begin{figure}
		\centering
		\input{figs/ex1dom}
		\caption{Domain composed of two materials $\surf = \surf_1 \cup \surf_2$. Zero deflection $\widetilde{w} = 0$ and zero bending movements $\widetilde{\bm{M}} = 0$ are imposed on the boundary of the domain $\partial \surf$.}
		\label{fig:ex1dom}
\end{figure}

\begin{figure}
		\centering
		\begin{subfigure}{0.48\linewidth}
			\centering
			\includegraphics[width=1\linewidth]{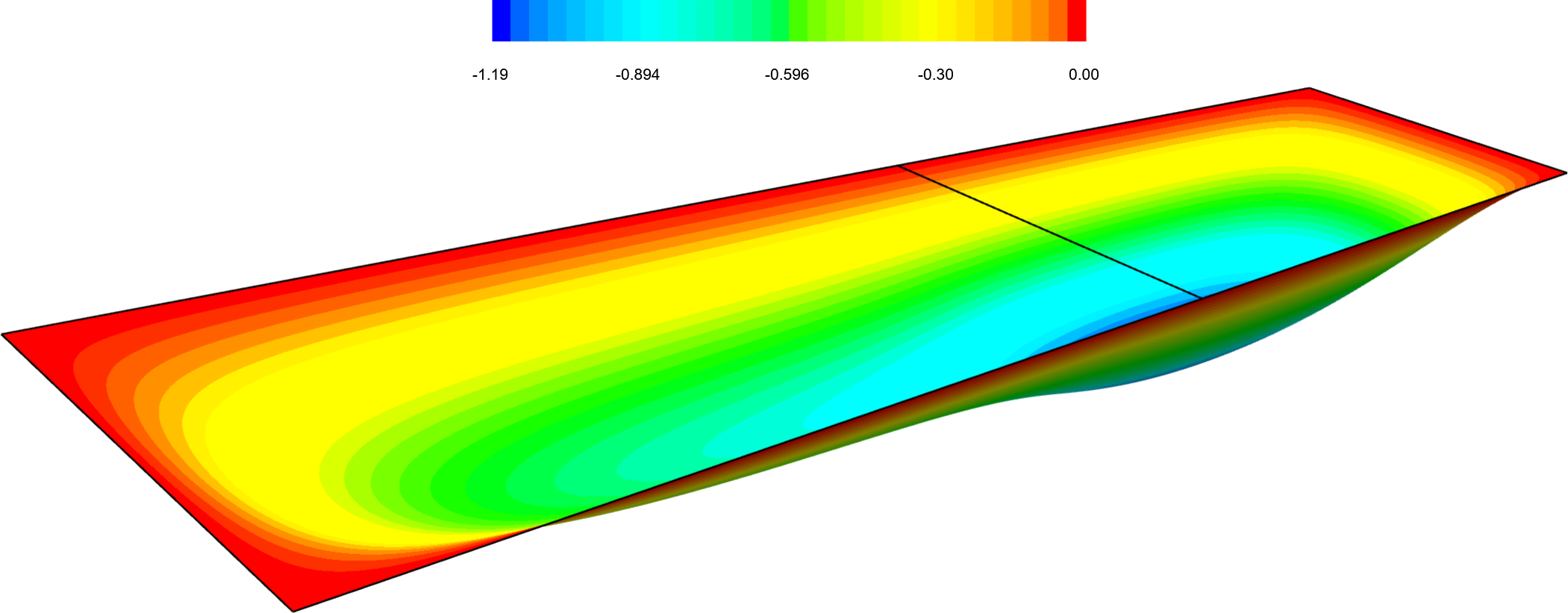}
			\caption{}
		\end{subfigure}
	    \begin{subfigure}{0.48\linewidth}
	    	\centering
	    	\includegraphics[width=1\linewidth]{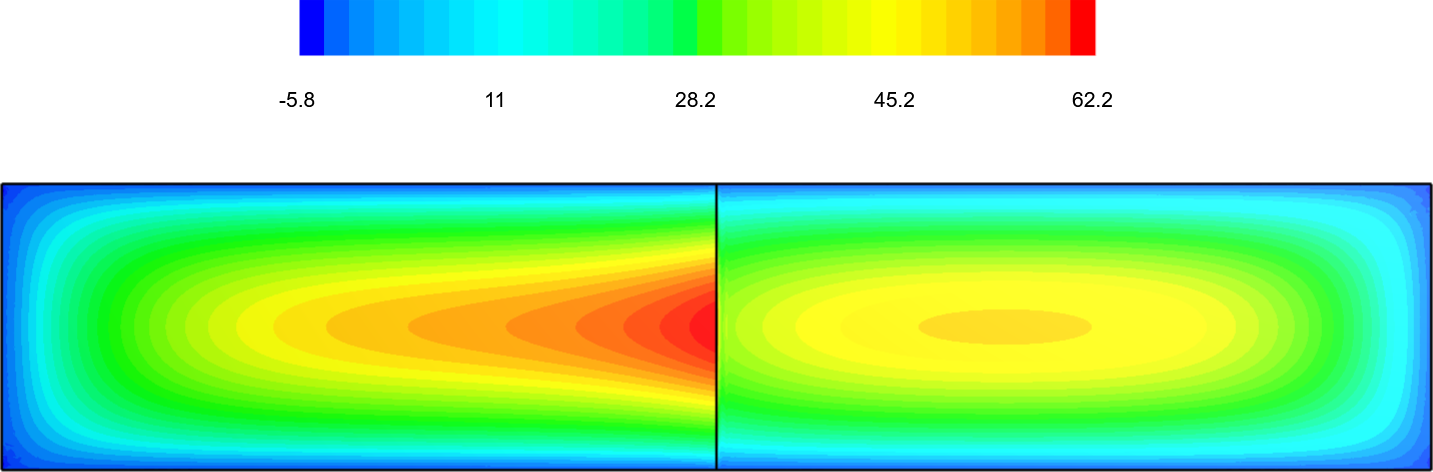}
	    	\caption{}
	    \end{subfigure}
		\caption{Solution for the deflection field $w$ (a) and the component of the bending moments $M_{yy}$ (b).}
		\label{fig:prmex}
\end{figure}

We consider the thin-plate domain $\overline{\surf} = [-5,5] \times [-1,1]$ with the constant thickness $t = 0.1$. 
We split the domain across two materials such that $\surf_1 = [-5,0) \times [-1, 1]$ and $\surf_2 = (0,5] \times [-1, 1]$ as depicted in \cref{fig:ex1dom}. In $\surf_1$ we define the material parameters $E = 300$ and $\nu = 0.5$, whereas the material parameters in $\surf_2$ are $E = 150$ and $\nu = 0.5$. On the boundary of the domain $\partial \surf$ we impose the boundary condition $\widetilde{w}|_{\partial \surf} = 0$, while leaving the rotation $\bm{\phi}$ free, implying a complete Neumann boundary for the rotations. The latter represents a Dirichlet boundary condition of zero bending moments on the boundary $\widetilde{\bm{M}}|_{\partial \surf} =  0$ in the mixed formulations. Finally, we apply the constant force $g = -100$.
We compute a high-fidelity solution using the primal formulation with a very fine mesh of $8282$ elements and polynomial order $p = 3$. The resulting deflection field $w$ and the $y-y$-component of the bending moments $M_{yy}$ are depicted in \cref{fig:prmex}.
The differing material parameters in the domain cause the deflection $w$ to tend towards $\surf_2$, which is softer. More importantly, it is directly observable that the sudden change in material coefficients induces a jump in the bending moment $M_{yy}$. Since in the primal formulation the bending moments are retrieved via the constitutive relation from the symmetrised gradient of the rotations, their discontinuity is correctly captured by the formulation. To clarify, the field $\bm{M} = \mathbb{D} \sym \D \bm{\phi}$ is inherently discontinuous for $\bm{\phi} \in [\CG^p(\surf)]^2 \subset [\C^0(\surf)]^2$, since the symmetrised gradient is only guaranteed to be tangential-tangential continuous, and the application of the material tensor $\mathbb{D}$ further mixes the tensorial components such that even tangential-tangential continuity is no longer given.   

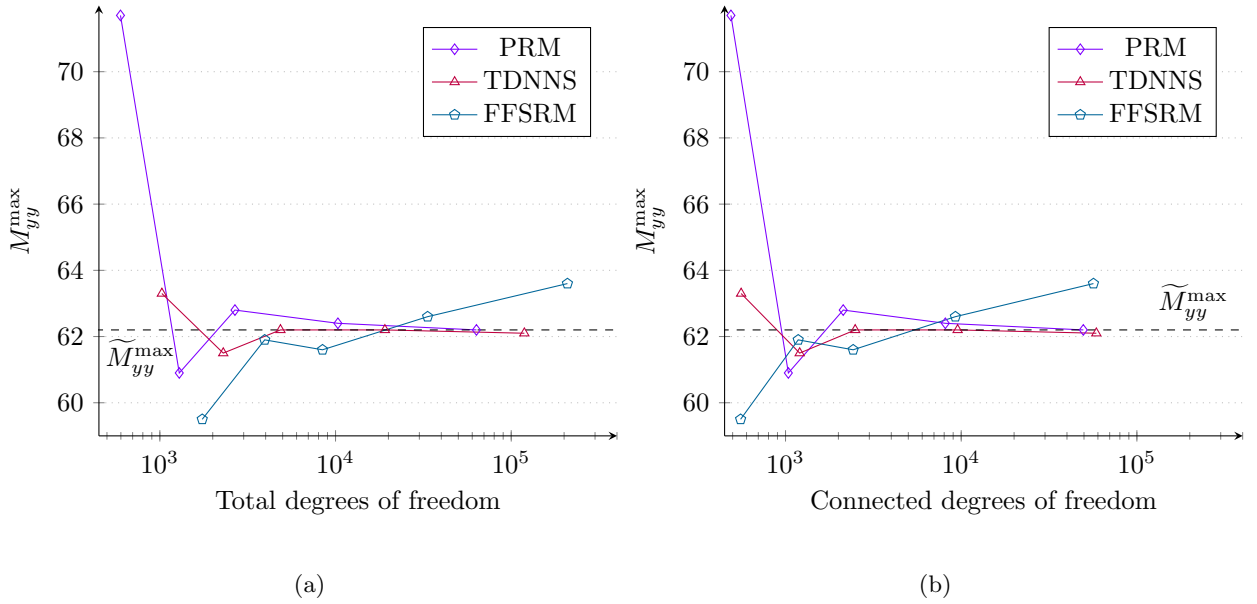
\begin{figure}
    	\centering
    	\begin{subfigure}{0.48\linewidth}
    		\centering
    		\input{figs/ex11}
    		\caption{}
    	\end{subfigure}
            \begin{subfigure}{0.48\linewidth}
    		\centering
    		\input{figs/ex12}
    		\caption{}
    	\end{subfigure}
    	\caption{Convergence towards the maximal bending moments in the plate $M_{yy}^\mathrm{max}$. In (a) the total number of degrees of freedom is given while in (b) only the connectivity-type degrees of freedom are considered, such that all cell-type degrees of freedom are eliminated. The graph in (a) represents the total computational effort whereas the graph in (b) considers only the effort to solve the final global system.}
    \label{fig:ex1res}
\end{figure}

\begin{figure}
		\centering
		\begin{subfigure}{0.48\linewidth}
			\centering
		\includegraphics[width=1\linewidth]{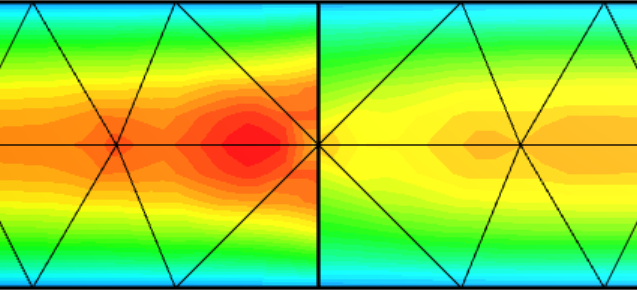}
        \includegraphics[width=1\linewidth]{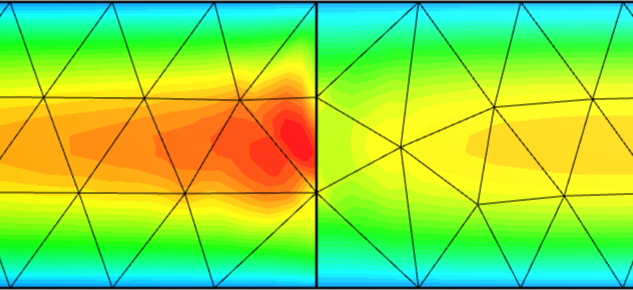}
        \includegraphics[width=1\linewidth]{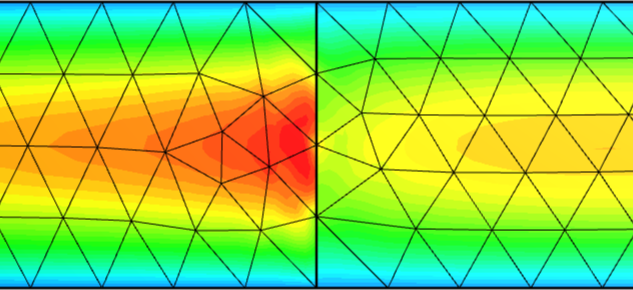}
			\caption{}
		\end{subfigure}
	    \begin{subfigure}{0.48\linewidth}
	    	\centering
	    \includegraphics[width=1\linewidth]{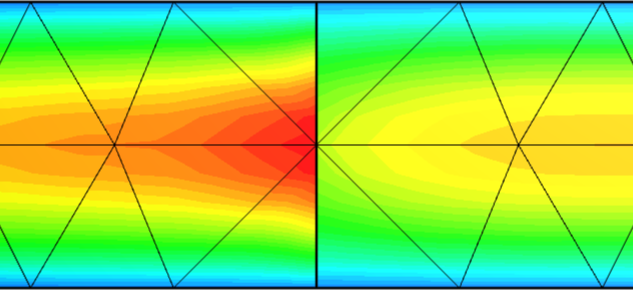}
        \includegraphics[width=1\linewidth]{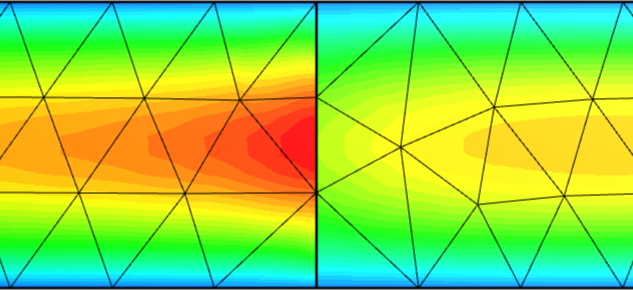}
        \includegraphics[width=1\linewidth]{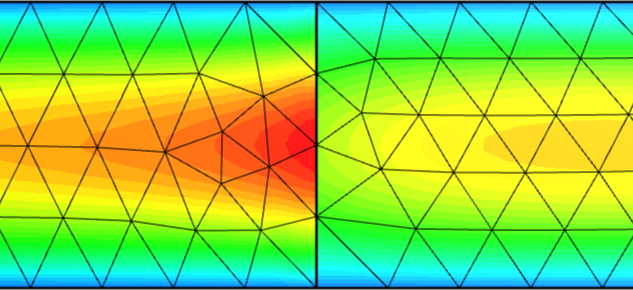}
	    	\caption{}
	    \end{subfigure}
		\caption{Solution of the bending moment $M_{yy}$ for FFSRM (a) and TDNNS (b) over meshes with $36$, $84$ and $182$ elements. The TDNNS formulation satisfies the jump condition completely whereas FFSRM must compensate for the increased continuity at the vertices on the interface.}
		\label{fig:myy}
\end{figure}

In the following we consider five refinements of the domain on unstructured meshes with $36$, $84$, $182$, $732$ and $4616$ elements, and compare all three formulations for the cubic polynomial order $p = 3$. The results are depicted in \cref{fig:ex1res}.
The high-fidelity solution of the maximal bending moment is $\widetilde{M}_{yy}^\mathrm{max} = 62.2$. We observe that the primal (PRM) and TDNNS formulations both converge towards the correct maximal value. In terms of the total amount of degrees of freedom, the primal formulation is in general cheaper, see \cref{fig:ex1res} (a). However, cell-type degrees of freedom are efficiently eliminated element-wise via static condensation on small dense matrices in the assembly procedure, such that the main computational effort is the solution of the global matrix. Consequently, a higher level of accuracy in the solution of the bending moments $\bm{M}$ at similar cost to the primal formulation can be achieved, as shown in \cref{fig:ex1res} (b) where the cell-type degrees of freedom are disregarded. Notably, the FFSRM formulation oscillates around the maximal value, such that the placement of vertices on the interface play a crucial role, see \cref{fig:myy}. In other words, the higher continuity at the vertices disrupts the computation from finding the maximal bending moment.  

\subsubsection{Curved mappings and boundary conditions}
The strong symmetry of the Hu--Zhang element in the FFSRM formulation implies that at the vertices, no simple separation of the tensorial components according to their projections on edges is generally possible. In other words, when enforcing a Dirichlet boundary condition, one automatically also controls the tangent-tangent component of the bending moments $M_{tt}$ at the vertices. This can only be avoided in special cases, where the tangent-tangent component is clear from the geometry of the domain, and the construction of the vertex base functions rotates the underlying Cartesian tensorial basis accordingly, compare \cite{sky2023reissnermindlin}. This increased control is problematic, since $M_{tt} = 0$ implies a hinge (free tangential rotation) where one might not be envisaged. Simply put, the lack of a clear association with edges may lead to the enforcement of incorrect boundary conditions.  

In order to simultaneously, demonstrate this phenomenon and the applicability of our mappings to curved geometries, we consider the curved L-shaped domain $\overline{\surf} = \{ (x,y) \in [1,-1]^1 \; | \; x^2 + y^2 \leq 1  \} \setminus (-1,1]^2$. We enforce the boundary conditions $\widetilde{w}|_{\partial \surf} = 0$ and $\tr_{*} \widetilde{M}|_{\partial \surf} = 0$, and set the material parameters to $t = 10^{-3}$, $E = 240$ and $\nu = 0.3$, while the force reads $g = -1000$. Note that $\tr_*(\cdot)$ is defined by the connectivity of the respective Hu--Zhang or Hellan--Herrmann--Johnson constructions, differing between the two. In particular, for the Hu--Zhang element the operator controls also vertex values.   
\begin{figure}
		\centering
		\begin{subfigure}{0.48\linewidth}
			\centering
		\includegraphics[width=0.7\linewidth]{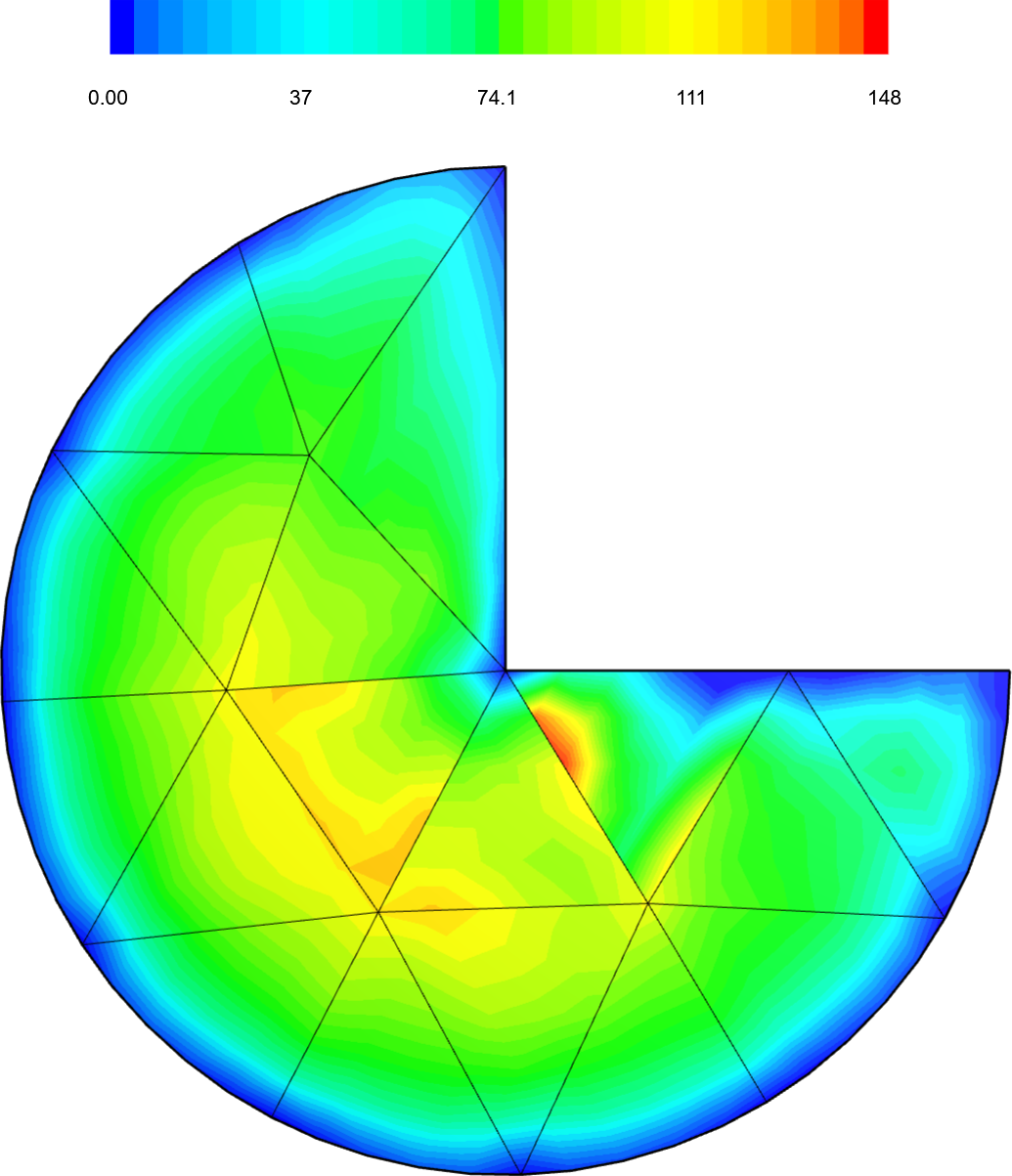}
        \includegraphics[width=0.7\linewidth]{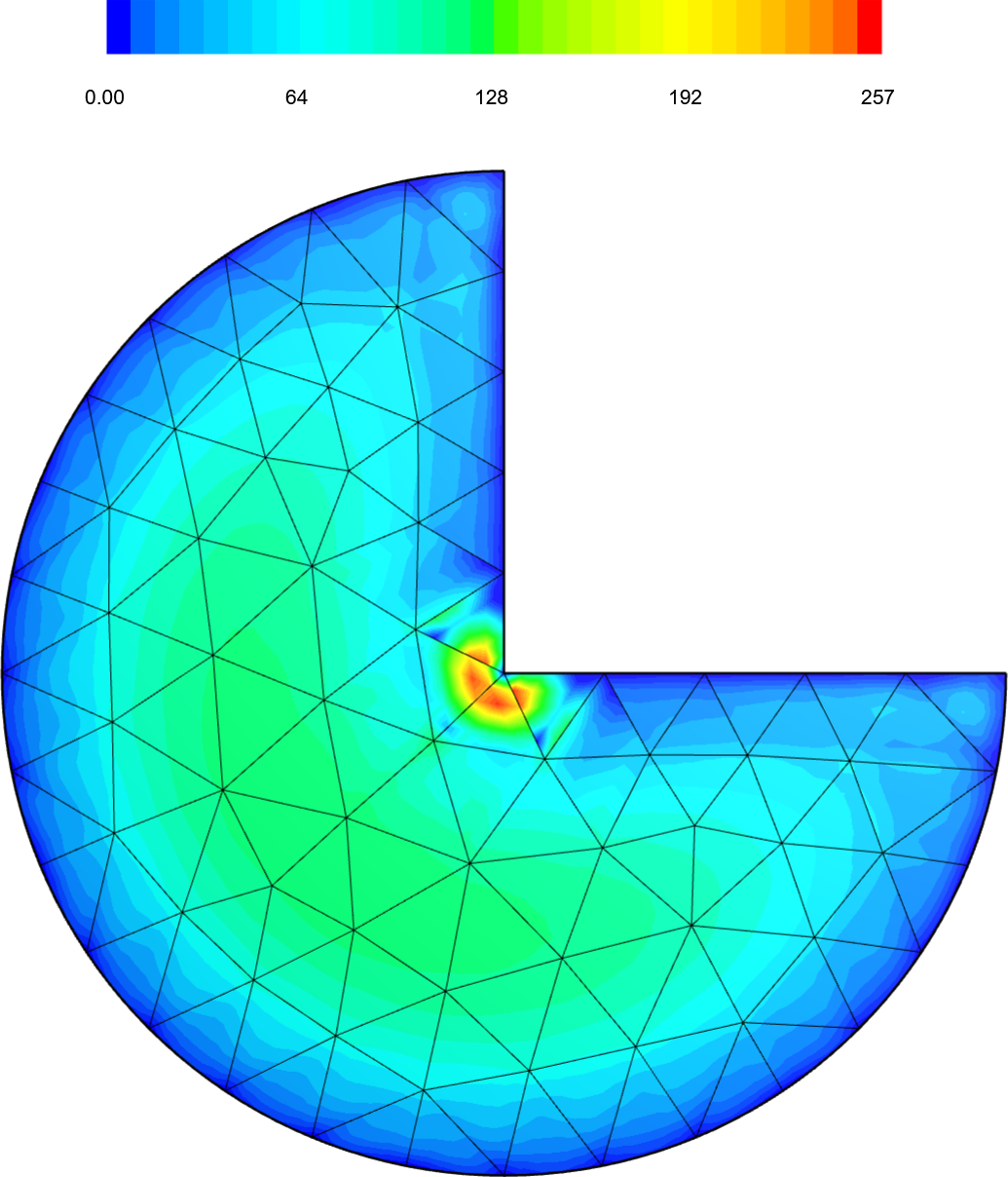}
        \includegraphics[width=0.7\linewidth]{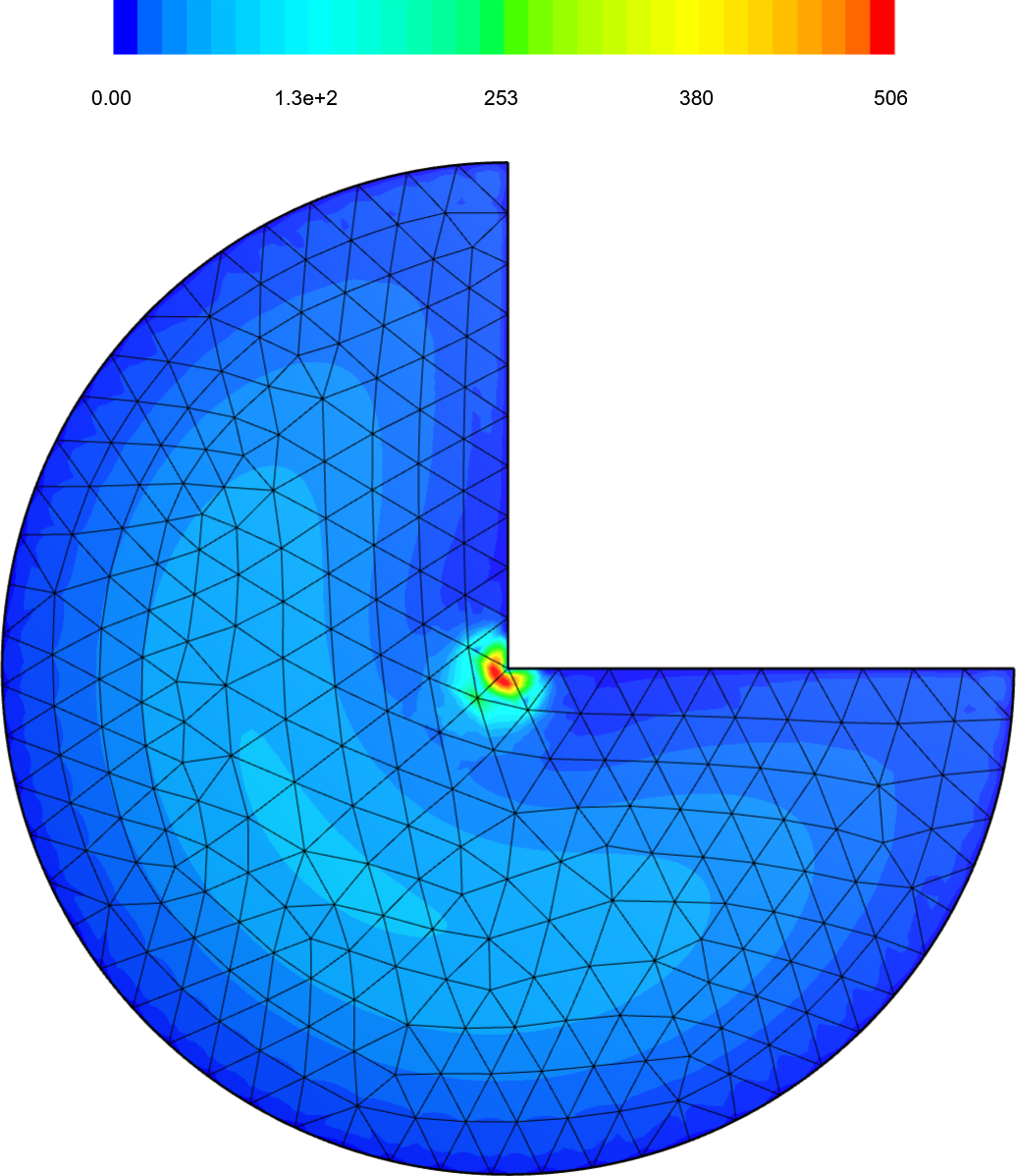}
			\caption{}
		\end{subfigure}
	    \begin{subfigure}{0.48\linewidth}
	    	\centering
	    \includegraphics[width=0.7\linewidth]{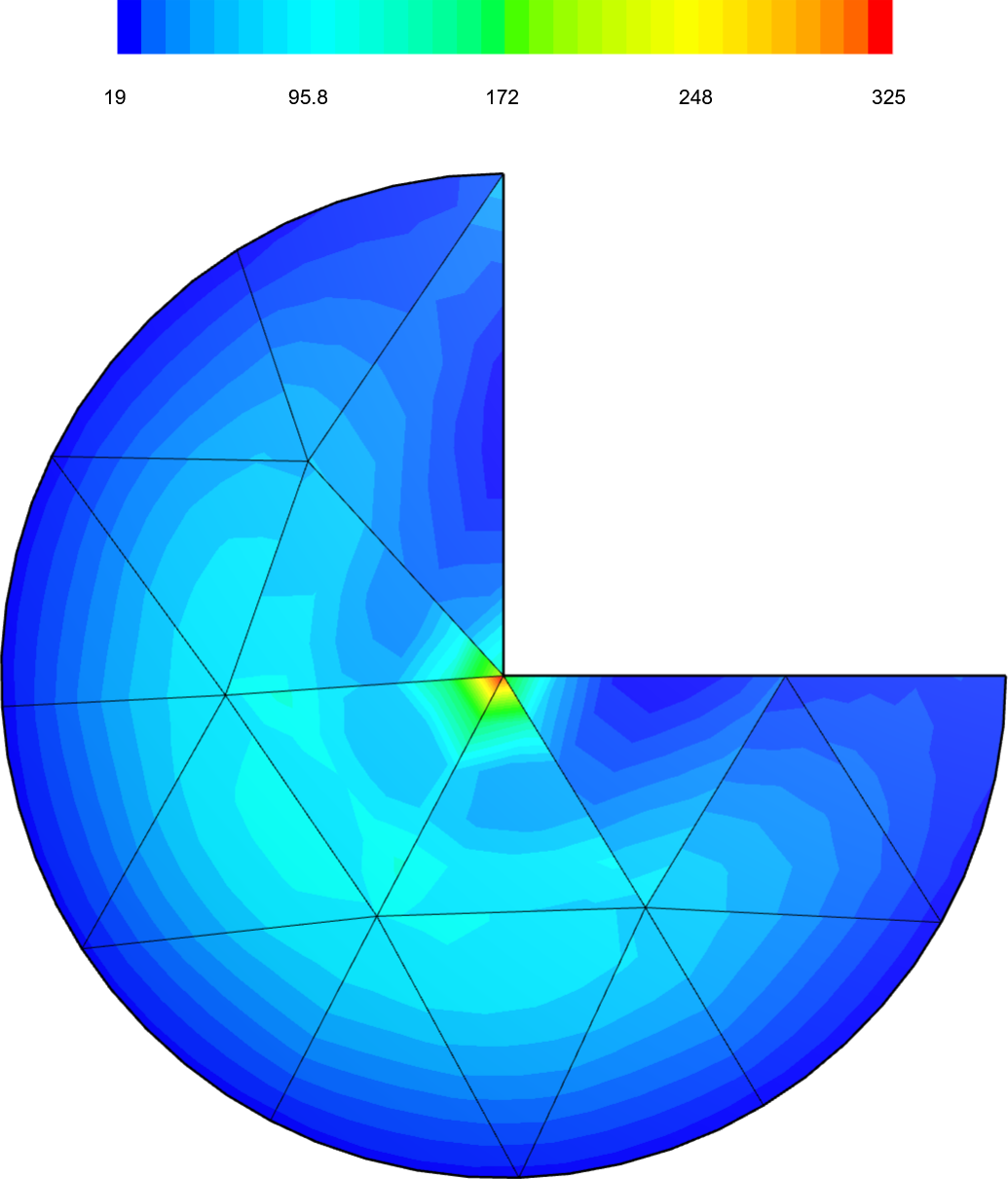}
        \includegraphics[width=0.7\linewidth]{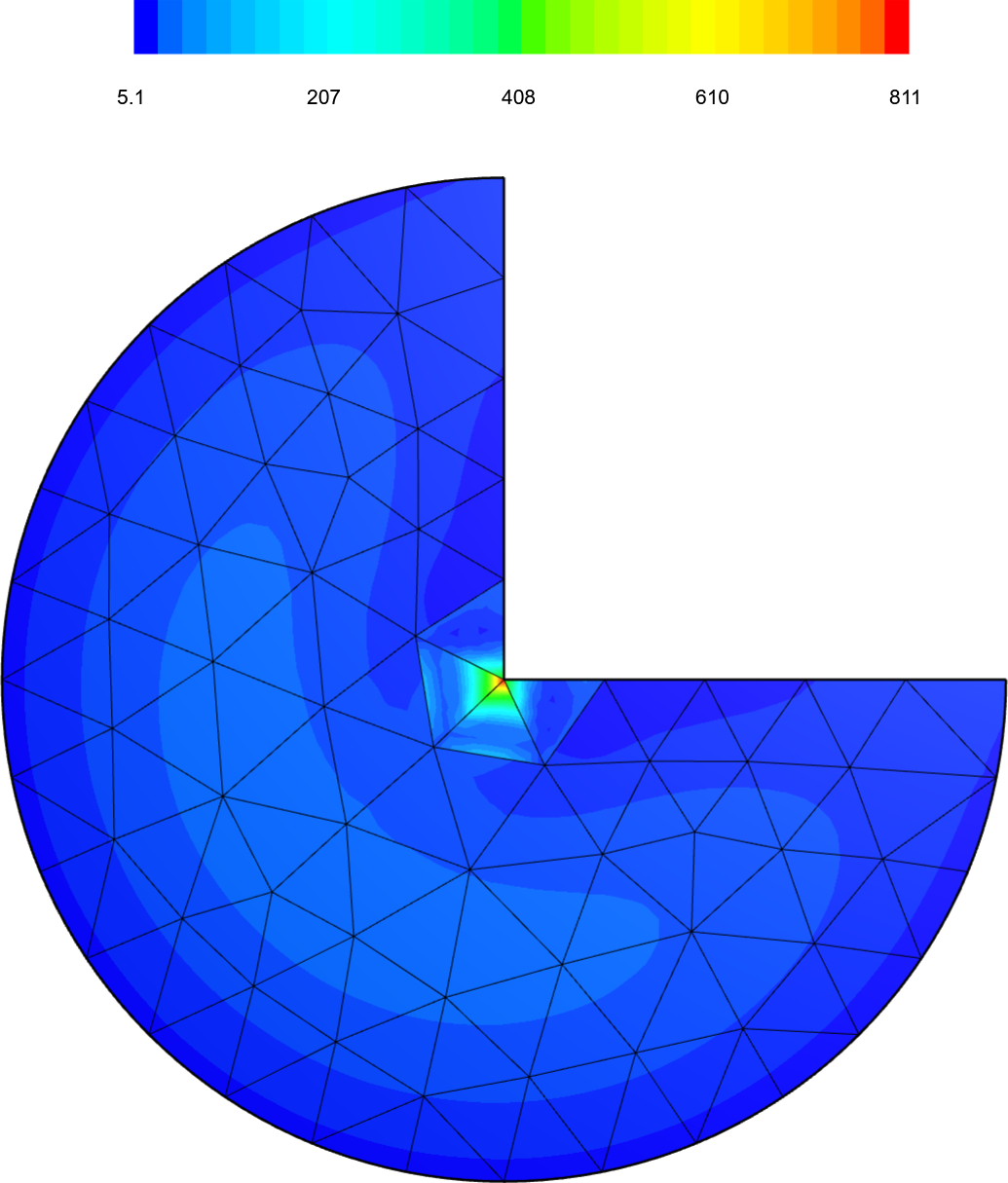}
        \includegraphics[width=0.7\linewidth]{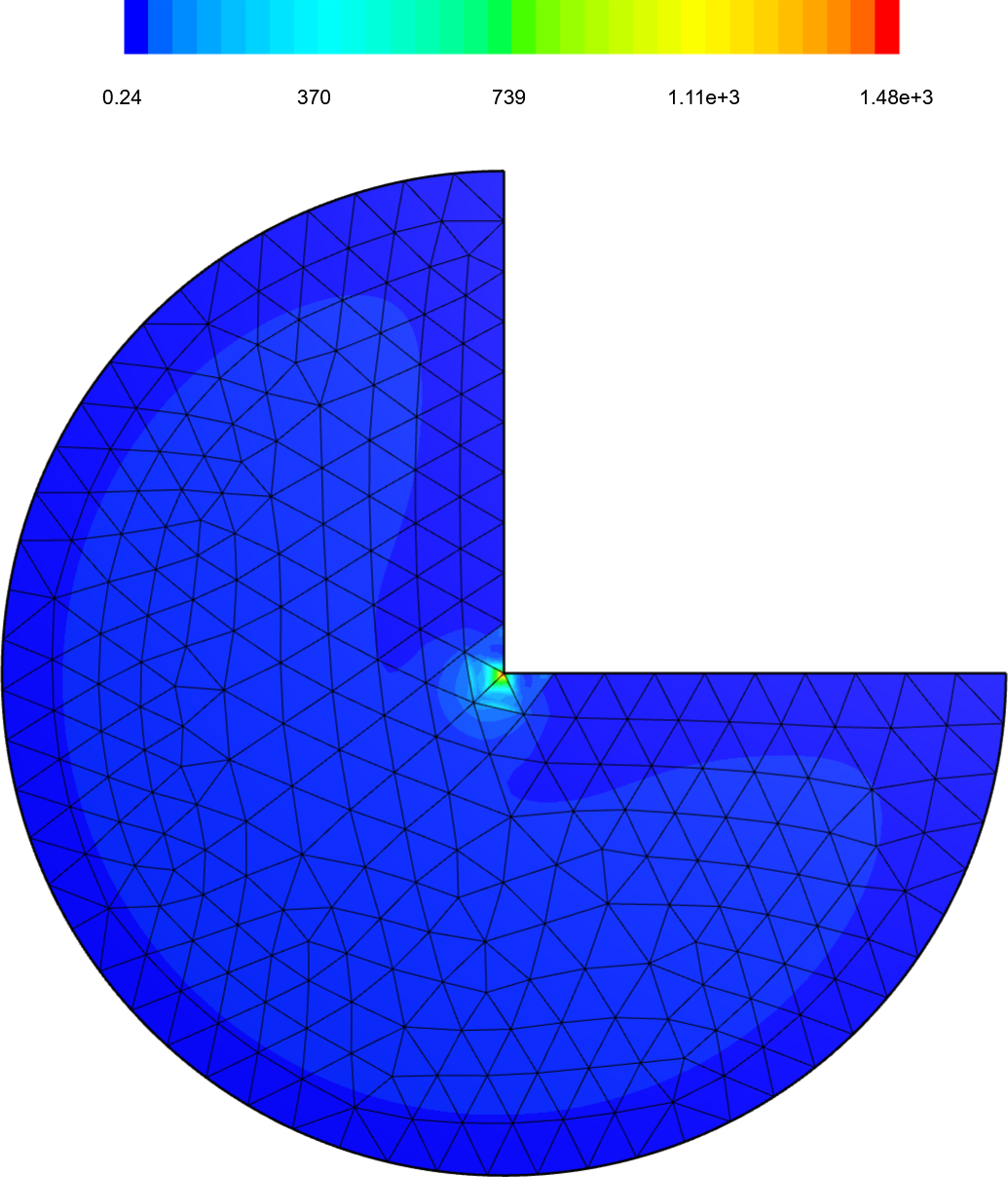}
	    	\caption{}
	    \end{subfigure}
		\caption{Solution of the norm of the bending moment $\norm{\bm{M}}$ for FFSRM (a) and TDNNS (b) over meshes with $18$, $118$ and $505$ cubic finite elements.}
		\label{fig:lshapenorm}
\end{figure} 

The results in the norm of the bending moment $\norm{\bm{M}}$ are depicted in \cref{fig:lshapenorm}. The re-entrant corner induces a singularity in the bending moments. From a purely geometric perspective, this singularity should be focused directly at the corner. However, the FFSRM formulation enforces a complete zero bending moment $\widetilde{\bm{M}} = 0$ at the re-entrant vertex. Consequently, the singularity in the bending moments is pushed inwards into the domain, yielding a very different solution than with TDNNS. In contrast, TDNNS allows for non-vanishing bending moments at the re-entrant vertex since its degrees of freedom are inherently edge-type and the tangent-tangent component of the bending moment $M_{tt}$ is not forced to vanish on the boundary.  
We note that the peak values of the norm of the bending moments greatly differ, highlighting the importance of correctly applying the boundary conditions. Finally, both formulations are correctly mapped from the reference simplex to simplices in the physical mesh, validating the discussed transformations. 

\section{Conclusions and outlook} \label{sec:conc}
This work introduced a unified method of constructing the base functions for a myriad of tensor-valued finite element spaces. Specifically, we presented the construction for the Regge, Hellan--Herrmann--Johnson, Pechstein--Sch\"oberl, Hu--Zhang, Hu--Ma--Sun, and Gopalakrishnan--Lederer--Sch\"oberl elements. 
In particular, we demonstrated that vectorial templates for tangential or normal continuity can be defined on the reference simplex, such that tensorial templates for tangential-tangential, normal-normal or tangential-normal continuity can be naturally derived from them. The latter concept was further augmented for symmetric normal-continuous tensors using the global Cartesian basis and interface-unique orthogonal vector triples built from the shared vector of the respective interface (edge tangent or face normal).       

The definition of the base functions by their association with their respective polytopes on the reference simplex allowed us to intrinsically introduce conforming transformations to physical simplices in the finite element mesh. For the Hu--Zhang and Hu--Ma--Sun elements these deviate from the classical double Piola transformations and represent an additional novelty of this work.  
Finally, the polytopal approach allows one to choose the underlying scalar polynomial $\Hone(\body)$-conforming subspace independently, such that the construction is generalised and can be used in conjunction with any standard $\C^0(\body)$-continuous Lagrange-type subspace $\CG^p(\body)$, e.g., Lagrange, Bernstein, Legendre or others.  

We demonstrated the usefulness of the method with two examples of the Reissner--Mindlin plate problem of linear elasticity. The first example showed the importance of minimal regularity of the finite element space in the case of jumping material coefficients in order to correctly capture field values at material interfaces.
Notably, the maximal value of the bending moment could not be satisfactorily approached without the minimal regularity of the approximation space.  
The second example emphasised the importance of correctly enforcing boundary conditions, and indicated their relation to the regularity of the finite element space. In particular, the effect of a singularity in the bending moments in conjunction with erroneous boundary conditions was explored. 

By its polytopal nature, the method can be extended to non-simplex-type domains. Latter is a subject for future works.

\section*{Acknowledgements}
Michael Neunteufel acknowledges support by the National Science Foundation (USA) grant DMS-2409900.

\bibliographystyle{spmpsci}   

\footnotesize{
\bibliography{Ref}   
}

\normalsize
\appendix

\section{Piola transformations} \label{ap:a}
Consistent transformations are employed to map the base functions from the reference to the physical element \cite{Mon03}. Effectively, every element in the physical domain is constructed from one reference element. If the mapping of the physical space is affine and achieved via barycentric functions for each element, the polynomial degree is maintained across transformations. 
	
	Scalar base functions transform according to
	\begin{align}
		&n(\vb{x}) = n \circ [\vb{x}^{-1}(\bm{\xi})] \, ,  &&
		\nabla_x n = \bm{J}^{-T} \nabla_\xi n \, ,
	\end{align}
	where the result concerning the Jacobi matrix is a direct consequence of the chain rule.
	
	N\'ed\'elec elements are defined via their action on the tangent vectors of the element. Consequently, a consistent transformation is given by the equality 
	\begin{align}
		\langle \bm{\theta}, \, \vb{t} \rangle \dd \curv = \langle \bm{\theta}, \, \dd \vb{s} \rangle = \langle \bm{\theta},  \, \bm{J} \dd \bm{\mu} \rangle = 
		\langle \bm{\vartheta}, \, \dd \bm{\mu} \rangle = \langle \bm{\vartheta}, \, \bm{\tau} \rangle \dd \mu \quad \iff \quad \bm{\theta} = \bm{J}^{-T} \bm{\vartheta} \, ,
	\end{align}
	known as the covariant Piola transformation.
	This is the same transformation as for gradients, thus respecting the commuting property of the de Rham complex.
	Further, vectors undergoing the latter transformation exhibit the following transformation of the curl operator 
	\begin{align}
		\mathrm{curl}_x \bm{\theta} = \nabla_x \times \bm{\theta} = (\bm{J}^{-T} \nabla_\xi) \times (\bm{J}^{-T} \bm{\vartheta}) = \cof(\bm{J}^{-T}) (\nabla_\xi \times \bm{\vartheta}) = \dfrac{1}{\det\bm{J}} \bm{J} \mathrm{curl}_\xi \bm{\vartheta} \, ,
	\end{align}
	being the so called contravariant Piola transformation.
	The result is won by observing that
	\begin{align}
		\nabla_x \times \bm{J}^{-T} = \nabla_x \times \nabla_x \bm{\xi} = 0 \, .
	\end{align}
    For two-dimensional domains the formula reduces to
    \begin{align}
    	\rot_x \bm{\theta} = \mathrm{div}_x (\bm{R} \bm{\theta}) = \dfrac{1}{\det \bm{J}} \mathrm{div}_\xi (\bm{R} \bm{\vartheta}) = \rot_\xi \bm{\vartheta} \, ,
    \end{align}
    since the curl operator produces a scalar.
	The contravariant Piola transformation is compatible with the commuting diagram and preserves normal projections on the element's boundary. To see this characteristic define the base function $\bm{\phi}$ in the reference domain and $\bm{\varphi}$ in the physical domain and equate their normal projections on the outer surface of both domains
	\begin{align}
		\langle \bm{\varphi}, \, \vb{n} \rangle \dd \surf = \langle \bm{\varphi} , \, \dd \vb{A} \rangle = \langle \bm{\varphi} , \, \cof(\bm{J}) \dd \bm{\Gamma} \rangle = \langle \bm{\phi}, \, \dd \bm{\Gamma} \rangle = \langle \bm{\phi} , \, \bm{\nu} \rangle \dd \Gamma \quad \iff \quad \bm{\varphi} = \dfrac{1}{\det \bm{J}} \bm{J} \bm{\phi} \, .  
	\end{align}
The divergence of functions mapped by a contravariant Piola transformation is given by
\begin{align}
	\int_\vol q \, \mathrm{div}_x\bm{\varphi} \, \dd  \vol &= 
	\oint_{\partial \vol} q \, \langle \bm{\varphi} , \, \vb{n} \rangle \, \dd \surf - \int_\vol \langle \nabla_x q , \, \bm{\varphi} \rangle \, \dd \vol \notag \\
	&= \oint_{\partial \Omega} \hat{q} \, \langle \dfrac{1}{\det \bm{J}} \bm{J} \, \bm{\phi} , \, \det(\bm{J}) \, \bm{J}^{-T} \bm{\nu} \rangle \, \dd \Gamma - \int_\Omega \langle \bm{J}^{-T} \nabla_\xi \hat{q} , \, \dfrac{1}{\det \bm{J}} \bm{J} \, \bm{\phi} \rangle \, \det \bm{J} \, \dd \Omega \notag \\
	&= \oint_{\partial \Omega} \hat{q} \, \langle \bm{\phi} , \, \bm{\nu} \rangle \, \dd \Gamma - \int_\Omega \langle \nabla_\xi \hat{q} , \, \bm{\phi} \rangle \, \dd \Omega \notag \\
	&= \int_\Omega \hat{q} \, \mathrm{div}_\xi \bm{\phi} \, \dd \Omega = \int_\vol q \, \mathrm{div}_\xi(\bm{\phi}) \, \dfrac{1}{\det \bm{J}} \dd \vol \qquad \forall\, q \, \in \mathit{C}^\infty(\overline{\vol}) \, ,
\end{align}
where $\hat{q} = q \circ \vb{x}$. Consequently, there holds
\begin{align}
	\mathrm{div}_x\bm{\varphi}= 
	\dfrac{1}{\det \bm{J}}\, \mathrm{div}_\xi \bm{\phi} \, .
\end{align}

\section{Consistent orientations} \label{ap:orient}
\begin{figure}
    \definecolor{asl}{rgb}{0.4980392156862745,0.,1.}
     \definecolor{asb}{rgb}{0.,0.4,0.6}
	\centering
	\begin{subfigure}{0.45\linewidth}
		\centering
		\input{figs/orient}
	\end{subfigure}
    \begin{subfigure}{0.45\linewidth}
    	\centering
    	\textcolor{asb}{\begin{align}
    			&T_1 = \{v_1, \,v_2, \,v_3, \,v_4\} \notag \\[2ex]
    			&T_2 = \{v_1, \,v_2, \,v_3, \,v_5\} \notag 
    	\end{align}}
        \vspace{0.7cm}
    \end{subfigure}
	\caption{Consistent orientations of neighbouring elements using vertex sequences. Edges are oriented from the lower to the higher vertex and faces according to the left-hand rule starting from the lowest vertex across the middle to the highest. As a result, the tangential and normal vectors at interfaces match.~\ref{fn:ccby}}
	\label{fig:consist}
\end{figure}
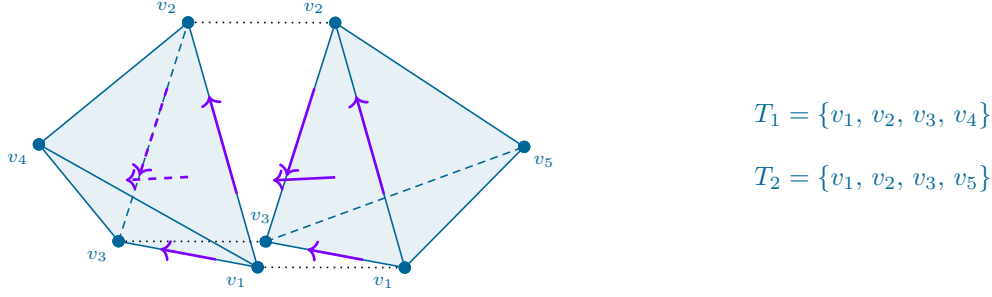

The co- and contravariant Piola transformations are not enough to ensure a consistent orientation of tangential or normal projections on edges and faces of neighbouring elements. The transformations control the size of the projections, but not whether these are parallel or anti-parallel with respect to the projection of neighbouring elements at the interface. However, consistent projections is a key requirement in ensuring no jumps occur in the trace of the respective space. In this work we rely on a solution based on the sequencing of vertices and the separation of orientational data, directly in the construction.  
We define the following rule for the orientation of edges
\begin{align}
	e_{ij} = \{v_i, \, v_j \} \qquad \text{s.t.} \quad i < j \, .
\end{align}
This means each edge starts at the lower vertex index and ends at the higher vertex index. 
This definition determines the orientation of the edge tangent vector, see \cref{fig:consist}~\ref{fn:ccby}. 
Analogously, for faces we define
\begin{align}
	f_{ijk} = \{v_i, \, v_j , \, v_k\}  \qquad \text{s.t.} \quad i < j < k  \, ,
\end{align}
such that each face is given by a sequence of increasing vertex indices. The orientation of the surface normal is given according to the left-hand rule. In other words, the direction of the normal is determined by the cross product of the vectors arising from the edges $\{v_i , \, v_j\}$ and $\{v_i, \, v_k\}$
\begin{align}
	\vb{n}_{ijk} \parallel \vb{t}_{ij} \times \vb{t}_{ik} \, . 
\end{align}
Consequently, in order to map each triangle (in 2D) or tetrahedron (in 3D) in the mesh to this orientation, we define each element as an increasing vertex-index sequence
\begin{align}
	T = \{ v_i , \, v_j , \, v_k , \, v_l \} \qquad \text{s.t.} \quad i < j < k < l \, ,
\end{align}
as depicted in \cref{fig:consist}~\ref{fn:ccby}.
The latter ensures the consistent projections of the base functions, since they are all mapped from the same reference element. However, integration in the reference element is determined by the determinant of the Jacobi matrix
\begin{align}
	\int_{\vol_e} \dd \vol = \int_{\Omega} \det \bm{J} \, \dd \Omega \, ,
\end{align}
which may be negative due to a reflection of the element in the mapping from the reference to the physical domain. The error is circumvented by taking only the absolute value of the determinant
\begin{align}
	\int_{\vol_e} \dd \vol = \int_{\Omega} |\det \bm{J}| \, \dd \Omega \, .
\end{align}
Thus, consistency is guaranteed by mapping from a single reference element. 
\begin{remark}
	The absolute value of $\det \bm{J}$ is only used for the integration over the element. In all other use-cases, the information of the sign is necessary.
\end{remark}
   
\section{Example elements} \label{ap:c}
In order to effectuate the construction methodology we give examples to the base functions of two element spaces in this section. For simplicity, we make use of the barycentric coordinates for the underlying $\CG^p(\body)$-space. On the reference triangle these read
\begin{align}
    &\lambda_0 = 1 -\xi - \eta \, ,
    && \lambda_1 = \eta \, , 
    && \lambda_2 = \xi \, .
\end{align}
On the reference tetrahedron they are
\begin{align}
    &\lambda_0 = 1 -\xi - \eta - \zeta \, ,
    && \lambda_1 = \zeta \, , 
    && \lambda_2 = \eta \, ,
    && \lambda_3 = \xi \, .
\end{align}

\subsection{The quadratic Hellan--Herrmann--Johnson element} \label{sec:quadhhj}
We recall the definitions $\bm{\iota}^\nu_1 = \vb{e}_1 - \vb{e}_2$, $\bm{\iota}^\nu_2 = -(1/2) (\vb{e}_1 + \vb{e}_2)$ and the polytopal template sets for normal continuity \cref{eq:tem2Dn}. With the definitions of the new tensorial sets as per \cref{eq:temhhj} we can define the base functions of the $\HHJ^2(\Gamma)$ element.
\begin{definition}[Quadratic Hellan--Herrmann--Johnson base functions]
    The functions are defined per polytope.
    \begin{itemize}
        \item On each edge $e_{ij}$ we have the base functions
        \begin{subequations}
            \begin{align}
            &e_{01}: & \bm{\Upsilon}_0 &= \lambda_0 \vb{e}_1 \otimes \vb{e}_1 \, , &
            \bm{\Upsilon}_1 &= \lambda_1 \bm{\iota}^\nu_1 \otimes \bm{\iota}^\nu_1 \, , &
            \bm{\Upsilon}_{01} &= \lambda_0 \lambda_{1} \vb{e}_1 \otimes \vb{e}_1 \, , \\
            &e_{02}: & \bm{\Upsilon}_0 &= \lambda_0 \vb{e}_2 \otimes \vb{e}_2 \, , &
            \bm{\Upsilon}_2 &= \lambda_2 \bm{\iota}^\nu_1 \otimes \bm{\iota}^\nu_1 \, , &
            \bm{\Upsilon}_{02} &= \lambda_0 \lambda_{2} \vb{e}_2 \otimes \vb{e}_2 \, , \\
            &e_{12}: & \bm{\Upsilon}_1 &= \lambda_1 \vb{e}_2 \otimes \vb{e}_2 \, , &
            \bm{\Upsilon}_2 &= \lambda_1 \vb{e}_1 \otimes \vb{e}_1 \, , &
            \bm{\Upsilon}_{12} &= \lambda_1 \lambda_{2} \bm{\iota}^\nu_2 \otimes \bm{\iota}^\nu_2 \, ,
        \end{align}
        \end{subequations}
        where on each edge the first two are the vertex-edge functions.
        \item The cell base functions read
        \begin{subequations}
            \begin{align}
            \bm{\Upsilon}_0 &= -\lambda_0 \sym(\vb{e}_1 \otimes \vb{e}_2) \, , & \bm{\Upsilon}_1 &= -\lambda_1 \sym(\bm{\iota}^\nu_1 \otimes \vb{e}_2) \, , &
            \bm{\Upsilon}_2 &= -\lambda_2 \sym(\bm{\iota}^\nu_1 \otimes \vb{e}_1) \, , \\
            \bm{\Upsilon}_{01}^{\tau \nu} &= \lambda_0 \lambda_1 \sym(\vb{e}_1 \otimes \vb{e}_2) \, , & 
            \bm{\Upsilon}_{01}^{\tau \tau} &= \lambda_0 \lambda_1 \vb{e}_2 \otimes \vb{e}_2 \, , & 
            \bm{\Upsilon}_{02}^{\tau \nu} &= -\lambda_0 \lambda_2 \sym(\vb{e}_1 \otimes \vb{e}_2) \, , \\ 
            \bm{\Upsilon}_{02}^{\tau \tau} &= \lambda_0 \lambda_2 \vb{e}_1 \otimes \vb{e}_1 \, , & 
            \bm{\Upsilon}_{12}^{\tau \nu} &= \lambda_1 \lambda_2 \sym(\bm{\iota}^\nu_1 \otimes \bm{\iota}^\nu_2) \, , & 
            \bm{\Upsilon}_{12}^{\tau \tau} &= \lambda_1 \lambda_2 \bm{\iota}^\nu_1 \otimes \bm{\iota}^\nu_1 \, ,
        \end{align}
        \end{subequations}
        where the first three are the vertex-cell functions and the remainder are the edge-cell base functions.
    \end{itemize}
\end{definition}

\subsection{The linear tetrahedral Regge element} \label{sec:linreg}
We recall the definition $\bm{\iota}^\tau = \vb{e}_1 + \vb{e}_2 + \vb{e}_3$ the vertex polytopal sets from \cref{eq:tem3Dtansets} and the procedure to derive the tensorial template sets in \cref{eq:toreg}, we can now define the linear tetrahedral Regge element $\Reg^1(\Omega)$.   
\begin{definition}[Linear Regge tetrahedral base functions]
    The base functions are defined on vertices with edge- or face-connectivity.
    \begin{itemize}
        \item The vertex-edge base functions read
            \begin{align}
                &e_{01}: & \bm{\Upsilon}_{0} &= \lambda_0 \vb{e}_3 \otimes \vb{e}_3 \, , & 
                 \bm{\Upsilon}_{1} &= \lambda_1 \bm{\iota}^\tau \otimes \bm{\iota}^\tau \, , \notag \\ 
                 &e_{02}: & \bm{\Upsilon}_{0} &= \lambda_0 \vb{e}_2 \otimes \vb{e}_2 \, , & 
                 \bm{\Upsilon}_{2} &= \lambda_2 \bm{\iota}^\tau \otimes \bm{\iota}^\tau \, ,  
                 \notag \\
                 &e_{03}: & \bm{\Upsilon}_{0} &= \lambda_0 \vb{e}_1 \otimes \vb{e}_1 \, , & 
                 \bm{\Upsilon}_{3} &= \lambda_3 \bm{\iota}^\tau \otimes \bm{\iota}^\tau \, ,  \\
                 &e_{12}: & \bm{\Upsilon}_{1} &= \lambda_1 \vb{e}_2 \otimes \vb{e}_2 \, , & 
                 \bm{\Upsilon}_{2} &= \lambda_2 \vb{e}_3 \otimes \vb{e}_3 \, ,  
                 \notag \\
                 &e_{13}: & \bm{\Upsilon}_{1} &= \lambda_1 \vb{e}_1 \otimes \vb{e}_1 \, , & 
                 \bm{\Upsilon}_{3} &= \lambda_3 \vb{e}_3 \otimes \vb{e}_3 \, , 
                 \notag \\
                 &e_{23}: & \bm{\Upsilon}_{2} &= \lambda_2 \vb{e}_1 \otimes \vb{e}_1 \, , & 
                 \bm{\Upsilon}_{3} &= \lambda_3 \vb{e}_3 \otimes \vb{e}_3 \, , \notag
            \end{align}
        for each edge $e_{ij}$.
        \item The vertex-face base functions are given by the remaining symmetric tensorial components 
        \begin{align}
            &f_{012}: & \bm{\Upsilon}_{0} &= \lambda_0 \sym(\vb{e}_2 \otimes \vb{e}_3) \, , & \bm{\Upsilon}_{1} &= \lambda_1 \sym(\vb{e}_2 \otimes \bm{\iota}^\tau) \, , & \bm{\Upsilon}_{2} &= -\lambda_2 \sym(\vb{e}_3 \otimes \bm{\iota}^\tau) \, ,
            \notag \\
            &f_{013}: & \bm{\Upsilon}_{0} &= \lambda_0 \sym(\vb{e}_1 \otimes \vb{e}_3) \, , & \bm{\Upsilon}_{1} &= \lambda_1 \sym(\vb{e}_1 \otimes \bm{\iota}^\tau) \, , & \bm{\Upsilon}_{3} &= -\lambda_3 \sym(\vb{e}_3 \otimes \bm{\iota}^\tau) \, ,
             \\
            &f_{023}: &  \bm{\Upsilon}_{0} &= \lambda_0 \sym(\vb{e}_1 \otimes \vb{e}_2) \, , & \bm{\Upsilon}_{2} &= \lambda_2 \sym(\vb{e}_1 \otimes \bm{\iota}^\tau) \, , & \bm{\Upsilon}_{3} &= -\lambda_3 \sym(\vb{e}_2 \otimes \bm{\iota}^\tau)  \, , 
            \notag \\
            &f_{123}: & \bm{\Upsilon}_{1} &= \lambda_1 \sym(\vb{e}_1 \otimes \vb{e}_2) \, , & \bm{\Upsilon}_{2} &= -\lambda_2 \sym(\vb{e}_1 \otimes \vb{e}_3) \, , & \bm{\Upsilon}_{3} &= \lambda_3 \sym(\vb{e}_2 \otimes \vb{e}_3) \, , \notag
        \end{align}
        for each face $f_{ijk}$.
    \end{itemize}
\end{definition}

\end{document}

%% file: figs/permut.tex
\definecolor{asl}{rgb}{0.4980392156862745,0.,1.}
\definecolor{asb}{rgb}{0.,0.4,0.6}
\begin{tikzpicture}[line cap=round,line join=round,>=triangle 45,x=1.0cm,y=1.0cm]
	\clip(-0.5,-0.5) rectangle (15.5,3.5);
	
	\draw (0,0) node[circle,fill=asb,inner sep=1.5pt] {};
	\draw (0,2) node[circle,fill=asb,inner sep=1.5pt] {};
	\draw (2,0) node[circle,fill=asb,inner sep=1.5pt] {};
	\draw [color=asb,line width=.6pt] (0,0) -- (0,2) -- (2,0) -- (0,0);
	\fill[opacity=0.1, asb] (0,0) -- (0,2) -- (2,0) -- cycle;
	\draw (0.6,0.6) node[color=asb] {$\Gamma$};
	\draw (0,0) node[color=asb,anchor=north east] {$_{v_0}$};
	\draw (2,0) node[color=asb,anchor=north west] {$_{v_2}$};
	\draw (0,2) node[color=asb,anchor=south east] {$_{v_1}$};
	\draw [-to,color=black,line width=1.pt, dashed] (0,0) -- (3,0);
	\draw [-to,color=black,line width=1.pt, dashed] (0,0) -- (0,3);
	\draw (3,0) node[color=black,anchor=west] {$\xi$};
	\draw (0,3) node[color=black,anchor=south] {$\eta$};
	
	\draw [-to,color=asl,line width=1pt] (-0.2,0.1) -- (-0.2,0.85);
	\draw (-0.2,0.85) node[color=asl,anchor=south] {$_{\tv}$};
	
	\draw (6,0) node[circle,fill=asb,inner sep=1.5pt] {};
	\draw (6,2) node[circle,fill=asb,inner sep=1.5pt] {};
	\draw (8,0) node[circle,fill=asb,inner sep=1.5pt] {};
	\draw [color=asb,line width=.6pt] (6,0) -- (6,2) -- (8,0) -- (6,0);
	\fill[opacity=0.1, asb] (6,0) -- (6,2) -- (8,0) -- cycle;
	\draw (6.6,0.6) node[color=asb] {$\Gamma$};
	\draw (6,0) node[color=asb,anchor=north east] {$_{v_0}$};
	\draw (8,0) node[color=asb,anchor=north west] {$_{v_1}$};
	\draw (6,2) node[color=asb,anchor=south east] {$_{v_2}$};
	\draw [-to,color=black,line width=1.pt, dashed] (6,0) -- (9,0);
	\draw [-to,color=black,line width=1.pt, dashed] (6,0) -- (6,3);
	\draw (9,0) node[color=black,anchor=west] {$\xi$};
	\draw (6,3) node[color=black,anchor=south] {$\eta$};
	
	\draw [-to,color=asl,line width=1pt] (6.1, -0.2) -- (6.85,-0.2);
	\draw (6.85,-0.2) node[color=asl,anchor=west] {$_{\tv}$};
	
	\draw (12,0) node[circle,fill=asb,inner sep=1.5pt] {};
	\draw (12,2) node[circle,fill=asb,inner sep=1.5pt] {};
	\draw (14,0) node[circle,fill=asb,inner sep=1.5pt] {};
	\draw [color=asb,line width=.6pt] (12,0) -- (12,2) -- (14,0) -- (12,0);
	\fill[opacity=0.1, asb] (12,0) -- (12,2) -- (14,0) -- cycle;
	\draw (12.6,0.6) node[color=asb] {$\Gamma$};
	\draw (12,0) node[color=asb,anchor=north east] {$_{v_2}$};
	\draw (14,0) node[color=asb,anchor=north west] {$_{v_1}$};
	\draw (12,2) node[color=asb,anchor=south east] {$_{v_0}$};
	\draw [-to,color=black,line width=1.pt, dashed] (12,0) -- (15,0);
	\draw [-to,color=black,line width=1.pt, dashed] (12,0) -- (12,3);
	\draw (15,0) node[color=black,anchor=west] {$\xi$};
	\draw (12,3) node[color=black,anchor=south] {$\eta$};
	
	\draw (2,1.5) node[color=black,anchor=south west] {$_{\{v_0,v_1,v_2\} \mapsto \{v_0,v_2,v_1\}}$};
	\draw [-Triangle,color=black,line width=1.pt] (2,1.5) -- (5.5,1.5);
	
	\draw [-to,color=asl,line width=1pt] (12.2,2) -- (13.05,2);
	\draw (13.05,2) node[color=asl,anchor=west] {$_{\tv}$};
	
	\draw (8,1.5) node[color=black,anchor=south west] {$_{\{v_0,v_1,v_2\} \mapsto \{v_2,v_0,v_1\}}$};
	\draw [-Triangle,color=black,line width=1.pt] (8,1.5) -- (11.5,1.5);
\end{tikzpicture}

%% file: figs/tri_nii.tex
\definecolor{asl}{rgb}{0.4980392156862745,0.,1.}
		\definecolor{asb}{rgb}{0.,0.4,0.6}
		\begin{tikzpicture}[line cap=round,line join=round,>=triangle 45,x=1.0cm,y=1.0cm]
			\clip(-2,-1.5) rectangle (12.5,4.5);
			\draw (-0.5,-0.5) node[circle,fill=asb,inner sep=1.5pt] {};
			\draw (-0.5,4) node[circle,fill=asb,inner sep=1.5pt] {};
			\draw (4,-0.5) node[circle,fill=asb,inner sep=1.5pt] {};
			\draw [color=asb,line width=.6pt] (-0.5,0) -- (-0.5,3);
			\draw [color=asb,line width=.6pt] (0,-0.5) -- (3,-0.5);
			\draw [color=asb,line width=.6pt] (0.3,3.3) -- (3.3,0.3);
			\draw [dotted,color=asb,line width=.6pt] (0,0) -- (0,3) -- (3,0) -- (0,0);
			\fill[opacity=0.1, asb] (0,0) -- (0,3) -- (3,0) -- cycle;
			\draw (-0.5,-0.5) node[color=asb,anchor=north east] {$_{v_0}$};
			\draw (4,-0.5) node[color=asb,anchor=north west] {$_{v_2}$};
			\draw (-0.5,4) node[color=asb,anchor=south east] {$_{v_1}$};
			
			\draw (-0.58,1.5) node[color=asb,anchor=west] {$_{e_{01}}$};
			\draw (1.5,-.5) node[color=asb,anchor=south] {$_{e_{02}}$};
			\draw (1.94,1.94) node[color=asb,anchor=north east] {$_{e_{12}}$};
			
			\draw [-to,color=asl,line width=1pt] (-0.75,0) -- (-0.75,0.75);
			\draw [-to,color=asl,line width=1pt] (-1.5,2.25) -- (-0.75,3);
			
			\draw [-to,color=asl,line width=1pt] (0,-0.75) -- (0.75,-0.75);
			\draw [-to,color=asl,line width=1pt] (2.25,-1.5) -- (3,-0.75);
			
			\draw [-to,color=asl,line width=1pt] (0.45,3.45) -- (1.2,3.45);
			\draw [to-,color=asl,line width=1pt] (3.45,0.45) -- (3.45,1.2);
			
			\draw [to-,color=asl,line width=1pt, dashdotted] (-1.5,1.5) -- (-0.75,1.5);
			\draw [-to,color=asl,line width=1pt, densely dashed] (-0.75,1.125) -- (-0.75,1.875);
			
			\draw [-to,color=asl,line width=1pt, dashdotted] (1.5,-1.5) -- (1.5,-0.75);
			\draw [-to,color=asl,line width=1pt, densely dashed] (1.125,-0.75) -- (1.875,-0.75);
			
			\draw [-to,color=asl,line width=1pt, dashdotted] (2,2) -- (2.75,2.75);
			\draw [-to,color=asl,line width=1pt, densely dashed] (1.7,2.3) -- (2.3,1.7);
			
			\draw [-to,color=asl,line width=1pt, dotted] (0.65,0.65) -- (1.35,0.65);
			\draw [-to,color=asl,line width=1pt, dotted] (0.65,0.65) -- (0.65,1.35);
			
			\draw (0.65,0.7)
			node[color=asb,anchor=south west] {$_{c_{012}}$};
			
			\draw [-to,color=asl,line width=1pt] (6,3) -- (7,3);
			\draw (7,3)
			node[color=asl,anchor=west] {Vertex-edge template vectors};
			\draw [-to,color=asl,line width=1pt, densely dashed] (6,2.25) -- (7,2.25);
			\draw (7,2.25)
			node[color=asl,anchor=west] {Edge template vectors};
			\draw [-to,color=asl,line width=1pt, dashdotted] (6,1.5) -- (7,1.5);
			\draw (7,1.5)
			node[color=asl,anchor=west] {Edge-cell template vectors};
			\draw [-to,color=asl,line width=1pt, dotted] (6,0.75) -- (7,0.75);
			\draw (7,0.75)
			node[color=asl,anchor=west] {Cell template vectors};
		\end{tikzpicture}

%% file: figs/tri_bdm.tex
\definecolor{asl}{rgb}{0.4980392156862745,0.,1.}
		\definecolor{asb}{rgb}{0.,0.4,0.6}
		\begin{tikzpicture}[line cap=round,line join=round,>=triangle 45,x=1.0cm,y=1.0cm]
			\clip(-2,-1.5) rectangle (12.5,4.5);
			\draw (-0.5,-0.5) node[circle,fill=asb,inner sep=1.5pt] {};
			\draw (-0.5,4) node[circle,fill=asb,inner sep=1.5pt] {};
			\draw (4,-0.5) node[circle,fill=asb,inner sep=1.5pt] {};
			\draw [color=asb,line width=.6pt] (-0.5,0) -- (-0.5,3);
			\draw [color=asb,line width=.6pt] (0,-0.5) -- (3,-0.5);
			\draw [color=asb,line width=.6pt] (0.3,3.3) -- (3.3,0.3);
			\draw [dotted,color=asb,line width=.6pt] (0,0) -- (0,3) -- (3,0) -- (0,0);
			\fill[opacity=0.1, asb] (0,0) -- (0,3) -- (3,0) -- cycle;
			\draw (-0.5,-0.5) node[color=asb,anchor=north east] {$_{v_0}$};
			\draw (4,-0.5) node[color=asb,anchor=north west] {$_{v_2}$};
			\draw (-0.5,4) node[color=asb,anchor=south east] {$_{v_1}$};
			
			\draw (-0.58,1.5) node[color=asb,anchor=west] {$_{e_{01}}$};
			\draw (1.5,-.5) node[color=asb,anchor=south] {$_{e_{02}}$};
			\draw (1.94,1.94) node[color=asb,anchor=north east] {$_{e_{12}}$};
			
			\draw [-to,color=asl,line width=1pt] (-1.5,0) -- (-0.75,0);
			\draw [-to,color=asl,line width=1pt] (-1.5,3.75) -- (-0.75,3);
			
			\draw [-to,color=asl,line width=1pt] (0,-0.75) -- (0,-1.5);
			\draw [to-,color=asl,line width=1pt] (3.75,-1.5) -- (3,-0.75);
			
			\draw [to-,color=asl,line width=1pt] (0.45,3.45) -- (0.45,4.2);
			\draw [to-,color=asl,line width=1pt] (3.45,0.45) -- (4.2,0.45);
			
			\draw [-to,color=asl,line width=1pt, densely dashed] (-1.5,1.5) -- (-0.75,1.5);
			\draw [-to,color=asl,line width=1pt, dashdotted] (-0.75,1.125) -- (-0.75,1.875);
			
			\draw [to-,color=asl,line width=1pt, densely dashed] (1.5,-1.5) -- (1.5,-0.75);
			\draw [-to,color=asl,line width=1pt, dashdotted] (1.125,-0.75) -- (1.875,-0.75);
			
			\draw [to-,color=asl,line width=1pt, densely dashed] (2,2) -- (2.75,2.75);
			\draw [-to,color=asl,line width=1pt, dashdotted] (1.7,2.3) -- (2.3,1.7);
			
			\draw [-to,color=asl,line width=1pt, dotted] (0.65,0.65) -- (1.35,0.65);
			\draw [-to,color=asl,line width=1pt, dotted] (0.65,0.65) -- (0.65,1.35);
			
			\draw (0.65,0.7)
			node[color=asb,anchor=south west] {$_{c_{012}}$};
			
			\draw [-to,color=asl,line width=1pt] (6,3) -- (7,3);
			\draw (7,3)
			node[color=asl,anchor=west] {Vertex-edge template vectors};
			\draw [-to,color=asl,line width=1pt, densely dashed] (6,2.25) -- (7,2.25);
			\draw (7,2.25)
			node[color=asl,anchor=west] {Edge template vectors};
			\draw [-to,color=asl,line width=1pt, dashdotted] (6,1.5) -- (7,1.5);
			\draw (7,1.5)
			node[color=asl,anchor=west] {Edge-cell template vectors};
			\draw [-to,color=asl,line width=1pt, dotted] (6,0.75) -- (7,0.75);
			\draw (7,0.75)
			node[color=asl,anchor=west] {Cell template vectors};
		\end{tikzpicture}

%% file: figs/tet_nii.tex
\definecolor{asl}{rgb}{0.4980392156862745,0.,1.}
		\definecolor{asb}{rgb}{0.,0.4,0.6}
		\begin{tikzpicture}
			\begin{axis}
				[
				width=30cm,height=25cm,
				view={50}{15},
				enlargelimits=true,
				xmin=-1,xmax=2,
				ymin=-1,ymax=2,
				zmin=-1,zmax=2,
				domain=-10:10,
				axis equal,
				hide axis
				]
				\draw (-0.2, -0.2, -0.2) node[circle,fill=asb,inner sep=1.5pt] {};
				\draw (-0.2, -0.2, 1.2) node[circle,fill=asb,inner sep=1.5pt] {};
				\draw (1.2, -0.2, -0.2) node[circle,fill=asb,inner sep=1.5pt] {};
				\draw (-0.2, 1.2, -0.2) node[circle,fill=asb,inner sep=1.5pt] {};
				
				\addplot3[color=asl][line width=1pt,mark=o]
				coordinates {(0.1, 0.1, 0.1)};
				
				\addplot3[color=asb][line width=0.6pt,dotted]
				coordinates {(0,0,0)(0.5,0,0)(0,0.5,0)(0,0,0)};
				\addplot3[color=asb][line width=0.6pt,dotted]
				coordinates {(0,0,0)(0,0,0.5)};
				\addplot3[color=asb][line width=0.6pt,dotted]coordinates {(0.5,0,0)(0,0,0.5)};
				\addplot3[color=asb][line width=0.6pt,dotted]coordinates {(0,0.5,0)(0,0,0.5)};
				\fill[opacity=0.1, asb] (axis cs: 0,0,0) -- (axis cs: 0.5,0,0) -- (axis cs: 0,0.5,0) -- (axis cs: 0,0,0.5) -- cycle;
				
				\draw[color=asb] (-0.2, -0.2, -0.2) node[anchor=north east] {$_{v_{0}}$};
				\draw[color=asb] (-0.2, -0.2, 1.2) node[anchor=south east] {$_{v_{1}}$};
				\draw[color=asb] (-0.2, 1.2, -0.2) node[anchor=south west] {$_{v_{2}}$};
				\draw[color=asb] (1.2, -0.2, -0.2) node[anchor=north west] {$_{v_{3}}$};
				
				\draw[line width=.6pt, color=asb](-0.2, -0.2, 0)--(-0.2,-0.2,1);
				\draw[line width=.6pt, color=asb](0, -0.2, -0.2)--(1,-0.2,-0.2);
				\draw[line width=.6pt, color=asb](-0.2, 0, -0.2)--(-0.2,1,-0.2);
				\draw[line width=.6pt, color=asb](-0, 1.0, -0.2)--(1,0,-0.2);
				\draw[line width=.6pt, color=asb](0,-0.2,1)--(1,-0.2,0);
				\draw[line width=.6pt, color=asb](-0.2,0,1)--(-0.2,1,0);
				
				\draw[-to, line width=1.pt, color=asl](-0.2, -0.2, 0)--(-0.2,-0.2,0.2);
				\draw[-to, line width=1.pt, color=asl](0, -0.2, -0.2)--(0.2,-0.2,-0.2);
				\draw[-to, line width=1.pt, color=asl](-0.2, 0, -0.2)--(-0.2,0.2,-0.2);
				
				\draw[to-, line width=1.pt, color=asl](-0.2, -0.2, 1)--(-0.3,-0.3,0.9);
				\draw[to-, line width=1.pt, color=asl](1,-0.2,-0.2)--(0.9,-0.3,-0.3);
				\draw[to-, line width=1.pt, color=asl](-0.2,1,-0.2)--(-0.3,0.9,-0.3);
				
				\draw[-to, line width=1.pt, color=asl](0,-0.2,1)--(0.2,-0.2,1);
				\draw[-to, line width=1.pt, color=asl](-0.2,0,1)--(-0.2,0.2,1);
				\draw[-to, line width=1.pt, color=asl](1,-0.2,0.2)--(1,-0.2,0);
				\draw[-to, line width=1.pt, color=asl](-0.2,1,0.2)--(-0.2,1,0);
				
				\draw[-to, line width=1.pt, color=asl](1,0.2,-0.2)--(1,0,-0.2);
				\draw[-to, line width=1.pt, color=asl](0, 1.0, -0.2)--(0.2, 1.0, -0.2);
				
				\draw[-to, line width=1.pt, color=asl,densely dashed](-0.2,-0.2,0.4)--(-0.2,-0.2,0.6);
				\draw[-to, line width=1.pt, color=asl,densely dashed](0.4,-0.2,-0.2)--(0.6,-0.2,-0.2);
				\draw[-to, line width=1.pt, color=asl, densely dashed](0.5,-0.2,0.5)--(0.7,-0.2,0.5);
				\draw[-to, line width=1.pt, color=asl, densely dashed](-0.2,0.5,0.5)--(-0.2,0.7,0.5);
				\draw[-to, line width=1.pt, color=asl, densely dashed](0.5,0.5,-0.2)--(0.7,0.5,-0.2);
				
				\draw[-to, line width=1.pt, color=asl, dashdotted](-0.2,-0.2,0.5)--(-0.4,-0.2,0.5);
				\draw[-to, line width=1.pt, color=asl, dashdotted](-0.2,-0.2,0.5)--(-0.2,-0.4,0.5);
				\draw[to-, line width=1.pt, color=asl, dashdotted](0.5,-0.2,-0.2)--(0.5,-0.4,-0.2);
				\draw[to-, line width=1.pt, color=asl, dashdotted](0.5,-0.2,-0.2)--(0.5,-0.2,-0.4);
				\draw[to-, line width=1.pt, color=asl, dashdotted](-0.2,0.5,-0.2)--(-0.2,0.5,-0.4);
				
				\draw[-to, line width=1.pt, color=asl, dashdotted](0.5,-0.2,0.5)--(0.5,0,0.5);
				\draw[-to, line width=1.pt, color=asl, dashdotted](0.5,-0.2,0.5)--(0.6,-0.1,0.6);
				
				\draw[-to, line width=1.pt, color=asl, dashdotted](-0.2,0.5,0.5)--(-0.1,0.6,0.6);
				\draw[to-, line width=1.pt, color=asl, dashdotted](-0.2,0.5,0.5)--(0,0.5,0.5);
				
				\draw[-to, line width=1.pt, color=asl, dashdotted](0.5,0.5,-0.2)--(0.6,0.6,-0.1);
				\draw[to-, line width=1.pt, color=asl, dashdotted](0.5,0.5,-0.2)--(0.5,0.5,0);
				
				\draw[-to, line width=1.pt, color=asl, dotted](0.1,0,0.1)--(0.1,0,0.3);
				\draw[-to, line width=1.pt, color=asl, dotted](0.1,0,0.1)--(0.3,0,0.1);
				
				\draw[to-, line width=1.pt, color=asl, dashdotdotted](0.1,0,0.1)--(0.1,-0.2,0.1);
				
				\draw[-to, line width=1.pt, color=asl, dotted](0.2,0.2,0.1)--(0.2,0.4,0.1);
				\draw[-to, line width=1.pt, color=asl, dotted](0.2,0.2,0.1)--(0.4,0.2,0.1);
				
				\draw[-to, line width=1.pt, color=asl, dashdotdotted](0.2,0.2,0.1)--(0.3,0.3,0.2);
				
				\draw[-to,line width=1pt, color=asl](1, 1, 1)--(1.2, 1.2, 1.007);
				\draw[color=asl] (1.2, 1.2, 1.007) node[anchor=west] {Vertex-edge template vectors};
				\draw[-to,line width=1pt, color=asl, densely dashed](1, 1, 1-0.15)--(1.2, 1.2, 1.007-0.15);
				\draw[color=asl] (1.2, 1.2, 1.007-0.15) node[anchor=west] {Edge template vectors};
				\draw[-to,line width=1pt, color=asl, dashdotted](1, 1, 1-0.3)--(1.2, 1.2, 1.007-0.3);
				\draw[color=asl] (1.2, 1.2, 1.007-0.3) node[anchor=west] {Edge-face template vectors};
				\draw[-to,line width=1pt, color=asl, dotted](1, 1, 1-0.45)--(1.2, 1.2, 1.007-0.45);
				\draw[color=asl] (1.2, 1.2, 1.007-0.45) node[anchor=west] {Face template vectors};
				\draw[-to,line width=1pt, color=asl,dashdotdotted](1, 1, 1-0.6)--(1.2, 1.2, 1.007-0.6);
				\draw[color=asl] (1.2, 1.2, 1.007-0.6) node[anchor=west] {Face-cell template vectors};
				
				\addplot3[color=asl][line width=1pt,mark=o]
				coordinates {(1.1, 1.1, 1-0.75)};
				\draw[color=asl] (1.2, 1.2, 1.007-0.75) node[anchor=west] {Cell-Cartesian template vectors};
				
				\draw[color=asb] (-0.2, -0.2, 0.5) node[anchor=west] {$_{e_{01}}$};
				\draw[color=asb] (0.5, -0.2, -0.2) node[anchor=south] {$_{e_{03}}$};
				\draw[color=asb] (-0.2,0.5,0.5) node[anchor=east] {$_{e_{12}}$};
				\draw[color=asb] (0.5,0.5,-0.2) node[anchor=east] {$_{e_{23}}$};
				\draw[color=asb] (0.5,-0.2,0.6) node[anchor=south] {$_{e_{13}}$};
				\draw[color=asb] (0.1,0.1,-0.05) node[anchor=south] {$_{f_{013}}$};
			\end{axis}
		\end{tikzpicture}

%% file: figs/tet_bdm.tex
\definecolor{asl}{rgb}{0.4980392156862745,0.,1.}
		\definecolor{asb}{rgb}{0.,0.4,0.6}
		\begin{tikzpicture}
			\begin{axis}
				[
				width=30cm,height=25cm,
				view={50}{15},
				enlargelimits=true,
				xmin=-1,xmax=2,
				ymin=-1,ymax=2,
				zmin=-1,zmax=2,
				domain=-10:10,
				axis equal,
				hide axis
				]
				\draw (-0.2, -0.2, -0.2) node[circle,fill=asb,inner sep=1.5pt] {};
				\draw (-0.2, -0.2, 1.2) node[circle,fill=asb,inner sep=1.5pt] {};
				\draw (1.2, -0.2, -0.2) node[circle,fill=asb,inner sep=1.5pt] {};
				\draw (-0.2, 1.2, -0.2) node[circle,fill=asb,inner sep=1.5pt] {};
				
				\addplot3[color=asl][line width=1pt,mark=o]
				coordinates {(0.1, 0.1, 0.1)};
				
				\addplot3[color=asb][line width=0.6pt,dotted]
				coordinates {(0,0,0)(0.5,0,0)(0,0.5,0)(0,0,0)};
				\addplot3[color=asb][line width=0.6pt,dotted]
				coordinates {(0,0,0)(0,0,0.5)};
				\addplot3[color=asb][line width=0.6pt,dotted]coordinates {(0.5,0,0)(0,0,0.5)};
				\addplot3[color=asb][line width=0.6pt,dotted]coordinates {(0,0.5,0)(0,0,0.5)};
				\fill[opacity=0.1, asb] (axis cs: 0,0,0) -- (axis cs: 0.5,0,0) -- (axis cs: 0,0.5,0) -- (axis cs: 0,0,0.5) -- cycle;
				
				\draw[color=asb] (-0.2, -0.2, -0.2) node[anchor=north east] {$_{v_{0}}$};
				\draw[color=asb] (-0.2, -0.2, 1.2) node[anchor=east] {$_{v_{1}}$};
				\draw[color=asb] (-0.2, 1.2, -0.2) node[anchor=west] {$_{v_{2}}$};
				\draw[color=asb] (1.2, -0.2, -0.2) node[anchor=west] {$_{v_{3}}$};
				
				\draw[line width=.6pt, color=asb](-0.2, -0.2, 0)--(-0.2,-0.2,1);
				\draw[line width=.6pt, color=asb](0, -0.2, -0.2)--(1,-0.2,-0.2);
				\draw[line width=.6pt, color=asb](-0.2, 0, -0.2)--(-0.2,1,-0.2);
				\draw[line width=.6pt, color=asb](-0, 1.0, -0.2)--(1,0,-0.2);
				\draw[line width=.6pt, color=asb](0,-0.2,1)--(1,-0.2,0);
				\draw[line width=.6pt, color=asb](-0.2,0,1)--(-0.2,1,0);
				
				\draw[-to, line width=1.pt, color=asl](-0.2, -0.2, -0.2)--(-0.4,-0.2,-0.2);
				\draw[-to, line width=1.pt, color=asl](-0.2, -0.2, -0.2)--(-0.2,0,-0.2);
				\draw[-to, line width=1.pt, color=asl](-0.2, -0.2, -0.2)--(-0.2,-0.2,-0.4);
				
				\draw[-to, line width=1.pt, color=asl](-0.2, -0.2, 1.2)--(-0.4,-0.2,1.4);
				\draw[-to, line width=1.pt, color=asl](-0.2,-0.2,1.2)--(-0.2,0,1);
				\draw[-to, line width=1.pt, color=asl](-0.2,-0.2,1.2)--(-0.2,-0.2,1);
				
				\draw[-to, line width=1.pt, color=asl](-0.2, 1.2, -0.2)--(-0.4, 1.4, -0.2);
				\draw[-to, line width=1.pt, color=asl](-0.2, 1.2, -0.2)--(-0.2, 1.4, -0.4);
				\draw[-to, line width=1.pt, color=asl](-0.2, 1.2, -0.2)--(-0.2, 1, -0.2);
				
				\draw[-to, line width=1.pt, color=asl](1.2, -0.2, -0.2)--(1., -0., -0.2);
				\draw[-to, line width=1.pt, color=asl](1.2, -0.2, -0.2)--(1.4, -0.2, -0.4);
				\draw[-to, line width=1.pt, color=asl](1.2, -0.2, -0.2)--(1., -0.2, -0.2);
				
				\draw[-to, line width=1.pt, color=asl,dashdotted](-0.2,-0.2,0.4)--(-0.2,-0.2,0.6);
				\draw[-to, line width=1.pt, color=asl,dashdotted](0.4,-0.2,-0.2)--(0.6,-0.2,-0.2);
				\draw[-to, line width=1.pt, color=asl, dashdotted](0.4,-0.2,0.6)--(0.6,-0.2,0.4);
				\draw[-to, line width=1.pt, color=asl, dashdotted](-0.2,0.4,0.6)--(-0.2,0.6,0.4);
				\draw[-to, line width=1.pt, color=asl, dashdotted](0.4,0.6,-0.2)--(0.6,0.4,-0.2);
				
				\draw[-to, line width=1.pt, color=asl, densely dashed](-0.2,-0.2,0.5)--(-0.4,-0.2,0.5);
				\draw[to-, line width=1.pt, color=asl, densely dashed](-0.2,-0.2,0.5)--(-0.2,-0.4,0.5);
				\draw[to-, line width=1.pt, color=asl, densely dashed](0.5,-0.2,-0.2)--(0.5,-0.4,-0.2);
				\draw[-to, line width=1.pt, color=asl, densely dashed](0.5,-0.2,-0.2)--(0.5,-0.2,-0.4);
				\draw[-to, line width=1.pt, color=asl, densely dashed](-0.2,0.5,-0.2)--(-0.2,0.5,-0.4);
				
				\draw[to-, line width=1.pt, color=asl, densely dashed](0.5,-0.2,0.5)--(0.5,-0.4,0.7);
				\draw[to-, line width=1.pt, color=asl, densely dashed](0.5,-0.2,0.5)--(0.5,-0.2,0.7);
				
				\draw[to-, line width=1.pt, color=asl, densely dashed](-0.2,0.5,0.5)--(-0.2,0.5,0.7);
				\draw[-to, line width=1.pt, color=asl, densely dashed](-0.2,0.5,0.5)--(-0.4,0.5,0.7);
				
				\draw[-to, line width=1.pt, color=asl, densely dashed](0.5,0.5,-0.2)--(0.5,0.7,-0.4);
				\draw[-to, line width=1.pt, color=asl, densely dashed](0.5,0.7,-0.2)--(0.5,0.5,-0.2);
				
				\draw[-to, line width=1.pt, color=asl, dashdotdotted](0.1,0,0.1)--(0.1,0,0.3);
				\draw[-to, line width=1.pt, color=asl, dashdotdotted](0.1,0,0.1)--(0.3,0,0.1);
				
				\draw[to-, line width=1.pt, color=asl, dotted](0.1,0,0.1)--(0.1,-0.2,0.1);

				\draw[-to,line width=1pt, color=asl](1, 1, 1)--(1.2, 1.2, 1.007);
				\draw[color=asl] (1.2, 1.2, 1.007) node[anchor=west] {Vertex-face template vectors};
				\draw[-to,line width=1pt, color=asl, densely dashed](1, 1, 1-0.15)--(1.2, 1.2, 1.007-0.15);
				\draw[color=asl] (1.2, 1.2, 1.007-0.15) node[anchor=west] {Edge-face template vectors};
				\draw[-to,line width=1pt, color=asl, dashdotted](1, 1, 1-0.3)--(1.2, 1.2, 1.007-0.3);
				\draw[color=asl] (1.2, 1.2, 1.007-0.3) node[anchor=west] {Edge-cell template vectors};
				\draw[-to,line width=1pt, color=asl, dotted](1, 1, 1-0.45)--(1.2, 1.2, 1.007-0.45);
				\draw[color=asl] (1.2, 1.2, 1.007-0.45) node[anchor=west] {Face template vectors};
				\draw[-to,line width=1pt, color=asl,dashdotdotted](1, 1, 1-0.6)--(1.2, 1.2, 1.007-0.6);
				\draw[color=asl] (1.2, 1.2, 1.007-0.6) node[anchor=west] {Face-cell template vectors};
				
				\addplot3[color=asl][line width=1pt,mark=o]
				coordinates {(1.1, 1.1, 1-0.75)};
				\draw[color=asl] (1.2, 1.2, 1.007-0.75) node[anchor=west] {Cell-Cartesian template vectors};
				
				\draw[color=asb] (-0.2, -0.2, 0.5) node[anchor=west] {$_{e_{01}}$};
				\draw[color=asb] (0.5, -0.2, -0.2) node[anchor=south] {$_{e_{03}}$};
				\draw[color=asb] (-0.2,0.5,0.5) node[anchor=east] {$_{e_{12}}$};
				\draw[color=asb] (0.5,0.5,-0.2) node[anchor=east] {$_{e_{23}}$};
				\draw[color=asb] (0.5,-0.2,0.6) node[anchor=south] {$_{e_{13}}$};
				\draw[color=asb] (0.1,0.1,-0.05) node[anchor=south] {$_{f_{013}}$};
			\end{axis}
		\end{tikzpicture}

%% file: figs/ex1dom.tex
\definecolor{wwqqcc}{rgb}{0.4,0,0.8}
\definecolor{qqwwzz}{rgb}{0,0.4,0.6}
\begin{tikzpicture}[line cap=round,line join=round,>=triangle 45,x=1cm,y=1cm]
\clip(-5.5,-1.5) rectangle (5.5,2);
\fill[line width=1pt,color=qqwwzz,fill=qqwwzz,fill opacity=0.1] (-5,-1) -- (-5,1) -- (0,1) -- (0,-1) -- cycle;
\fill[line width=1pt,color=wwqqcc,fill=wwqqcc,fill opacity=0.1] (0,1) -- (5,1) -- (5,-1) -- (0,-1) -- cycle;
\draw [line width=0.7pt,color=qqwwzz] (-5,-1)-- (-5,1);
\draw [line width=0.7pt,color=qqwwzz] (-5,1)-- (0,1);
\draw [line width=0.7pt,dashed,color=qqwwzz] (0,1)-- (0,-1);
\draw [line width=0.7pt,color=qqwwzz] (0,-1)-- (-5,-1);
\draw [line width=0.7pt,color=wwqqcc] (0,1)-- (5,1);
\draw [line width=0.7pt,color=wwqqcc] (5,1)-- (5,-1);
\draw [line width=0.7pt,color=wwqqcc] (5,-1)-- (0,-1);
\draw [line width=0.7pt,dashed,color=wwqqcc] (0,-1)-- (0,1);
\draw [line width=0.7pt,color=qqwwzz] (-5,-1)-- (-5,1);
\draw [line width=0.7pt,color=wwqqcc] (5,-1)-- (5,1);
\draw [-to,line width=1pt] (0,0) -- (1.5,0);
\draw [-to,line width=1pt] (0,0) -- (0,1.5);
\draw (1.5,0) node[anchor=west] {$x$};
\draw (0,1.5) node[anchor=south] {$y$};
\draw (-2.5,0) node[color=qqwwzz] {$A_1$};
\draw (2.5,0) node[color=wwqqcc] {$A_2$};
\end{tikzpicture}

%% file: figs/ex11.tex
\begin{tikzpicture}[scale = 1]
    			\definecolor{asl}{rgb}{0.4980392156862745,0.,1.}
    			\definecolor{asb}{rgb}{0.,0.4,0.6}
    			\begin{semilogxaxis}[
    				/pgf/number format/1000 sep={},
    				axis lines = left,
    				xlabel={Total degrees of freedom},
    				ylabel={$M_{yy}^{\mathrm{max}}$ },
    				xmin=450, xmax=4e5,
    				ymin=59, ymax=72,
    				xtick={1e2,1e3,1e4,1e5,1e6},
    				ytick={60,62,64,66,68,70},
    				legend style={at={(0.95,0.95)},anchor= north east},
    				ymajorgrids=true,
    				grid style=dotted,
    				]
    				\addplot[color=asl, mark=diamond] coordinates {
    					(597, 71.7)
                            (1290, 60.9)
                            (2676, 62.8)
                            (10317, 62.4)
                            (63399, 62.2)
    				};
    				\addlegendentry{PRM}
    				
    				\addplot[color=purple, mark=triangle] coordinates {
                        (1027, 63.3)
                        (2296, 61.5)
                        (4858, 62.2)
                        (19099, 62.2)
                        (118789, 62.1)
    				};
    				\addlegendentry{TDNNS}
    				
    				\addplot[color=asb, mark=pentagon] coordinates {
    					(1743, 59.5)
                            (3953, 61.9)
                            (8433, 61.6)
                            (33423, 62.6)
                            (208923, 63.6)
    				};
    				\addlegendentry{FFSRM}

    				\addplot[dashed,color=black, mark=none]
    				coordinates {
    					(1, 62.2)
    					(5e+6, 62.2)
    				};
    			
    			\end{semilogxaxis}

                \draw (0.55,1.4) 
    			node[anchor=north]{$\widetilde{M}_{yy}^\mathrm{max}$};

\end{tikzpicture}

%% file: figs/ex12.tex
\begin{tikzpicture}[scale = 1]
    			\definecolor{asl}{rgb}{0.4980392156862745,0.,1.}
    			\definecolor{asb}{rgb}{0.,0.4,0.6}
    			\begin{semilogxaxis}[
    				/pgf/number format/1000 sep={},
    				axis lines = left,
    				xlabel={Connected degrees of freedom},
    				ylabel={$M_{yy}^{\mathrm{max}}$ },
    				xmin=450, xmax=4e5,
    				ymin=59, ymax=72,
    				xtick={1e2,1e3,1e4,1e5,1e6},
    				ytick={60,62,64,66,68,70},
    				legend style={at={(0.95,0.95)},anchor= north east},
    				ymajorgrids=true,
    				grid style=dotted,
    				]
    				\addplot[color=asl, mark=diamond] coordinates {
    					(597-108, 71.7)
                            (1290-252, 60.9)
                            (2676-546, 62.8)
                            (10317-2196, 62.4)
                            (63399-13848, 62.2)
    				};
    				\addlegendentry{PRM}
    				
    				\addplot[color=purple, mark=triangle] coordinates {
                        (1027 - 468, 63.3)
                        (2296 - 1092, 61.5)
                        (4858 - 2366, 62.2)
                        (19099 - 9555, 62.2)
                        (118789 - 60008, 62.1)
    				};
    				\addlegendentry{TDNNS}
    				
    				\addplot[color=asb, mark=pentagon] coordinates {
    					(1743 - 1188, 59.5)
                            (3953 - 2772, 61.9)
                            (8433 - 6006, 61.6)
                            (33423 - 24156, 62.6)
                            (208923 - 152328, 63.6)
    				};
    				\addlegendentry{FFSRM}

    				\addplot[dashed,color=black, mark=none]
    				coordinates {
    					(1, 62.2)
    					(5e+6, 62.2)
    				};
    			
    			\end{semilogxaxis}

                \draw (6.2,1.4) 
    			node[anchor=south]{$\widetilde{M}_{yy}^\mathrm{max}$};

\end{tikzpicture}

%% file: figs/orient.tex
\begin{tikzpicture}
			\begin{axis}
				[
				width=30cm,height=30cm,
				view={20}{5},
				enlargelimits=true,
				xmin=-1,xmax=2,
				ymin=-1,ymax=2,
				zmin=-1,zmax=2,
				domain=-10:10,
				axis equal,
				hide axis
				]
				\addplot3[color=asb][line width=0.6pt,mark=*]
				coordinates {(0.5,0,0)(0,0.5,0)};
				\addplot3[color=asb][line width=0.6pt, densely dashed,mark=*]
				coordinates {(0,0.5,0)(0.25,0.25,0.5)};
				\addplot3[color=asb][line width=0.6pt,mark=*]
				coordinates {(0.25,0.25,0.5)(0.5,0,0)};
				\addplot3[color=asb][line width=0.6pt,mark=*]
				coordinates {(0.5,0,0)(0,0,0.25)};
				\addplot3[color=asb][line width=0.6pt,mark=*]
				coordinates {(0.25,0.25,0.5)(0,0,0.25)};
				\addplot3[color=asb][line width=0.6pt,mark=*]
				coordinates {(0,0.5,0)(0,0,0.25)};
				\fill[opacity=0.1, asb] (axis cs: 0.5,0,0) -- (axis cs: 0,0.5,0) -- (axis cs: 0,0,0.25) -- (axis cs: 0.25,0.25,0.5) -- cycle;
				\draw[color=asb] (0.5,0,0) node[anchor=north east] {$_{v_{1}}$};
				\draw[color=asb] (0,0.5,0) node[anchor=north east] {$_{v_{3}}$};
				\draw[color=asb] (0,0,0.25) node[anchor=north east] {$_{v_{4}}$};
				\draw[color=asb] (0.25,0.25,0.5) node[anchor=south east] {$_{v_{2}}$};
				\draw[->, line width=1.pt, color=asl](0.35,0.15,0)--(0.15,0.35,0);
				\draw[->, line width=1.pt, color=asl, dashed](0.175,0.325,0.35)--(0.075,0.425,0.15);
				\draw[->, line width=1.pt, color=asl](0.425,0.075,0.15)--(0.325,0.175,0.35);
				\draw[->, line width=1.pt, color=asl, dashed](0.25,0.25,1/6)--(0.15,7/30-0.1,1/6);
				
				\addplot3[color=asb][line width=0.6pt,mark=*]
				coordinates {(0.8,0.1,0)(0.3,0.6,0)(0.55,0.35,0.5)(0.8,0.1,0)};	
				\addplot3[color=asb][line width=0.6pt,mark=*]
				coordinates {(0.8,0.1,0)(1.0,0.3,0.25)};
				\addplot3[color=asb][line width=0.6pt,densely dashed,mark=*]
				coordinates {(0.3,0.6,0)(1.0,0.3,0.25)};
				\addplot3[color=asb][line width=0.6pt,mark=*]
				coordinates {(0.55,0.35,0.5)(1.0,0.3,0.25)};
				\addplot3[color=black][line width=0.6pt,dotted]
				coordinates {(0.5,0,0)(0.8,0.1,0)};
				\addplot3[color=black][line width=0.6pt,dotted]
				coordinates {(0,0.5,0)(0.3,0.6,0)};
				\addplot3[color=black][line width=0.6pt,dotted]
				coordinates {(0.25,0.25,0.5)(0.55,0.35,0.5)};
				\fill[opacity=0.1, asb] (axis cs: 0.3,0.6,0) -- (axis cs: 0.8,0.1,0) -- (axis cs: 1.0,0.3,0.25) -- (axis cs: 0.55,0.35,0.5) -- cycle;
				\draw[color=asb] (0.29,0.6,0.02) node[anchor=south] {$_{v_{3}}$};
				\draw[color=asb] (0.8,0.1,0) node[anchor=north east] {$_{v_{1}}$};
				\draw[color=asb] (1.0,0.3,0.25) node[anchor=north west] {$_{v_{5}}$};
				\draw[color=asb] (0.55,0.35,0.5) node[anchor=south east] {$_{v_{2}}$};
				\draw[->, line width=1.pt, color=asl](0.65,0.25,0)--(0.45,0.45,0);
				\draw[->, line width=1.pt, color=asl, ](0.475,0.425,0.35)--(0.375,0.525,0.15);
				\draw[->, line width=1.pt, color=asl](0.725,0.175,0.15)--(0.625,0.275,0.35);
				\draw[->, line width=1.pt, color=asl](0.55,0.35,1/6)--(0.45,7/30,1/6);
			\end{axis}
		\end{tikzpicture}